\documentclass[12pt,sort&compress]{elsarticle}
\usepackage{bm}
\usepackage{mathrsfs}
\usepackage{amsmath}
\usepackage{amsfonts}
\usepackage{amsthm}
\usepackage{color}
\usepackage{multirow}
\usepackage{booktabs}
\usepackage{arydshln}
\usepackage[export]{adjustbox}

\usepackage{caption}
\usepackage{subcaption}
\usepackage{tikz}
\usepackage{pgfplots}
\pgfplotsset{compat=1.18}

\usepackage{float}
\usepackage{url}
\usepackage{multicol}
\usepackage[top=1in,bottom=1in,left=1in,right=1in]{geometry}
\usepackage[pdfborder={0 0 0},colorlinks,allcolors=blue]{hyperref}
\usepackage{setspace}
\usepackage{enumitem}
\usepackage{lipsum}
\usepackage{gensymb}
\usepackage{newtxtext, newtxmath}

\usepackage{nth}

\usepackage{xfrac}

\theoremstyle{definition}
\newtheorem{remark}{Remark}

\graphicspath{{figures/}}

\allowdisplaybreaks

\def\thickness{H^\text{sh}_\text{th}}

\begin{document}
\begin{frontmatter}

\title{Immersogeometric fluid--structure interaction modeling of the human mitral valve and transcatheter edge-to-edge repair}

\author[ut]{Keon Ho~Kim\fnref{fn1}\corref{cor}}
\ead{kkeonho@utexas.edu}
\author[isu]{Ashton~M.~Corpuz\fnref{fn1}}
\author[isu]{Monu~Jaiswal}
\author[ut]{Michael~S.~Sacks}
\author[isu]{Ming-Chen~Hsu\corref{cor}}
\ead{jmchsu@iastate.edu}

\cortext[cor]{Corresponding authors}
\fntext[fn1]{These two authors contributed equally to this work.}

\address[ut]{Oden Institute for Computational Engineering and Sciences, University of Texas at Austin, Austin, TX 78712, USA}
\address[isu]{Department of Mechanical Engineering, Iowa State University, Ames, IA 50011, USA}


\begin{abstract}
Transcatheter edge-to-edge repair (TEER) treats mitral regurgitation with a clip that grasps and coapts the mitral valve (MV) leaflets.
Predicting its outcome requires resolving the interaction between the blood flow and the complete MV apparatus, in which highly deformable leaflets and chordae tendineae undergo large deformation and complex contact.
In this work, we introduce an immersogeometric fluid--structure interaction (FSI) formulation in which the coupling among the valvular components, the surrounding flow, and the repair is enforced weakly, without requiring conforming discretizations.
A fully parameterized reconstruction pipeline is developed to convert segmented images into analysis-ready isogeometric meshes.
The leaflets and chordae are modeled as Kirchhoff--Love shells and cables, and the shell formulation is modified to include an in-vivo prestrain tensor that accounts for the loading already present in the imaged geometry.
We propose a subcell approach that improves the quadrature for the interfacial constraint and eliminates the over-refinement otherwise needed to prevent unphysical leakage through the leaflets.
The MV repair is modeled by a penalty formulation that enforces leaflet coaptation over a prescribed region, allowing clip number and placement to be varied without remeshing or modeling the device.
A lumped-parameter model at the outflow tract of the left heart imposes a flow-dependent pressure rather than a prescribed waveform, so that the split between forward and regurgitant flow matches physiological behavior.
We validate the reconstruction against the segmented imaging data, demonstrate mesh convergence of the proposed formulations, and apply the framework to a prolapsed MV and its single- and double-clip repairs to quantify pre- and post-operative changes in valve dynamics and left-heart hemodynamics.
\end{abstract}


\end{frontmatter}

\newpage
\tableofcontents

\newpage
\section{Introduction}\label{s:intro}
The mitral valve (MV) is one of the two atrioventricular valves in the heart and directs the unidirectional flow of blood from the left atrium to the left ventricle. 
Dysfunction of the MV can lead to mitral regurgitation (MR), in which the valve leaflets fail to coapt completely, allowing backflow of blood into the left atrium during systole~\cite{nishimura2016mitral, el2018mitral}.
Transcatheter edge-to-edge repair (TEER), in which a clip device is delivered percutaneously to grasp and coapt the anterior and posterior leaflets, has emerged as a minimally invasive treatment option particularly for elderly patients and those at high risk for open-heart surgery~\cite{feldman2005percutaneous, hausleiter2023mitral, meulengracht2025outcomes}. 
Its efficacy, however, remains variable across patients~\cite{samimi2023predictors, yokoyama2024outcomes, hatab2024echocardiographic}. 
In particular, residual MR and elevated transvalvular gradients are common following TEER and have been linked to negative clinical outcomes and reduced long-term survival~\cite{neuss2017elevated, higuchi2021impact, singh2024combined}.
Additionally, the number and placement of clips are typically determined heuristically during the procedure, and their influence on valve function and left-heart hemodynamics remains poorly understood~\cite{puls2020implantation, kassar2022anatomical}. 
Consequently, a predictive framework capable of quantifying post-TEER hemodynamic outcomes could provide valuable guidance for treatment planning and optimization~\cite{ooida2024silico}.

Computational modeling has been widely employed to investigate MV mechanics and to evaluate various repair strategies~\cite{kunzelman1998flexible, maisano2005annular, votta2007geoform, votta20023}.
For example, structure-only finite element (FE) analyses have provided detailed insights into leaflet and chordal functions~\cite{kunzelman1993finite, kunzelman1997annular, cochran1998effect, salgo2002effect, lim2005three, khalighi2017mitral, drach2018comprehensive}.
More recently, patient-specific simulations reconstructed from pre-operative three-dimensional transesophageal echocardiography (3D TEE) have been used to predict outcomes of interventions such as undersized annuloplasty rings and TEER~\cite{liu2023computational, simonian2023patient, simonian2025patient}.
However, structure-only models cannot reproduce hemodynamic quantities such as transvalvular pressure gradients, regurgitant jets, or ventricular filling patterns, which clinical studies have shown to correlate strongly with mortality following the procedure~\cite{sato2022hemodynamic, shibahashi2024mismatch, sammour2024impact}.
Computational fluid dynamics (CFD) and fluid--structure interaction (FSI) models of the MV~\cite{peskin1980modeling, kunzelman2007fluid, griffith2009simulating, lau2010mitral, ma2013image, gao2017coupled, toma2017fluid, biffi2019workflow, feng2019analysis, cai2019some} and of repair techniques~\cite{lau2011fluid, kamakoti2019numerical, caballero2020comprehensive, vellguth2022effect, dabiri2022simulation, kim2025diastolic} have therefore been proposed.
These methods are based on formulations such as arbitrary Lagrangian--Eulerian~\cite{bavo2016fluid, joda2019comparison}, space--time~\cite{Takizawa17a, terahara2020heart}, immersed boundary~\cite{griffith2012immersed, borazjani2013fluid}, fictitious domain~\cite{dehart2003three, dehart2003computational}, and smoothed particle hydrodynamics~\cite{mao2016fluid, laha2024smoothed}.
Many of these MV FSI studies rely on idealized or manually constructed valve geometries, with flow domains that are often reduced to simple tubular channels or built from patient anatomy through case-specific manual procedures that are time-consuming and difficult to reproduce.

In the present study, we propose a new computational FSI method for modeling the human MV and the TEER procedure to investigate post-repair valve function and hemodynamics.
Constructing a high-quality MV geometry from medical imaging typically involves manual, model-specific processing that is labor-intensive, operator-dependent, and difficult to reproduce~\cite{biffi2019workflow, vellguth2019user, pouch2012semi}.
To address this challenge, we develop a fully parameterized surface reconstruction pipeline that converts segmented in-vivo MV imaging data directly into an analysis-ready isogeometric mesh.
A partial differential equation (PDE)-based parameterization is proposed to avoid the geometric distortions and coordinate singularities associated with conventional coordinate-dependent mappings.
With the parametric coordinates of the segmented points fixed in advance, surface reconstruction is then formulated as a single regularized linear least-squares problem that is solved directly in one step.
This avoids the iterative distance-minimization schemes that treat the parametric locations as additional unknowns, require initialization and case-by-case parameter tuning, and may converge to suboptimal stationary points on noisy or incomplete data~\cite{pouch2012semi, khalighi2018multi, moola2026valvefit}.
The direct solve yields smooth B-spline leaflet surfaces that are combined with functionally equivalent chordae tendineae~\cite{simonian2025patient} to form complete MV geometric models.
Because these surfaces are reconstructed from in-vivo 3D TEE data, the imaged geometry is already under physiological load and is not stress-free.
Treating it as stress-free places the tissue in the compliant part of its nonlinear response and overpredicts leaflet deformation.
We account for this effect by modifying the isogeometric Kirchhoff--Love shell formulation to include an in-vivo prestrain tensor that composes with the deformation gradient and shifts the leaflet stress--strain response into its physiological range~\cite{grashow2006biaixal, amini2012vivo, liu2023computational, simonian2023patient}.
The prestrained constitutive law and its consistent linearization are derived within the plane-stress and incompressibility reduction of the shell, which allows the imaged geometry to serve directly as the reference configuration.

To capture the interaction between the valvular structures and the surrounding blood flow in the left heart, we adopt immersogeometric analysis (IMGA)~\cite{Kamensky15ch, Hsu14dh, Hsu15fb, Kamensky15gr, Kamensky17bi}, which immerses isogeometric analysis (IGA) geometries directly into a non-body-fitted fluid domain.
Although IMGA has been applied successfully to aortic valve FSI~\cite{Xu18id, Wu18gv, Yu18ga, Hsu18el, wu2019immersogeometric, terahara2020heart, johnson2020thinner, xu2021computational, johnson2022effects, neighbor2023leveraging}, simulating full human MV FSI remains challenging. 
The MV apparatus combines highly deformable leaflets with multiple tendinous cords and undergoes large deformation, coaptation, and contact interactions over the cardiac cycle. 
To accommodate this complexity, we formulate a numerical framework that couples these structural components consistently while maintaining stable enforcement of the contact conditions.
The fluid--structure coupling requires additional treatment.
In conventional IMGA, the coupling is enforced at quadrature points associated with the structural discretization.
As the background fluid mesh is refined, the fixed density of these points becomes insufficient to evaluate the interfacial integrals accurately, and the resulting loss of constraint enforcement produces unphysical flow leakage through the leaflets unless the structural mesh is over-refined.
We remove this dependence using a subcell quadrature strategy that recursively subdivides the quadrature cells within each shell element to generate a denser set of coupling points.
Since these points are independent of the structural quadrature points, the coupling is enforced effectively without adding structural degrees of freedom or increasing the quadrature order beyond what the shell discretization requires.

Modeling the repair presents further complexity.
Existing computational studies of TEER either model the clip device explicitly or impose nodal constraints~\cite{lau2011fluid, kamakoti2019numerical, caballero2020comprehensive, vellguth2022effect, dabiri2022simulation, kim2025diastolic}.
Explicit representation of the device requires additional geometric modeling and involved contact treatment, whereas nodal constraints act only at mesh nodes, tying the extent of the grasp to the underlying discretization.
To address this challenge, we introduce a patch-based penalty formulation, defined over the clipped region within IMGA, that enforces coaptation of the grasped leaflet regions without modeling the device geometry.
Because the grasped patch is specified parametrically and independently of the mesh, the number of clips and their placement can be varied without modifying the MV discretization.
This allows different repair configurations to be compared consistently within the same model, which we use to assess the effectiveness of single- and double-clip placements.

For the fluid domain, patient-specific left heart models are rarely available and do not generally match an MV reconstructed from another subject, while a simplified conduit such as a straight tube does not reproduce the ventricular filling and ejection that load the leaflets.
We therefore propose a parameterized left heart model fitted to the reconstructed MV annulus and driven by a single prescribed volumetric flowrate, which sets both the ventricular wall motion and the atrial pressure waveform.
During systole, the volume displaced by the ventricle leaves along two competing paths: forward through the outflow tract and backward through the incompetent MV, with the split governed by their relative impedance.
Prescribing a fixed outlet pressure waveform, as is common in MV FSI studies~\cite{caballero2020comprehensive, pasta2023silico}, sets the forward-path load independently of the computed flow and thus prescribes part of what a TEER model should predict.
We instead couple the fluid domain to a lumped-parameter model~\cite{Westerhof09, gao2017coupled} that returns a flow-dependent pressure at the outflow tract, with the outlet switching to a no-slip state to represent aortic valve closure.
A single calibrated parameter set is used across all cases, and the drop in peak left ventricular pressure under moderate-to-severe MR, together with its recovery after repair, is therefore a result of the coupled problem rather than an imposed condition.

\begin{figure}[!t]\centering
    \includegraphics[width=\textwidth]{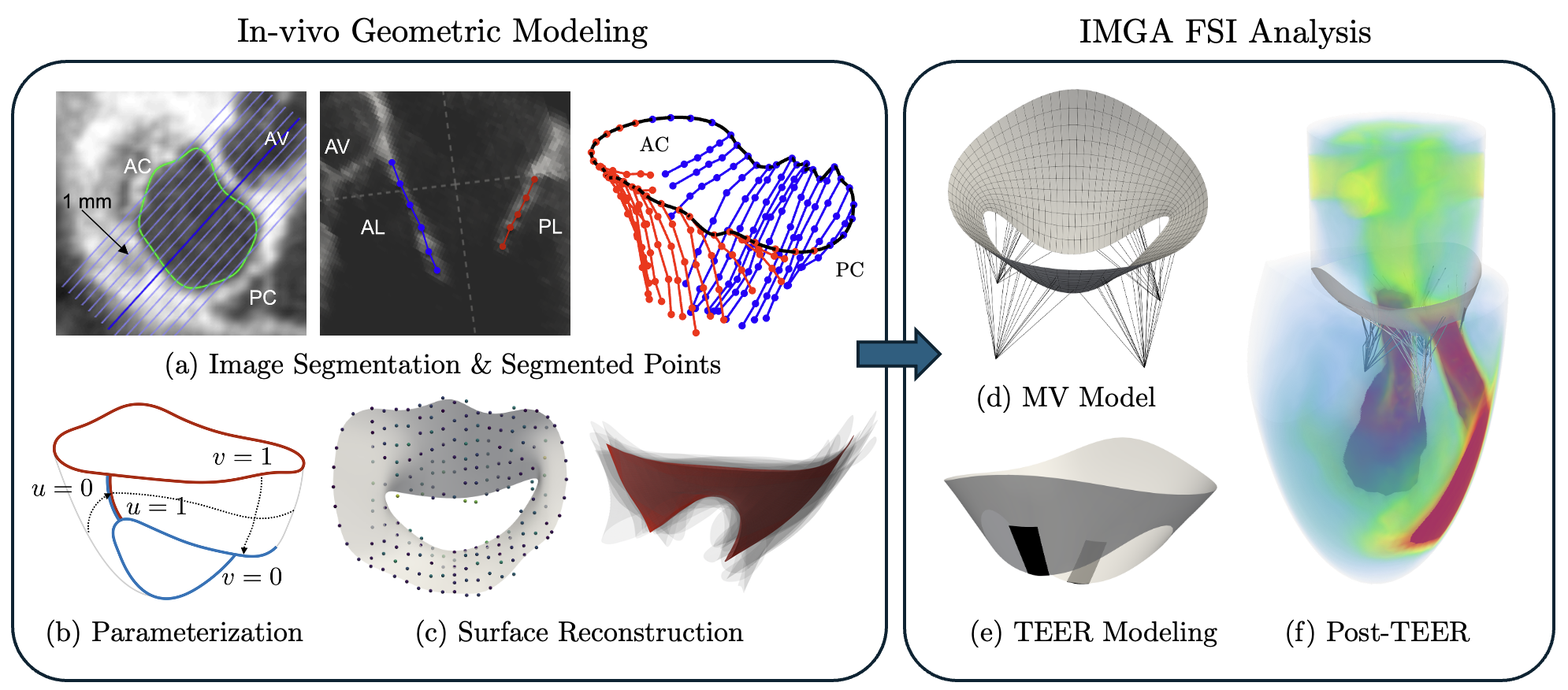}
  \caption{Overall schematic of the proposed framework for in-vivo geometric modeling and immersogeometric FSI analysis after MV TEER. 
  (a) 3D TEE images showing the parallel cross-sections used to segment the leaflets from commissure to commissure, and obtained segmented points representing MV leaflets (AL: anterior leaflet; AC: anterior commissure; PL: posterior leaflet; PC: posterior commissure; AV: aortic valve; Green curve: MV annulus; Blue dotted curve: anterior slice; Red dotted curve: posterior slice). The black curve represents the B-spline interpolated annulus.
  (b) PDE-based parameterizations in both $(u,v)$ parametric directions.
  (c) In-vivo surface reconstruction showing patient-specific and population-averaged leaflet surfaces.
  (d) Analysis-ready IGA MV model including valve leaflets and chordae tendineae.
  (e) Clip patches for modeling TEER.
  (f) Post-TEER hemodynamics at peak diastole.}
  \label{f:overview_figure}
\end{figure}

We assess the accuracy and robustness of the proposed model through convergence studies of the valvular dynamics, the TEER penalty formulation, the subcell quadrature, and the coupled FSI problem.
We then apply the framework to a series of clinically relevant configurations, including the pre-operative prolapse and post-operative single- and double-clip repairs, and quantify the resulting changes in MV dynamics and left-heart hemodynamics.
Because both the geometric and FSI models are fully parameterized, the complete pipeline from medical image reconstruction to FSI analysis can be applied to different patient anatomies and TEER configurations with minimal manual intervention.
Figure~\ref{f:overview_figure} shows a complete overview of the proposed image-to-FSI modeling framework, including leaflet segmentation from 3D TEE images, parameterization, reconstruction of in-vivo MV geometries, clip-patch placement, and FSI analysis of the repaired valve.
Overall, this work establishes a mathematically rigorous FSI framework that relaxes the geometric and discretization constraints of conventional methods for simulating patient-specific atrioventricular valves and their repairs under large deformation, multibody contact, and extended leaflet coaptation.

This paper is organized as follows. 
Section~\ref{s:MV_geo} outlines the geometric modeling pipeline, including the segmentation, parameterization, surface reconstruction, and regularization of the leaflets, as well as the modeling of the chordae tendineae.
Section~\ref{s:simFk} presents the continuous and discrete immersogeometric FSI formulations for simulating MV dynamics, including the prestrained Kirchhoff--Love shell formulation for the leaflets, the penalty coupling for TEER, and the subcell quadrature strategy for fluid--structure coupling.
Section~\ref{s:appl} validates the in-vivo surface reconstruction from segmented medical images, evaluates the developed FSI model through convergence tests, and investigates pre- and post-operative hemodynamic responses for several clinically relevant configurations.
Section~\ref{s:discussion} concludes with final remarks and future directions.

\section{Geometric modeling}\label{s:MV_geo}
The leaflets in our FSI model are thin shells discretized by IGA, which requires a smooth B-spline surface, at least $C^1$-continuous, that also serves as the basis for analysis.
In this section, we present a fully parameterized human MV geometric modeling pipeline that can directly convert segmented medical images into analysis-ready IGA meshes. 
To ensure physiological realism and temporal consistency, we consider the entire MV apparatus at the same cardiac time point, specifically at end-diastole, including valve leaflets, chordae tendineae (CT) structures, and papillary muscle locations. 

\subsection{Leaflets}
\label{s:leaflet}
Advanced in-vivo imaging techniques, such as cardiac computed tomography and magnetic resonance imaging, provide an accessible method of creating patient-specific valve reconstructions~\cite{wang2013finite}. 
Additionally, 3D TEE, which is the clinical gold standard for MV imaging, is widely used for acquiring in-vivo data with high temporal resolution and contrast~\cite{pouch2012semi}.
However, fitting a smooth, analysis-ready spline surface to the noisy and sparsely sampled point clouds produced by image segmentation remains difficult~\cite{moola2026valvefit}.
Here, we propose a reconstruction method in which the segmented leaflet data are parameterized, fitted, and regularized to produce the MV surfaces used in our analysis.

\subsubsection{Image segmentation}\label{s:segmentation}
Three-dimensional TEE images acquired as standard of care during the TEER procedure are pre-processed and segmented by interactive cross-sectional tracing, following the established protocol of~\citet{rego2018noninvasive}.
Figure~\ref{f:overview_figure}a shows the parallel cross-sections of the 3D TEE images used to segment the leaflets from commissure to commissure. 
The mitral annulus is oriented into a short-axis view, and the geometric center of the orifice is aligned with the intersection of the septolateral and intercommissural axes. 
Septolateral cross-sections covering the full intercommissural span are generated at approximately $1 \, \text{mm}$ intervals. 
In each cross-section, the anterior and posterior leaflets are delineated as distinct midsurfaces extending from the annulus to the free edge, which are illustrated as blue and red curves in Figure~\ref{f:overview_figure}a. 
The commissural regions are segmented to the extent that the anterior and posterior leaflets can be distinguished in the TEE image.
As a result, the segmented point cloud shown in Figure~\ref{f:overview_figure}a contains unsegmented regions near the commissures. 
In general, such an incomplete segmentation near the commissures requires substantial interpolation across the missing regions. 
To minimize the introduction of artificially interpolated data, we reconstruct the annular contour by using cubic B-spline interpolation, which provides a continuous boundary for subsequent surface reconstruction.
Figure~\ref{f:overview_figure}a shows the resulting segmented leaflet point cloud along with the continuous annular curve reconstructed using cubic B-spline interpolation.
Note that we do not introduce any additional points within the commissural regions except along the annular points.

\subsubsection{Parameterization}
\label{s:parameterization}
To construct an analysis-ready IGA mesh, we first generate a smooth parameterization of the segmented point cloud, which enables high-quality surface reconstruction.
Although common parameterizations such as chord-length mapping or cylindrical/spherical coordinates are straightforward, they often fail to capture complex anatomical geometries and may introduce distortions, uneven point distributions, or coordinate singularities unless the geometry naturally aligns with the chosen coordinates~\cite{ma1995parameterization, FLOATER1997231}.
Instead, we develop a PDE-based method that yields a robust and smooth parameterization without relying on a specific coordinate system.
The parametric coordinates of the segmented points are obtained from a harmonic problem on the embedded MV surface with Dirichlet boundary conditions as follows:
\begin{equation}\label{eq:parameterization_laplace}
\nabla^2 \mathbf{U} \left(\mathbf{X}\right) = \mathbf{0}\text{ ,}
\end{equation}
in which $\mathbf{U} = \left(u,v\right) \in \mathbb{R}^2$ is a parametric variable and $u$ and $v$ are the circumferential and radial parametric coordinates, respectively,
$\mathbf{X} \in \Omega_{\text{MV}} \subset \mathbb{R}^3$ is a segmented point, and $\nabla^2 = \nabla \cdot \nabla$ denotes the surface Laplacian on the embedded MV geometry.
The Laplace equation defines a harmonic mapping that minimizes the surface Dirichlet energy subject to the prescribed boundary conditions.
For $u$, we first identify one seam location in the MV commissure area to unwrap the periodic circumferential direction and prescribe Dirichlet values along the commissural boundaries. 
Similarly, Dirichlet boundary conditions are imposed along the free edge and annulus to define a normalized radial coordinate $v$.
Since the harmonic solution propagates boundary values smoothly through the interior, the resulting parametric coordinates avoid coordinate singularities associated with prescribed cylindrical or spherical coordinate systems.
Figure~\ref{f:overview_figure}b illustrates the harmonic parameterization.
Such a PDE-based mapping offers a smooth distribution of parameters over the unstructured point cloud and provides a suitable foundation for the subsequent IGA surface reconstruction, ensuring that the resulting mesh accurately captures the MV anatomical characteristics.

The segmented MV geometry is represented by a sparse point cloud without explicit surface connectivity.
To discretize this harmonic problem on the segmented MV point cloud, we approximate the Laplace problem by using a symmetrized $k$-nearest-neighbor graph Laplacian~\cite{belkin2003laplacian,hein2007graph,calder2022improved}. 
For the circumferential parameterization, the anterior and posterior leaflets are treated independently by constructing a separate graph Laplacian for each leaflet, which prevents artificial connections between spatially adjacent points across different leaflets.
For the anterior leaflet, the two commissural boundaries are assigned $u=0$ and $u=0.5$, whereas for the posterior leaflet, they are assigned $u = 0.5$ and $u = 1$. 
Note that $u = 0$ and $u = 1$ correspond to the same commissural seam in the periodic circumferential parameterization.
In contrast, the radial parameterization uses a single nearest-neighbor graph-harmonic system with $v = 0$ along the free edge and $v = 1$ along the annular boundary. 
The resulting graph-harmonic systems are solved subject to the prescribed Dirichlet boundary values to obtain the parametric coordinates of the interior points.

\subsubsection{Surface reconstruction}
\label{s:surface_fitting}
Conventional surface fitting methods minimize the distances between the data points and the fitted surface, and since the parametric location of each point on the surface is unknown, the resulting least-squares distance problem is nonlinear and must be solved iteratively. Such iterative methods are useful when enforcing nonlinear constraints. However, they require time-consuming parameter tuning and may converge to suboptimal stationary points rather than the global minimizer, particularly in the presence of noise or incomplete data~\cite{hoschek1988intrinsic, speer1998global,  wang2006fitting, merchel2023adaptive}. 
Since the PDE-based parameterization fixes the parametric coordinates of the segmented points in advance, the reconstruction is linear in the control points, and a single-shot direct method can yield a smooth and analysis-ready surface by solving a single regularized linear system.
Because this direct formulation avoids the convergence issues, sensitivity to initialization, and computational cost associated with iterative optimization, we adopt it to reconstruct the MV surfaces from the segmented imaging data.

The parametric coordinates $\left(u,v\right)$ obtained from the PDE-based parameterization in Section~\ref{s:parameterization} are used to construct a three-dimensional MV surface in the form of a tensor product of B-splines. 
This parametric surface maps the regular computational domain to the physical MV geometry, yielding a smooth description suitable for complex deformations in IGA.

Let $N_{i,p} \left(u\right)$ and $N_{j,q} \left(v\right)$ be univariate B-spline basis functions of degree $p$ and degree $q$ in the $u$- and $v$-parametric directions, respectively, defined over the knot vectors $\mathcal{U} = \{ u_0, u_1, \cdots, u_{n+p+1} \}$ and $\mathcal{V} = \{ v_0, v_1, \cdots, v_{m+q+1} \}$.
A three-dimensional surface $\mathscr{S}$ defined by $n+1$ and $m+1$ control points in the $u$- and $v$-parametric directions, respectively, can then be written as follows:
\begin{equation}
    \mathscr{S} \left(u,v\right) = \sum_{i=0}^{n} \sum_{j=0}^m \mathbf{Q}_{ij}\, N_{i,p}\left(u\right)\, N_{j,q} \left(v\right)\text{ ,}
\end{equation}
in which $\mathbf{Q}_{ij} \in \mathbb{R}^3$ are control points.
A uniform periodic knot vector is used to enforce continuity across the MV surface in the circumferential direction:
\begin{equation}
    u_i = \frac{i-p}{n-p+1}, \quad \text{for }\, i = 0,\cdots,n+p+1\text{ .}
\end{equation}
Note that the domain of the surface is $[u_p,u_{n+1}]$, and there are $p$ wrapping control points in the $u$-direction.
The radial direction employs an open uniform knot vector:
\begin{equation}
    v_i = \begin{cases}
    0, &\text{ if }\, i = 0,\cdots,q\text{ ,} \\
    \frac{i-q}{m-q+1}, \, &\text{ if } i = q+1, \cdots, m\text{ ,}\\
    1, &\text{ if }\, i = m+1,\cdots, m+q+1\text{ .}
    \end{cases}
\end{equation}
We use $p=q=3$ in our surface reconstruction process, which offers global $C^2$ continuity of the resulting B-spline surface, providing a suitable geometric representation of the human MV.

Given $N$ segmented MV points $\mathbf{P} \in \mathbb{R}^{N \times 3}$ and their corresponding parametric coordinates $\left(u_k,v_k\right)$, for $k=1,\cdots,N$, from the PDE-based parameterization, the objective is to determine the set of $M=\left(n+1\right)\times\left(m+1\right)$ control points $\mathbf{Q} \in \mathbb{R}^{M \times 3}$ that best approximates the overall geometry.
The condition that the B-spline surface fits the segmented points can be written as
\begin{equation}\label{eq:surface fit equation}
    \mathbf{M}_\text{s} \mathbf{Q} = \mathbf{P}\text{ ,}
\end{equation}
in which $\mathbf{M}_\text{s} \in \mathbb{R}^{N \times M} $  is a B-spline basis matrix with entries
\begin{equation}
    \left\{\mathbf{M}_\mathrm{s}\right\}_{\left(k,i+(n+1)j\right)} = N_{i,p}\left(u_k\right)\, N_{j,q}\left(v_k\right), \quad \text{for }\, i = 0,\cdots,n\, \text{ and }\, j = 0,\cdots,m\text{ .}
\end{equation} 
Since the parametric coordinates $\left(u_k,v_k\right)$ are fixed by the parameterization, Eq.~\eqref{eq:surface fit equation} is linear in the control points.
The system is overdetermined when the number of data points exceeds the number of control points.
In practice, $\mathbf{Q}$ is obtained by the linear least-squares solution~\cite{kim2025vivo}:
\begin{equation}\label{eq:fitting equation}
    \mathbf{Q}= \left(\mathbf{M}_\text{s}^\top\, \mathbf{M}_\text{s} \right)^{-1} \mathbf{M}_\text{s}^\top\,  \mathbf{P}\text{ .}
\end{equation}
We emphasize that the least-squares solution here refers to the direct solution of the overdetermined linear system in Eq.~\eqref{eq:surface fit equation}, not to an iterative minimization of point-to-surface distances, and it is obtained in a single solve without initialization. This produces a smooth parametric control mesh of the human MV that preserves the anatomical features while satisfying the geometric regularity requirements for IGA. Note that we use a reduced matrix system only for unique control points in the actual surface reconstruction process.

\subsubsection{Regularization}\label{s:regularization}
Regularization is often added to B-spline surface fitting of cardiovascular structures to reduce the sensitivity of the resulting surface to noisy or incomplete imaging data~\cite{kim2025vivo, moola2026valvefit, Bayer2014cmmm}.
To establish a robust anatomical shape of the human MV, we incorporate two additional features into our B-spline surface reconstruction.
First, we constrain the reconstruction of the MV leaflets by fitting the annular curve as the initial step. 
As the most reliably identified valve structure in clinical images, the annulus provides a stable bound that ensures consistent leaflet attachment in the left heart, reduces geometric ambiguity, and defines a natural coordinate framework for surface reconstruction. 
Constraining both the annulus and the free edge provides a fully bounded and stable surface reconstruction problem. 
However, the complete free edge cannot be reliably recovered from our segmentation because of missing leaflet regions near the commissures. 
We therefore use only the closed annular curve as the primary boundary constraint.
For the annular curve reconstruction, we use the obtained $u$-parameters described in Section~\ref{s:parameterization}. 
The curve is then fitted using the corresponding reduced form for a curve of the B-spline surface formulation presented in Section~\ref{s:surface_fitting}, retaining only the basis functions in the $u$-direction.
Note that we use the points on the annular curve as shown in Figure~\ref{f:overview_figure}a, not the interpolated B-spline curve itself.
After computing annular control points, we reduce the entire system in Eq.~\eqref{eq:surface fit equation} by removing the columns of the B-spline basis matrix associated with the obtained annular control points and incorporating their contributions into the right-hand side of the reduced system.

The commissural regions are relatively narrow compared with the anterior and posterior leaflets, as illustrated by the anterior and posterior commissures in Figure~\ref{f:overview_figure}a.
Therefore, a specialized regularization is required to stabilize the surface reconstruction and accurately capture the anatomical characteristics of the human MV.
We employ first-order regularization to penalize steep variations in slope, with periodic coupling in the $u$-direction to maintain the continuity across the seam. 
In addition, we use directional weights to more accurately reflect the physiological shape of the mitral valve in each parametric direction.
Suppose $\mathbf{D}_u \in \mathbb{R}^{M\times M}$ and $\mathbf{D}_v \in \mathbb{R}^{M\times M}$ are first-order finite difference (FD) operators defined over the control net in the $u$- and $v$-directions.
The first-order penalty terms are defined as
\begin{equation}\label{first order}
    \lambda_u  \mathbf{D}_u^\top  \mathbf{D}_u + \lambda_v \mathbf{D}_v^\top  \mathbf{D}_v\text{ ,}
\end{equation}
in which $\lambda_u$ and $\lambda_v$ are weights corresponding to their relative contributions.
Large penalty parameters prioritize smoothness or uniformity over the given data, resulting in a poor approximation of the original shape. 
To avoid this, we use $\lambda_{(\cdot)} \approx c \,\|\mathbf{M}_{\text{s}}\|_2^2/\|\mathbf{D}_{(\cdot)}\|_2^2$, in which $c$ is a non-dimensional scalar coefficient.
This scaling balances the spectral magnitudes of the basis and FD operators, ensuring the dimensionless coefficient $c$ imposes a consistent regularization constraint independent of point cloud density and control mesh resolution.

To further regulate curvatures, we incorporate discrete second-order FD operators, $\mathbf{D}^2_u$, $\mathbf{D}_u\mathbf{D}_v$, and $\mathbf{D}^2_v$, yielding additional penalty terms
\begin{equation}\label{second order}
    \lambda_{uu} \left(\mathbf{D}^2_u\right)^\top \mathbf{D}^2_u + \lambda_{uv} \mathbf{D}_v^\top \mathbf{D}_u^\top \mathbf{D}_u \mathbf{D}_v  + \lambda_{vv} \left(\mathbf{D}^2_v\right)^\top \mathbf{D}^2_v\text{ ,}
\end{equation}
in which $\lambda_{uu}$, $\lambda_{uv}$, and $\lambda_{vv}$ are penalty parameters determined by the same scaling strategy as in the first-order operators.
Note that we impose the periodic condition along the seam in the circumferential direction.

\subsection{Chordae tendineae}
\label{s:CT}
The mitral valve chordae tendineae (MVCT) are the tendinous cords that connect the MV leaflets to the papillary muscles within the left ventricle, and their primary function is to prevent the valve leaflets from prolapsing into the left atrium during ventricular contraction. 
Without a proper CT structure in the model, the leaflets would lack the necessary mechanical tethering to withstand the high pressure from the left ventricle, leading to an incorrect simulation of valve closure~\cite{khalighi2017mitral,drach2018comprehensive,khalighi2019development}.
Consequently, it is crucial to model functionally equivalent MVCTs for accurate and physiologically realistic MV models.

\begin{figure}[t!]\centering
    \begin{subfigure}[t]{0.65\textwidth}\centering
    \includegraphics[width=\textwidth]{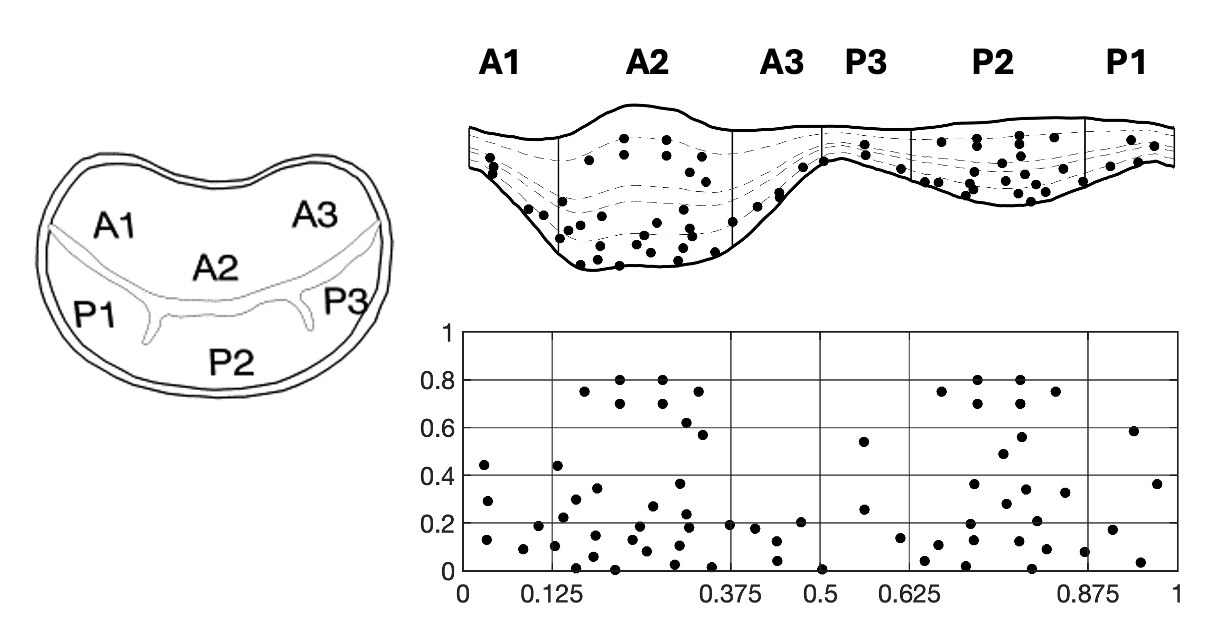}
    \caption{}\label{f:MVCT}
    \end{subfigure}
    \hspace{0.005\textwidth}
    \begin{subfigure}[t]{0.33\textwidth}\centering
    \includegraphics[width=\textwidth]{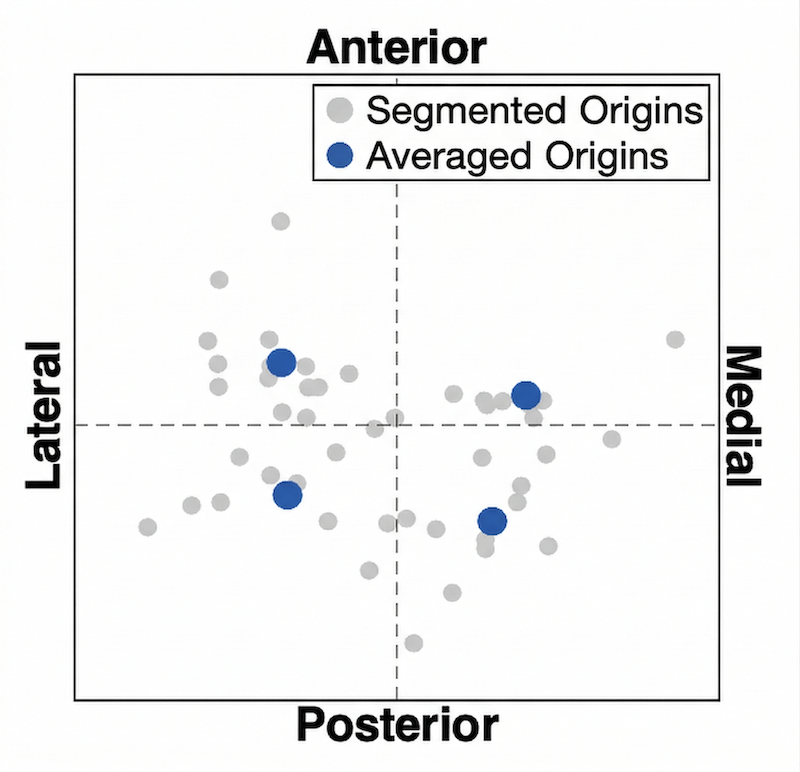}
    \caption{}\label{f:PM origins}
    \end{subfigure}

  \caption{(a) Population-averaged MVCT insertion locations on the leaflet surface and in the parametric domain. (b) Population-averaged CT origins on the normalized annular plane based on patient-specific TEE images. These normalized mappings are used to generate an anatomical distribution of MVCT insertions and CT origins for any given MV leaflet geometry.}
\end{figure}

CTs, however, are not fully visible in standard-of-care TEE images.
To model CT structures, we leverage the population-averaged CT distribution maps established in previous studies~\cite{liu2023computational, simonian2023patient, simonian2025patient}.
Since the proposed PDE-based leaflet reconstruction is inherently parameterized, we can seamlessly project normalized CT insertion locations to any reconstructed MV leaflet surface. 
The distribution of CT insertions on the leaflets is obtained from an anatomical analysis of excised human hearts~\cite{liu2023computational,simonian2025patient}.
We then normalize the distribution map to directly correspond to the parametric variables, which allows CT insertion locations to be generated directly on the MV leaflet surface.
Figure~\ref{f:MVCT} shows the MVCT insertion locations on both the physical MV surface and the parametric domain with denser chordal insertions in the coaptation zone and sparser insertions in the clear zone.

Similar to the CT insertions, CT origins are difficult to consistently visualize from TEE imaging. 
To overcome this limitation, we use a population-averaged CT origins database developed from anatomical analysis of patient-specific TEE images by~\citet{simonian2025patient}, which normalizes the position of each visible origin relative to its respective annulus and then averages these normalized positions by quadrant.
The vertical locations are determined by averaging their origins within each quadrant and measuring the distance from the annular plane.
Figure~\ref{f:PM origins} shows fully parameterized CT origins on the parameterized annular plane.
This provides stable and reproducible anchor points for the CTs, which are critical for ensuring that the physical model replicates the proper tethering and support of the mitral valve leaflets. 
The averaged CT origins are mapped to the reconstructed physical MV surface along with the CTs, allowing us to build a comprehensive and physiologically realistic MV model.
In addition, we model CTs as branchless fibers that reduce geometric complexity while maintaining the essential load-transfer role of the subvalvular apparatus~\cite{khalighi2017mitral}.

\section{Simulation framework}\label{s:simFk}
In this section, we present the IMGA framework used for modeling the human MV FSI and the MV repair. 
The FSI problem includes blood flow, the left heart chambers, and the complete MV apparatus, in which the leaflet model accounts for an in-vivo prestrain and contact is enforced through a nonlocal approach.
We then introduce a penalty coupling formulation for TEER and describe the immersogeometric discretization of the coupled system.

\subsection{Fluid--structure interaction problem}\label{subsec:FSI}
The continuous FSI problem considers $(\Omega_\text{s})_t$ and $(\Omega_\text{f})_t$ to be regions occupied by immersed hyperelastic structures and an incompressible Newtonian fluid at time $t$, respectively, which are coupled by kinematic and traction compatibility conditions at the shared fluid--structure interface $(\Gamma_\text{I})_t$.
The variational formulation of the FSI problem following the augmented Lagrangian approach~\cite{BazHsu12b} is defined as: find fluid velocity
$\mathbf{u}_\text{f}\in\mathcal{S}_u$, fluid pressure $p\in\mathcal{S}_p$, structure displacement $\mathbf{y}\in\mathcal{S}_y$, and a fluid--solid interface Lagrange multiplier $\pmb{\lambda}\in\mathcal{S}_\ell$ such that for all $\mathbf{w}_\text{f}\in\mathcal{V}_u$, $q\in\mathcal{V}_q$, $\mathbf{w}_\text{s}\in\mathcal{V}_y$, and $\delta\pmb{\lambda}\in\mathcal{V}_\ell$,
\begin{align}
    \label{general-al-variational}
    \nonumber &B_\text{f}\left(\{\mathbf{w}_\text{f},q\},\{\mathbf{u}_\text{f},p\}\right)-F_\text{f}\left(\{\mathbf{w}_\text{f},q\}\right) + B_\text{s}\left(\mathbf{w}_\text{s},\mathbf{y}\right)-F_\text{s}\left(\mathbf{w}_\text{s}\right) \\
    \nonumber & \quad + \int_{\left( \Gamma_\text{I} \right)_t}(\mathbf{w}_\text{f}-\mathbf{w}_\text{s})\cdot\pmb{\lambda} \, \text{d}\Gamma + \int_{\left( \Gamma_\text{I} \right)_t}\delta\pmb{\lambda}\cdot(\mathbf{u}_\text{f}-\mathbf{u}_\text{s}) \, \text{d}\Gamma\\
    & \quad + \int_{\left( \Gamma_\text{I} \right)_t}(\mathbf{w}_\text{f}-\mathbf{w}_\text{s})\cdot\beta \, (\mathbf{u}_\text{f}-\mathbf{u}_\text{s}) \, \text{d}\Gamma = 0\text{ ,}
\end{align}
in which $\mathcal{S}_{(\cdot)}$ and $\mathcal{V}_{(\cdot)}$ are trial solution and test function spaces,  $B_\text{f}$, $B_\text{s}$, $F_\text{f}$, and $F_\text{s}$ are variational forms defining the fluid and structure subproblems, $\mathbf{u}_\text{s}$ is the material time derivative of $\mathbf{y}$, and $\beta$ is a penalty parameter.
The fluid--structure coupling conditions along the fluid--structure interface are enforced by the integral terms over $\left( \Gamma_\text{I} \right)_t$.
The fluid and structure subproblems are defined in the subsequent sections.

\subsection{Fluid formulation}\label{s:fluid}
The fluid subproblem in Eq.~\eqref{general-al-variational} is governed by the incompressible Navier--Stokes equations posed on a moving domain and is given in the arbitrary Lagrangian--Eulerian (ALE) description~\cite{Hughes81a} as follows:
\begin{align}\label{eq:fluid-formulation}
    \nonumber B_\text{f}&\left(\left\{\mathbf{w}_\text{f},q\right\},\left\{\mathbf{u}_\text{f},p\right\}\right) - F_\text{f}(\{\mathbf{w}_\text{f},q\}) \\
    \nonumber =& \int_{(\Omega_\text{f})_t}\mathbf{w}_\text{f}\cdot\rho_\text{f}\left(\left.\frac{\partial\mathbf{u}_\text{f}}{\partial t}\right|_{\hat{\mathbf{x}}}+\left(\mathbf{u}_\text{f}-\hat{\mathbf{u}}\right)\cdot\nabla \mathbf{u}_\text{f}\right)~\text{d}\Omega + \int_{(\Omega_\text{f})_t}\pmb{\varepsilon}(\mathbf{w}_\text{f}):\pmb{\sigma}~\text{d}\Omega + \int_{(\Omega_\text{f})_t}q\nabla \cdot\mathbf{u}_\text{f}~\text{d}\Omega \\
    &- \gamma\int_{(\Gamma_\text{f}^\text{h})_t}\mathbf{w}_\text{f}\cdot\rho_\text{f}\{\left(\mathbf{u}_\text{f}-\hat{\mathbf{u}}\right)\cdot\mathbf{n}_\text{f}\}_{-}~\mathbf{u}_\text{f}~\text{d}\Gamma -\int_{(\Omega_\text{f})_t}\mathbf{w}_\text{f}\cdot\rho_\text{f}\mathbf{f}_\text{f}~\text{d}\Omega - \int_{(\Gamma_\text{f}^\text{h})_t}\mathbf{w}_\text{f}\cdot\mathbf{h}_\text{f}~\text{d}\Gamma\text{ ,}
\end{align}
in which $\rho_\text{f}$ is the fluid mass density, $\partial(\cdot)/\partial t|_{\hat{\mathbf{x}}}$ is the time derivative taken with respect to the fixed coordinates $\hat{\mathbf{x}}$ of the spatial configuration, $\hat{\mathbf{u}}$ is the (arbitrary) velocity with which the fluid subproblem domain ${(\Omega_\text{f})_t}$ deforms, $\nabla$ is the gradient taken with respect to the spatial coordinate $\mathbf{x}$ of the current configuration, $\pmb{\varepsilon}(\cdot)$ is the symmetric gradient operator given by $\pmb{\varepsilon}\left(\mathbf{w}\right) = \frac{1}{2}\left(\nabla \mathbf{w}+\nabla \mathbf{w}^\top\right)$, $\pmb{\sigma} = -p\mathbf{I} + 2 \mu_\text{f} \, \pmb{\varepsilon}\left(\mathbf{u}_\text{f}\right)$ is the fluid Cauchy stress, $\mu_\text{f}$ is the dynamic viscosity, $\mathbf{f}_\text{f}$ is a prescribed body force, $\mathbf{h}_\text{f}$ is a prescribed traction on $\Gamma_\text{f}^\text{h}$, and the $\gamma$ $(\geq 0)$ term, referred to as backflow stabilization~\cite{Moghadam11a}, improves the well-posedness of the problem when there is significant inflow through the Neumann boundary $\Gamma_\text{f}^\text{h}$. In this term, $\{\cdot\}_{-}$ isolates the negative part of its argument and $\mathbf{n}_\text{f}$ is the outward-facing normal vector to the fluid domain.

\subsection{Structural formulations}\label{s:structural}
This section describes the structural formulations for the MV leaflets and CT, the coupling between them, the treatment of leaflet prestrain and contact, and the penalty-based modeling of TEER.
We represent the valve leaflets as thin shells and the CT structures as IGA cables, based on their thickness relative to their other dimensions.
We use superscripts ``sh'' and ``ca'' to denote the shell and cable, respectively, and decompose the structural displacement and test function as $\mathbf{y} = \left\{\mathbf{y}^\text{sh}, \mathbf{y}^\text{ca}\right\}$ and $\mathbf{w}_\text{s} = \left\{\mathbf{w}^\text{sh}_\text{s}, \mathbf{w}^\text{ca}_\text{s}\right\}$, respectively.
The variational form of the structural subproblem is then written as
\begin{equation}
    \label{eq:bs-full}
    B_\text{s}\left(\mathbf{w}_\text{s},\mathbf{y}\right) = B^\text{sh}_\text{s}\left(\mathbf{w}^\text{sh}_\text{s},\mathbf{y}^\text{sh}\right) + B^\text{ca}_\text{s}\left(\mathbf{w}^\text{ca}_\text{s},\mathbf{y}^\text{ca}\right)\text{ ,}
\end{equation}
and
\begin{equation}
    \label{fs-full}
    F_\text{s}\left(\mathbf{w}_\text{s}\right) = F^\text{sh}_\text{s}\left(\mathbf{w}^\text{sh}_\text{s}\right) + F^\text{ca}_\text{s}\left(\mathbf{w}^\text{ca}_\text{s}\right)\text{ ,}
\end{equation}
in which $B^\text{sh}_\text{s}, F^\text{sh}_\text{s}$ and $B^\text{ca}_\text{s}, F^\text{ca}_\text{s}$ are the variational formulations of the shell and cable subproblems, respectively.

\subsubsection{Thin shell formulation for the leaflets}\label{s:shell}
We describe the MV leaflet surface as a hyperelastic isogeometric Kirchhoff--Love shell~\cite{Kiendl15ez}, which assumes that cross sections initially normal to the midsurface remain straight and normal after deformation. 
A material point $\mathbf{x}^{\text{sh}}$ in the shell continuum can be represented by a midsurface position $\mathbf{r}$ and a unit normal $\mathbf{a}_3$ as
\begin{equation}\label{eq:curvlinear}
\mathbf{x}^{\text{sh}}\left(\xi^1,\xi^2,\xi^3\right) = \mathbf{r}\left(\xi^1,\xi^2\right) + \xi^3\,\mathbf{a}_3\left(\xi^1,\xi^2\right)\text{ ,}
\end{equation}
in which $\left(\xi^1,\xi^2\right)$ are in-plane curvilinear coordinates, $\xi^3 \in [-\thickness/2,\thickness/2]$ is the through-thickness coordinate, and $\thickness$ denotes the shell thickness.  

By using the curvilinear coordinate description, the covariant base vectors are defined as $\mathbf{g}_{i}=\partial \mathbf{x}^{\text{sh}}/\partial \xi^i$, and the associated metric coefficients are $g_{ij}=\mathbf{g}_{i}\cdot\mathbf{g}_{j}$ for $i,j=1,2,3$.
The contravariant basis $\mathbf{g}^{i}$ is defined such that $\mathbf{g}^{i}\cdot\mathbf{g}_{j}=\delta^{i}_{j}$, and the contravariant metric follows as $\left[g^{ij}\right]=\left[g_{ij}\right]^{-1}$.
In the Kirchhoff--Love shell theory, we only consider the in-plane strain components and neglect transverse shear strains.
Under the assumption of a linear distribution of strain through the thickness, the in-plane metric tensor is expressed as 
$
g_{\alpha\beta}=a_{\alpha\beta}-2\,\xi^3 b_{\alpha\beta}
$,
in which $a_{\alpha\beta}=\mathbf{a}_{\alpha}\cdot\mathbf{a}_{\beta}$ and $b_{\alpha\beta}=\mathbf{a}_{\alpha,\beta}\cdot\mathbf{a}_3$, for $\alpha,\beta = 1,2$. 
Here, $\mathbf{a}_{\alpha}=\partial \mathbf{r}/\partial \xi^\alpha$ is the midsurface tangent vector, and $\mathbf{a}_3=(\mathbf{a}_1\times\mathbf{a}_2)/\|\mathbf{a}_1\times\mathbf{a}_2\|$ is the unit normal.
These relations apply to both reference and current configurations, and we denote the undeformed quantities by $\mathring{(\cdot)}$. 
The determinant of the deformation mapping is $J=\sqrt{|g_{ij}|/|\mathring{g}_{ij}|}$, and its in-plane counterpart is $J_o=\sqrt{|g_{\alpha\beta}|/|\mathring{g}_{\alpha\beta}|}$.

The weak form of the thin shell problem is defined as
\begin{align}\label{eq:b2-shell}
\nonumber B^\text{sh}_\text{s}(\mathbf{w}_\text{s}^{\text{sh}},\mathbf{y}^{\text{sh}}) - F^\text{sh}_\text{s}(\mathbf{w}_\text{s}^{\text{sh}})  
=& \int_{(\mathcal{S}^\text{sh})_0}\mathbf{w}_\text{s}^{\text{sh}}\cdot\rho^\text{sh}_\text{s} \,\thickness \left.\frac{\partial^2\mathbf{y}^{\text{sh}}}{\partial t^2}\right|_{\mathring{\mathbf{x}}}\,\text{d}\mathcal{S}
+ \int_{(\mathcal{S}^\text{sh})_0}\int_{-\thickness/2}^{\thickness/2}\delta \mathbf{E}:\mathbf{S}\,\text{d}\xi^3\,\text{d}\mathcal{S} \\
& - \int_{(\mathcal{S}^\text{sh})_0}\mathbf{w}_\text{s}^{\text{sh}}\cdot\rho^\text{sh}_\text{s} \,\thickness \,\mathbf{f}^\text{sh}_\text{s}\,\text{d}\mathcal{S}
- \int_{(\mathcal{S}^\text{sh})_t}\mathbf{w}_\text{s}^{\text{sh}}\cdot\mathbf{h}^{\text{net}}_\text{s}\,\text{d}\mathcal{S}\text{ ,}
\end{align}
in which $(\mathcal{S}^\text{sh})_0$ and $(\mathcal{S}^\text{sh})_t$ are the reference and deformed midsurfaces, $\mathbf{y}^{\text{sh}}$ is the midsurface displacement, $\rho^\text{sh}_\text{s}$ is the density, $\mathbf{f}_\text{s}^\text{sh}$ is a body force, and $\mathbf{h}^{\text{net}}_\text{s}$ is the net traction from both sides of the shell. 
The second Piola--Kirchhoff stress tensor $\mathbf{S}$ is obtained from a strain energy density function $\psi$: $\mathbf{S}=\partial_\mathbf{E}\psi$, with the Green--Lagrange strain $\mathbf{E}=\tfrac{1}{2}(\mathbf{C}-\mathbf{I})$ and the right Cauchy--Green tensor $\mathbf{C}$.
The variation $\delta\mathbf{E}$ is taken with respect to $\mathbf{y}^{\text{sh}}$ in the direction of $\mathbf{w}_\text{s}^{\text{sh}}$.
For incompressibility, we enforce $J=\sqrt{\det\mathbf{C}}=1$ by augmenting $\psi_\text{el}$ with a Lagrange multiplier $p^\text{sh}$: 
$\psi = \psi_\text{el} - p^\text{sh}(J-1)$.

The in-plane Green--Lagrange strain $E_{\alpha\beta}$ can be split into membrane and curvature contributions as
$
    E_{\alpha\beta} = \varepsilon_{\alpha\beta} + \xi^3\kappa_{\alpha\beta}
$,
in which $\varepsilon_{\alpha\beta} = \frac{1}{2}\left( a_{\alpha\beta} - \mathring{a}_{\alpha\beta} \right)$ and $ \kappa_{\alpha\beta} = \mathring{b}_{\alpha\beta} - b_{\alpha\beta}$ are the membrane strain and change of curvature tensors on the shell midsurface, respectively.
By applying the plane-stress constraint $S^{33}=0$ together with the incompressibility condition $C_{33}=J_o^{-2}$~\cite{Kiendl15ez}, the Lagrange multiplier is obtained as $p^\text{sh}=2J_o^{-2}\,\partial\psi_\text{el}/\partial C_{33}$, which yields a reduced, purely in-plane material tangent tensor $\widehat{\mathbb{C}}$ such that $\text{d}S^{\alpha\beta}=\widehat{\mathbb{C}}^{\alpha\beta\gamma\delta}\text{d}E_{\gamma\delta}$, in which $\text{d}\mathbf{S}$ and $\text{d}\mathbf{E}$ are total differentials of $\mathbf{S}$ and $\mathbf{E}$:
\begin{align}
  S^{\alpha\beta} &= 2\frac{\partial\psi_\text{el}}{\partial C_{\alpha\beta}}-2\frac{\partial\psi_\text{el}}{\partial C_{33}}J_o^{-2} g^{\alpha\beta} \text{ ,}\label{eq:Sinc}\\
 \widehat{\mathbb{C}}^{\alpha\beta\gamma\delta} 
 &=  4\frac{\partial^2\psi_\text{el}}{\partial C_{\alpha\beta}\partial C_{\gamma\delta}} 
  + 4 \frac{\partial^2 \psi_\text{el}}{\partial C_{33}^2} J_o^{-4}  g^{\alpha\beta}g^{\gamma\delta}
  -4 \frac{\partial^2 \psi_\text{el}}{\partial C_{33}\partial C_{\alpha\beta}} J_o^{-2} g^{\gamma\delta} 
  -4 \frac{\partial^2 \psi_\text{el}}{\partial C_{33}\partial C_{\gamma\delta}} J_o^{-2} g^{\alpha\beta} \nonumber \\
  &~~~~+2\frac{\partial\psi_\text{el}}{\partial C_{33}}J_o^{-2} \left(2 g^{\alpha\beta} g^{\gamma\delta}+g^{\alpha\gamma}g^{\beta\delta}+g^{\alpha\delta}g^{\beta\gamma}\right)\text{ .}
 \label{eq:CCinc}
\end{align}
With the above equations, 3D hyperelastic constitutive models providing the first and second derivatives of $\psi_\text{el}$ with respect to $\mathbf{C}$ can be directly used for shell analysis.

In this work, we employ a transversely isotropic (Lee--Sacks) model~\cite{Lee14a,Wu18gv} to represent the MV behavior, in which the strain energy function combines a neo-Hookean term for the extracellular matrix with a convex combination of fully isotropic and transversely isotropic exponential-type terms to capture the collagen fiber network:
\begin{equation} \label{eq:LS}
    \psi_\text{el}\left(I_1,I_4\right) = \frac{c_0}{2}\left(I_1-3\right)+\frac{c_1}{2}\left(w e^{c_2\left(I_1-3\right)^2} + \left(1-w\right) e^{c_3\langle I_4-1 \rangle_+^2}-1\right)\text{ ,}
\end{equation}
in which $c_0$, $c_1$, $c_2$, $c_3$ are material parameters, $w\in\lbrack 0,1\rbrack$ determines the strain energy contribution due to anisotropy, $I_1 = \text{tr}\,\mathbf{C}$ is the first invariant of $\mathbf{C}$, $I_4 = \mathbf{m}\cdot\mathbf{C}\,\mathbf{m}$ is the squared fiber stretch, $\mathbf{m}$ is a unit vector defining the collagen fiber direction in the reference configuration, and $\langle \cdot \rangle_+ = \max\left(\cdot,0 \right)$.
The Macaulay bracket enforces the tension-only response of the collagen fibers: the anisotropic term is constant for $I_4 \le 1$, and the first and second derivatives of $\psi_\text{el}$ with respect to $I_4$ vanish under fiber compression.

A numerical challenge of exponential-type models such as Eq.~\eqref{eq:LS} is that intermediate Newton iterates can reach unrealistic strain levels at which the exponential terms cause the nonlinear solve to diverge.
To improve convergence, we impose an upper bound $(I_1^{\text{ub}}-3)$ on the invariant term $(I_1-3)$, as proposed by~\citet{JOHNSON2021113960}.
The strain energy function in Eq.~\eqref{eq:LS} is unchanged below this bound, and the bound is set well above the physiological strain levels reached by the MV leaflets so that the leaflet material response is unaffected.

\subsubsection{Leaflet prestrain}
Since living tissues such as heart valves are not stress-free in their native states, accounting for prestrains is essential to accurately capture their physiological functions~\cite{grashow2006biaixal,amini2012vivo,rausch2013effect,rausch2013mechanics}.
To account for leaflet prestrain in the human MV model, we compose the deformation gradient tensor with a prestrain tensor that maps the stress-free configuration to the imaged reference configuration.

The leaflet surface is parameterized with $\xi^1$ as the circumferential coordinate and $\xi^2$ as the radial coordinate, which is consistent with the $(u,v)$ parameterization in Section~\ref{s:parameterization}.
Let $\left(\mathring{\mathbf{e}}_1,\mathring{\mathbf{e}}_2,\mathring{\mathbf{e}}_3\right)$ denote a local orthonormal material basis aligned with the circumferential, in-plane transverse, and through-thickness directions in the imaged reference configuration.
Using the covariant and contravariant base vectors defined in Section~\ref{s:shell}, we set
$\mathring{\mathbf{e}}_1=\mathring{\mathbf g}_1/\|\mathring{\mathbf{g}}_{1}\|$,
$\mathring{\mathbf{e}}_2=\mathring{\mathbf g}^2/\|\mathring{\mathbf{g}}^{2}\|$, and
$\mathring{\mathbf{e}}_3=\mathring{\mathbf{e}}_1\times\mathring{\mathbf{e}}_2$.
We assume that the prescribed prestrain consists of principal stretches along the local material directions without rotation, which allows the same local orthonormal basis to be used to represent material directions in the stress-free and imaged reference configurations.
In this basis, the prestrain tensor is expressed as
\begin{align}
\mathbf{F}_\text{pre}
= \lambda_\text{c}\,\mathring{\mathbf{e}}_1\!\otimes\!\mathring{\mathbf{e}}_1
+ \lambda_\text{t}\,\mathring{\mathbf{e}}_2\!\otimes\!\mathring{\mathbf{e}}_2
+ (\lambda_\text{c}\lambda_\text{t})^{-1}\,\mathring{\mathbf{e}}_3\!\otimes\!\mathring{\mathbf{e}}_3\text{ ,}
\end{align}
in which $\lambda_\text{c}$ and $\lambda_\text{t}$ are the circumferential and transverse stretches, respectively~\cite{liu2023computational}.
Let $\mathbf{F}$ be the deformation gradient tensor of the shell continuum relative to the imaged reference configuration.
The effective deformation gradient tensor is then defined as
\begin{align}
    \mathbf{F}_\text{eff} = \mathbf{F} \mathbf{F}_\text{pre}\text{ ,}
\end{align}
and the corresponding effective right Cauchy--Green tensor is 
\begin{equation}\label{eq:leaflet_prestrains}
    \mathbf{C}_\text{eff} = \left( \mathbf{F} \mathbf{F}_\text{pre} \right)^\top \left( \mathbf{F} \mathbf{F}_\text{pre} \right) = \mathbf{F}_\text{pre}^\top\mathbf{C} \mathbf{F}_\text{pre}\text{ ,}
\end{equation}
so that both invariants $I_1$ and $I_4$ entering the strain energy function are evaluated using $\mathbf{C}_\text{eff}$. 
The prescribed prestrain is isochoric, $\det\left(\mathbf{F}_\text{pre}\right)=1$, and the determinant of the effective deformation gradient satisfies
\begin{equation}
    J_\text{eff} = \det\left(\mathbf{F}_\text{eff}\right) = \det\left(\mathbf{F}\right) \det\left(\mathbf{F}_\text{pre}\right) = J\text{ .}
    \label{eq:leaflet_Jeff}
\end{equation}
Thus, the incompressibility constraint remains unchanged after incorporating the prescribed prestrain.
This shifts the stress--strain response of the leaflet material to reflect the prestrain while remaining fully compatible with the shell formulation.

Let $C_{ij} = \mathring{\mathbf{e}}_i \cdot \mathbf{C}\,\mathring{\mathbf{e}}_j$ denote the components of $\mathbf C$ in the local orthonormal material basis $\left(\mathring{\mathbf{e}}_1,\mathring{\mathbf{e}}_2,\mathring{\mathbf{e}}_3\right)$.
Kirchhoff--Love kinematics with the diagonal form of $\mathbf{F}_\text{pre}$ give
\begin{align}
    \mathbf{C} = 
    \begin{bmatrix}
    C_{11} & C_{12} & 0 \\
    C_{12} & C_{22} & 0 \\
    0      & 0      & C_{33}
    \end{bmatrix}, 
\quad
    \mathbf{C}_\text{eff} = 
    \begin{bmatrix}
    \lambda_\text{c}^2 \, C_{11}               & \lambda_\text{c}\,\lambda_\text{t}\,C_{12} & 0 \\
    \lambda_\text{c}\,\lambda_\text{t}\,C_{12} & \lambda_\text{t}^2 \, C_{22}               & 0 \\
    0                                          & 0                                          & (\lambda_\text{c}\,\lambda_\text{t})^{-2}\,C_{33}
    \end{bmatrix}\text{ .} 
\end{align}
We then evaluate two invariants entering the Lee--Sacks model as
\begin{equation}\label{eq:leaflet_I1eff}
I_{1,\text{eff}} = \text{tr} \left(\mathbf{C}_\text{eff}\right) = \lambda_\text{c}^{2}\,C_{11} + \lambda_\text{t}^{2}\,C_{22} + (\lambda_\text{c}\,\lambda_\text{t})^{-2}\,C_{33}\text{ ,}
\end{equation}
and
\begin{equation}\label{eq:leaflet_I4eff}
I_{4,\text{eff}} = \mathbf{m}_0\cdot \mathbf{C}_\text{eff} \, \mathbf{m}_0\text{ .}
\end{equation}
Here, $\mathbf{m}_0$ denotes the unit collagen fiber direction defined in the stress-free material configuration.
For a general in-plane fiber direction $\mathbf{m}_0= m_\text{c}\,\mathring{\mathbf{e}}_1 +  m_\text{t}\,\mathring{\mathbf{e}}_2 +  m_\text{n}\,\mathring{\mathbf{e}}_3$ with $m_\text{c}^2 + m_\text{t}^2 =1$ and $m_\text{n}=0$, the anisotropic invariant becomes
\begin{equation}\label{eq:leaflet_I4eff_components}
I_{4,\text{eff}} = \lambda_\text{c}^{2}\,m_\text{c}^{2}\,C_{11} + 2\lambda_\text{c}\,\lambda_\text{t}\,m_\text{c}\,m_\text{t}\,C_{12} +
 \lambda_\text{t}^{2}\,m_\text{t}^{2}\,C_{22}\text{ .}
\end{equation}
Consequently, the Lee--Sacks strain energy introduced in Section~\ref{s:shell} is evaluated as
\begin{equation}
\psi_\text{el} = \frac{c_0}{2}\left(I_{1,\text{eff}}-3\right) +\frac{c_1}{2}\left(we^{c_2\left(I_{1,\text{eff}}-3\right)^2}
+(1-w)e^{c_3\langle I_{4,\text{eff}}-1\rangle_+^2}-1\right)\text{ .}
    \label{eq:leaflet_prestrain_energy}
\end{equation}

We next consider the derivatives of the strain energy with respect to $\mathbf{C}$ to show that the prestrained constitutive law remains compatible with the shell formulation in Section~\ref{s:shell}.
Define the two deformation-independent structural tensors
\begin{equation}
    \mathbf{A}_\text{pre}
    =
    \mathbf{F}_\text{pre}
    \mathbf{F}_\text{pre}^\top\text{ ,}
    \qquad
    \mathbf{M}_\text{pre}
    =
    \left(
        \mathbf{F}_\text{pre}\mathbf{m}_0
    \right)
    \otimes
    \left(
        \mathbf{F}_\text{pre}\mathbf{m}_0
    \right)\text{ .}
    \label{eq:leaflet_prestrain_structural_tensors}
\end{equation}
Using the cyclic property of the trace, Eqs.~\eqref{eq:leaflet_I1eff} and
\eqref{eq:leaflet_I4eff} can be expressed in coordinate-independent form as
\begin{equation}
    I_{1,\text{eff}}
    =
    \mathbf{C}:\mathbf{A}_\text{pre}\text{ ,}
    \qquad
    I_{4,\text{eff}}
    =
    \mathbf{C}:\mathbf{M}_\text{pre}\text{ .}
    \label{eq:leaflet_prestrain_invariant_tensor_form}
\end{equation}
Note that we recover $I_1$ and $I_4$ in Section~\ref{s:shell} if $\lambda_\text{c}=\lambda_\text{t}=1$.
Since $\mathbf{F}_\text{pre}$ is prescribed and independent of the current deformation,
\begin{equation}
    \frac{\partial I_{1,\text{eff}}}{\partial\mathbf{C}}
    =
    \mathbf{A}_\text{pre}\text{ ,}
    \qquad
    \frac{\partial I_{4,\text{eff}}}{\partial\mathbf{C}}
    =
    \mathbf{M}_\text{pre}\text{ .}
    \label{eq:leaflet_prestrain_invariant_derivatives}
\end{equation}
Additionally, the derivatives of the strain energy with respect to prestrained invariants are
\begin{align}
    \frac{\partial\psi_\text{el}}{\partial I_{1,\text{eff}}}
    &=
    \frac{c_0}{2}
    +
    c_1\,c_2\,w
    \left(
        I_{1,\text{eff}}-3
    \right)
    e^{
        c_2
        \left(
            I_{1,\text{eff}}-3
        \right)^2
    }\text{ ,} 
    \label{eq:leaflet_psi1}
    \\
    \frac{\partial\psi_\text{el}}{\partial I_{4,\text{eff}}}
    &=
    \begin{cases}
    c_1\,c_3\,(1-w)
    \left(
        I_{4,\text{eff}}-1
    \right)
    e^
    {
        c_3
        \left(
            I_{4,\text{eff}}-1
        \right)^2
    }\,, &\text{ if } I_{4,\text{eff}} > 1\text{ ,}\\
    0\,, &\text{ if } I_{4,\text{eff}} \le 1\text{ .}
    \end{cases}
    \label{eq:leaflet_psi4}
\end{align}
Application of the chain rule then yields
\begin{equation}
    \frac{\partial\psi_\text{el}}
    {\partial\mathbf{C}}
    =
     \frac{\partial\psi_\text{el}}{\partial I_{1,\text{eff}}}\mathbf{A}_\text{pre}
    +
    \frac{\partial\psi_\text{el}}{\partial I_{4,\text{eff}}}\mathbf{M}_\text{pre}\text{ .}
    \label{eq:leaflet_prestrain_dpsi}
\end{equation}
Thus, the prestrain modifies the constitutive stress response through the structural tensors
$\mathbf{A}_\text{pre}$ and $\mathbf{M}_\text{pre}$ while leaving the shell kinematics unchanged.

Because the prescribed prestrain is isochoric and contains no coupling between the in-plane and through-thickness directions, the incompressibility constraint and the plane-stress reduction derived for the Kirchhoff--Love shell retain the same kinematic form.
In particular, $C_{33}=J_o^{-2}$, and the Lagrange multiplier associated with incompressibility is obtained from the condition
$S^{33}=0$.
Additionally, because the collagen fibers are confined to the leaflet plane, the prestrained fiber direction $\mathbf{F}_\text{pre}\mathbf{m}_0$ remains in the same plane:
\begin{equation}
    M_\text{pre}^{33}
    =
    0\text{ ,}
    \qquad
    A_\text{pre}^{33}
    =
    (\lambda_\text{c}\,\lambda_\text{t})^{-2}\text{ ,}
\end{equation}
which gives
\begin{equation}
    \frac{\partial\psi_\text{el}}{\partial C_{33}}
    =
    \frac{\partial\psi_\text{el}}{\partial I_{1,\text{eff}}}\,(\lambda_\text{c}\,\lambda_\text{t})^{-2}\text{ .}
    \label{eq:ppsi_pc33}
\end{equation}
Therefore, the reduced in-plane stress can be written explicitly as
\begin{equation}
    S^{\alpha\beta}
    =
    2\frac{\partial\psi_\text{el}}{\partial I_{1,\text{eff}}}
        A_\text{pre}^{\alpha\beta}
        +
    2\frac{\partial\psi_\text{el}}{\partial I_{4,\text{eff}}}
        M_\text{pre}^{\alpha\beta}
        -
    2\frac{\partial\psi_\text{el}}{\partial I_{1,\text{eff}}}\,(\lambda_\text{c}\,\lambda_\text{t})^{-2}\,
        J_o^{-2}\,\bar{C}^{\alpha\beta}\text{ ,}
    \label{eq:leaflet_prestrain_shell_stress_explicit}
\end{equation}
in which $\bar{\mathbf{C}}$ is the inverse of $\mathbf{C}$.
Note that $\bar{C}^{\alpha\beta} = g^{\alpha\beta}$ when the components are expressed in the curvilinear basis~\cite{Kiendl15ez}.

To construct the consistent material tangent, we evaluate the second derivatives of the strain energy with respect to the effective invariants as follows:
\begin{align}
    \frac{\partial^2\psi_\text{el}}{\partial I_{1,\text{eff}}^2}
    &=
    c_1\,c_2\,w
    e^
    {
        c_2
        \left(
            I_{1,\text{eff}}-3
        \right)^2
    }
    \left(
        1
        +
        2c_2
        \left(
            I_{1,\text{eff}}-3
        \right)^2
    \right)\text{ ,}
    \label{eq:leaflet_psi11}
    \\
    \frac{\partial^2\psi_\text{el}}{\partial I_{4,\text{eff}}^2}
    &=
    \begin{cases}
    c_1\,c_3\,(1-w)
    e^
    {
        c_3
        \left(
            I_{4,\text{eff}}-1
        \right)^2
    }
    \left(
        1
        +
        2c_3
        \left(
            I_{4,\text{eff}}-1
        \right)^2
    \right)\,, &\text{ if } I_{4,\text{eff}} > 1\text{ ,}\\
    0\,, &\text{ if } I_{4,\text{eff}} \le 1\text{ ,}
    \end{cases}\label{eq:leaflet_psi44}\\
    \frac{\partial^2\psi_\text{el}}{\partial I_{1,\text{eff}}\,\partial I_{4,\text{eff}}}&=0\text{ ,}
    \label{eq:leaflet_psi41}
\end{align}
since the Lee--Sacks strain energy contains no mixed term between $I_1$ and $I_4$.
Moreover, both invariants are linear with respect to $\mathbf{C}$ for a prescribed $\mathbf{F}_\text{pre}$.
The second derivative of the elastic strain energy is therefore
\begin{equation}
    \mathbb{H}
    =
    \frac{\partial^2\psi_\text{el}}
    {\partial\mathbf{C}^2}
    =
    \frac{\partial^2\psi_\text{el}}{\partial I_{1,\text{eff}}^2}\,
    \mathbf{A}_\text{pre}\otimes\mathbf{A}_\text{pre}
    +
    \frac{\partial^2\psi_\text{el}}{\partial I_{4,\text{eff}}^2}\,
    \mathbf{M}_\text{pre}\otimes\mathbf{M}_\text{pre}\text{ .}
    \label{eq:leaflet_prestrain_d2psi}
\end{equation}
The components required by the plane-stress reduction are thus
\begin{align}
    H^{\alpha\beta\gamma\delta}
    &=
    \frac{\partial^2\psi_\text{el}}{\partial I_{1,\text{eff}}^2}\,
    A_\text{pre}^{\alpha\beta}\,
    A_\text{pre}^{\gamma\delta}
    +
    \frac{\partial^2\psi_\text{el}}{\partial I_{4,\text{eff}}^2}\,
    M_\text{pre}^{\alpha\beta}\,
    M_\text{pre}^{\gamma\delta}\text{ ,}
    \\
    H^{33\alpha\beta}
    &=
    \frac{\partial^2\psi_\text{el}}{\partial I_{1,\text{eff}}^2}
    (\lambda_\text{c}\,\lambda_\text{t})^{-2}
    A_\text{pre}^{\alpha\beta}\text{ ,}
    \\
    H^{3333}
    &=
    \frac{\partial^2\psi_\text{el}}{\partial I_{1,\text{eff}}^2}(\lambda_\text{c}\,\lambda_\text{t})^{-4}\text{ ,}
    \label{eq:leaflet_prestrain_H_components}
\end{align}
in which the fiber terms vanish from the latter two expressions since $M_\text{pre}^{33}=0$.
Substitution of Eqs.~\eqref{eq:leaflet_prestrain_dpsi} and
\eqref{eq:leaflet_prestrain_d2psi}--\eqref{eq:leaflet_prestrain_H_components} into the reduced material tangent in Eq.~\eqref{eq:CCinc} gives
\begin{align}
    \widehat{\mathbb{C}}^{\alpha\beta\gamma\delta}
    ={}&
    4\frac{\partial^2\psi_\text{el}}{\partial I_{1,\text{eff}}^2}\,
        A_\text{pre}^{\alpha\beta}\,
        A_\text{pre}^{\gamma\delta}
        +
    4\frac{\partial^2\psi_\text{el}}{\partial I_{4,\text{eff}}^2}\,
        M_\text{pre}^{\alpha\beta}\,
        M_\text{pre}^{\gamma\delta}
    \nonumber\\
    &+
    4\frac{\partial^2\psi_\text{el}}{\partial I_{1,\text{eff}}^2}(\lambda_\text{c}\,\lambda_\text{t})^{-4}
    J_o^{-4}
    \bar{C}^{\alpha\beta}\bar{C}^{\gamma\delta}
    \nonumber\\
    &-
    4\frac{\partial^2\psi_\text{el}}{\partial I_{1,\text{eff}}^2}(\lambda_\text{c}\,\lambda_\text{t})^{-2}
    J_o^{-2}A_\text{pre}^{\alpha\beta}\,
        \bar{C}^{\gamma\delta}
    -
    4\frac{\partial^2\psi_\text{el}}{\partial I_{1,\text{eff}}^2}(\lambda_\text{c}\,\lambda_\text{t})^{-2}
    J_o^{-2}A_\text{pre}^{\gamma\delta}\,
        \bar{C}^{\alpha\beta}
    \nonumber\\
    &+
    2\frac{\partial\psi_\text{el}}{\partial I_{1,\text{eff}}}(\lambda_\text{c}\,\lambda_\text{t})^{-2}
    J_o^{-2}
    \left(
        2\bar{C}^{\alpha\beta}\bar{C}^{\gamma\delta}
        +
        \bar{C}^{\alpha\gamma}\bar{C}^{\beta\delta}
        +
        \bar{C}^{\alpha\delta}\bar{C}^{\beta\gamma}
    \right)\text{ .} 
    \label{eq:leaflet_prestrain_reduced_tangent}
\end{align}
Therefore, the same three-dimensional Lee--Sacks constitutive model and its consistent linearization~\cite{Wu18gv} can be used directly within the Kirchhoff--Love shell formulation after replacing the standard invariants by their prestrain-modified counterparts.

An important consequence of this formulation is that the reference configuration constructed by in-vivo imaging is not required to be stress-free.
When no additional deformation is applied relative to the imaged geometry, i.e., $\mathbf{F}=\mathbf{I}$, the effective deformation gradient tensor is generally different from the identity tensor.
The leaflet can therefore carry a nonzero initial stress in the reconstructed in-vivo configuration, while the imaged surface itself remains the computational reference geometry.
This shifts the material response into the physiological strain range without requiring an explicit reconstruction of the unknown stress-free leaflet configuration.

\subsubsection{Cable formulation for the chordae tendineae}\label{s:cable}
To accurately characterize the internal mechanics of chordal structures, one might initially consider comprehensive 3D solid elements. 
Physiologically, CTs are thin collagenous structures whose primary mechanical role is to transmit tensile loads between the valve leaflets and papillary muscles, and they are commonly idealized as cable or truss structures that carry load primarily in tension with negligible bending~\cite{lee2015effects, alleau2019use, gaidulis2022patient}.
Therefore, the structure can be accurately modeled using an extension of the IGA cable formulation~\cite{RAKNES2013127,JOHNSON2021113960}, which directly replicates the macroscopic mechanical function of the MV chordal structures without the computational burden of solid elements.

In this work, we employ a St.~Venant--Kirchhoff model with zero Poisson's ratio to represent the CT behavior, following the cable formulation of~\citet{JOHNSON2021113960}.
Although this model does not reproduce the nonlinear stress--strain response of soft tissues over the full strain range, it effectively captures the linear stress--strain relationship of the CT in its stiffer post-transition tensile regime.
By assuming the reference configuration already accounts for the pre-transition stretch, the tensile stiffness is governed entirely by the Young's modulus $E^\text{ca}$.

The IGA cable theory assumes that the mechanics are expressed entirely on a one-dimensional center curve parameterized by a single coordinate $\xi^1$. 
Then a material point on the cable centerline $(\mathcal{L}^\text{ca})_t$ in the deformed configuration is expressed by $\mathbf{x}^{\text{ca}}(\xi^1)$, and the displacement $\mathbf{y}^{\text{ca}}$ is defined relative to the reference configuration $\mathring{\mathbf{x}}^{\text{ca}}(\xi^1)$ of the center curve  $(\mathcal{L}^\text{ca})_0$ as
\begin{equation}
\mathbf{y}^{\text{ca}}\left(\mathring{\mathbf{x}}^{\text{ca}}\left(\xi^1\right)\right) = \mathbf{x}^{\text{ca}}\left(\xi^1\right) - \mathring{\mathbf{x}}^{\text{ca}}\left(\xi^1\right)\text{ .} 
\end{equation}
By using the 1D parametric description, the covariant base vector in the reference configuration is defined as $\mathring{\mathbf{g}}_1 = \partial \mathring{\mathbf{x}}^{\text{ca}}/\partial \xi^1$, and the corresponding contravariant base vector is defined as $\mathring{\mathbf{g}}^1 = \mathring{\mathbf{g}}_1/\left\Vert \mathring{\mathbf{g}}_1 \right\Vert^2$.
We denote the base vectors in the deformed configuration without $\mathring{\left(\cdot \right)}$.
Due to the small bending stiffness of the CT, we only consider the axial deformation and neglect both bending and torsional strains described in~\citet{RAKNES2013127}.
Then the axial strain is computed by
\begin{equation}
    \varepsilon = \frac{1}{2}\left(\mathbf{g}_1 \cdot \mathbf{g}_1 - \mathring{\mathbf{g}}_1 \cdot \mathring{\mathbf{g}}_1\right) \text{ .}
\end{equation}
The weak form of the cable subproblem is defined as
\begin{align}\label{eq:b2-cable}
\nonumber B^\text{ca}_\text{s}(\mathbf{w}_\text{s}^{\text{ca}},\mathbf{y}^{\text{ca}}) - F^\text{ca}_\text{s}(\mathbf{w}_\text{s}^{\text{ca}}) 
=& \int_{(\mathcal{L}^\text{ca})_0} \mathbf{w}_\text{s}^{\text{ca}} \cdot \rho^\text{ca}_\text{s} \ A_0^\text{ca} \ \frac{\partial^2 \mathbf{y}^{\text{ca}}}{\partial t^2} \,\text{d}\mathcal{L}
+ \int_{(\mathcal{L}^\text{ca})_0} \delta\varepsilon \ E^\text{ca} A_0^\text{ca}  \left\Vert \mathring{\mathbf{g}}^1 \right\Vert^4 \varepsilon \,\text{d}\mathcal{L} \\
\nonumber &+ \int_{(\mathcal{L}^\text{ca})_0} \mathbf{w}_\text{s}^{\text{ca}} \cdot c^\text{ca} \ A_0^\text{ca} \ \rho^\text{ca}_\text{s} \ \frac{\partial \mathbf{y}^{\text{ca}}}{\partial t} \,\text{d}\mathcal{L} \\
& 
- \int_{(\mathcal{L}^\text{ca})_0} \mathbf{w}_\text{s}^{\text{ca}}  \cdot \rho^\text{ca}_\text{s} \ A_0^\text{ca} \ \mathbf{f}^\text{ca}_\text{s} \,\text{d}\mathcal{L} 
 - \int_{(\mathcal{L}^\text{ca})_t} \mathbf{w}_\text{s}^{\text{ca}}  \cdot \mathbf{h}^\text{ca}_\text{s} \,\text{d}\mathcal{L}\text{ ,}
\end{align}
in which $\rho^\text{ca}_\text{s}$ is the density, $\delta\varepsilon$ is the variation of $\varepsilon$ with respect to $\mathbf{y}^{\text{ca}}$ in the direction of $\mathbf{w}_\text{s}^{\text{ca}}$, $E^\text{ca}$ is the Young's modulus of the cable structure, $A_0^\text{ca}$ is the reference cable cross-sectional area, $c^\text{ca}$ is a mass-proportional damping coefficient, $\mathbf{f}^\text{ca}_\text{s}$ is a body force, and $\mathbf{h}^\text{ca}_\text{s}$ is the traction prescribed on $\left(\mathcal{L}^\text{ca}\right)_t$.

\subsubsection{Penalty coupling for the shell and cables}\label{s:ss_coupling}
The MV leaflets and the MVCTs are modeled as separate structures, and the load transfer between them occurs through the chordal insertions.
To accommodate the non-conforming discretizations between shells and cables, we use a pointwise penalty coupling approach to enforce kinematic compatibility at the chordal insertion locations. 
The following term is added to Eq.~\eqref{eq:bs-full}:
\begin{equation}
\label{eq:displacement_penalty_sb}
+ \sum_{i=1}^{N_\text{I}}
\alpha_{i}^\text{sc} 
\left(\mathbf{w}_\text{s}^{\text{sh}}\left(\mathring{\mathbf{x}}^{\text{sh}}_i\right) - \mathbf{w}_\text{s}^{\text{ca}}\left(\mathring{\mathbf{x}}^{\text{ca}}_i\right)\right)\cdot\left(\mathbf{x}^{\text{sh}}\left(\mathring{\mathbf{x}}^{\text{sh}}_i\right) - \mathbf{x}^{\text{ca}}\left(\mathring{\mathbf{x}}^{\text{ca}}_i\right)\right)\text { .} 
\end{equation}
Eq.~\eqref{eq:displacement_penalty_sb} constrains the positions at the $N_\text{I}$ chordal insertions of a given MV, in which each pair $\left(\mathring{\mathbf{x}}_i^\text{sh},\mathring{\mathbf{x}}_i^\text{ca}\right)$ identifies which leaflet material point $\mathring{\mathbf{x}}_i^\text{sh}$ a given chordal end point $\mathring{\mathbf{x}}_i^\text{ca}$ is attached to.
The pairing is determined once in the reference configuration and is held fixed throughout the simulation, while the penalty constraint is evaluated using the current positions of the paired material points. 
This weakly enforces positional coincidence at each chordal insertion and ensures consistent load transfer between the MVCTs and the valve leaflets during deformation.
The penalty parameter $\alpha_{i}^\text{sc}$ is determined by
\begin{equation}
\label{eq:shell-beam-penalty}
\alpha_{i}^\text{sc} = c^\text{sc}\min\left(\frac{E^\text{sh} \, \thickness}{1-(\nu^\text{sh})^2},\,\frac{E^\text{ca}  A_0^\text{ca}}{h^\text{ca}_i}\right)\text{ ,}
\end{equation}
in which $c^\text{sc}$ is a dimensionless penalty coefficient, $E^\text{sh}$ is an effective material stiffness with units of stress, calculated based on the shear modulus of the neo-Hookean term for the hyperelastic material model of the leaflets, $\nu^\text{sh}$ is the Poisson's ratio of the shell, and $h^\text{ca}_i$ is the length of the cable element at insertion $i$.
By dimensional analysis, the shell contribution is the membrane stiffness of the leaflet and is insensitive to the element size, whereas the cable contribution is the axial stiffness of the attached element and scales inversely with its length.
Taking the minimum of the two stiffness-based values ensures proper constraint enforcement without causing ill-conditioning associated with the softer material.
In our numerical tests, we set the dimensionless penalty coefficient to $c^\text{sc}=4\times10^4$.

\subsubsection{Contact formulation}\label{s:contacts}
MV dynamics involve complex contact interactions, including leaflet-to-leaflet, leaflet-to-chordae, and chordae-to-chordae contacts.
To penalize interpenetration of various structural components, we adopt a nonlocal contact formulation~\cite{kamensky2018a, JOHNSON2021113960} that considers all structural parts to be a single body in the reference configuration $\Omega_0$. 
For contact between two material points, ${\mathbf{x}}^\text{A}$ and ${\mathbf{x}}^\text{B}$, the following contact term is added to Eq.~\eqref{eq:bs-full}:
\begin{equation}
\label{eq:contact}
    +\int_{\Omega_0\backslash B_R(\mathring{\mathbf{x}}^\text{A})} \int_{\Omega_0 } \delta \mathbf{r}^\text{AB}\cdot \phi'_\text{c}\left(r^\text{AB}\right) \frac{\mathbf{r}^\text{AB}}{r^\text{AB}} \text{d} \mathring{\mathbf{x}}^\text{A} \text{d} \mathring{\mathbf{x}}^\text{B} \text{ ,}
\end{equation}
in which $\mathring{\mathbf{x}}^\text{A}$ and $\mathring{\mathbf{x}}^\text{B}$ are the material points in the reference configuration, $B_R\left(\mathring{\mathbf{x}}^\text{A}\right)$ is the Euclidean ball of radius $R$ centered at $\mathring{\mathbf{x}}^\text{A}$,
$\mathbf{r}^\text{AB} = \mathbf{x}^\text{B} - \mathbf{x}^\text{A}$, $r^\text{AB} = \Vert \mathbf{r}^\text{AB} \Vert$, $\delta\mathbf{r}^\text{AB} = \mathbf{w}_\text{s}\left(\mathring{\mathbf{x}}^\text{B}\right) - \mathbf{w}_\text{s}\left(\mathring{\mathbf{x}}^\text{A}\right)$ is the variation of $\mathbf{r}^\text{AB}$ with respect to $\mathbf{y}$ in the direction of $\mathbf{w}_\text{s}$, and $\phi'_\text{c}\left(r^\text{AB}\right)$ is a contact kernel.

Let $s = r^{\mathrm{AB}} - r_{\mathrm{min}}$, $s_{\mathrm{out}} = r_{\mathrm{max}}- r_{\mathrm{min}}$, and $s_{\mathrm{in}} = \eta s_{\mathrm{out}}$, in which $r_{\mathrm{min}}$ is the minimum admissible separation, $r_{\mathrm{max}}$ is the contact cutoff, and $0 < \eta < 1$ controls the transition between two branches.
Then the contact kernel is defined as
\begin{equation}\label{eq:contact_kern}
    \phi'_{\mathrm{c}}\!\left(s\right)
=
-\begin{cases}
\infty,                                & s \leq 0\text{ ,} \\[6pt]
\frac{k_{\mathrm{c}}}{s^{p_\text{c}}}-C_2,        & 0<s<s_{\mathrm{in}}\text{ ,} \\[6pt]
C_1 \left(s-s_{\mathrm{out}}\right)^2, & s_{\mathrm{in}}\leq s<s_{\mathrm{out}}\text{ ,} \\[6pt]
0,                                     & s\geq s_{\mathrm{out}}\text{ ,}
\end{cases}
\end{equation}
in which $k_\text{c}$ is a penalty parameter and $p_\text{c}$ is the exponent controlling the inverse-power inner branch.
$C_1$ and $C_2$ are uniquely determined by smoothness and continuity constraints as follows:
\begin{equation}
\label{eq:contact_constants}
    C_1 = \frac{p_\text{c} k_{\mathrm{c}}} {2s_{\mathrm{in}}^{p_\text{c}+1} \left(s_{\mathrm{out}}-s_{\mathrm{in}}\right)}\text{ ,}
    \qquad
    C_2 = \frac{k_{\mathrm{c}}}{s_{\mathrm{in}}^{p_\text{c}}} - C_1 \left(s_{\mathrm{in}}-s_{\mathrm{out}}\right)^2\text{ .} 
\end{equation}
This ensures continuity of the force and its derivative at $s = s_{\mathrm{in}}$, and both vanish smoothly at $s = s_{\mathrm{out}}$.

In practice, the nonlocal integral in Eq.~\eqref{eq:contact} is evaluated at a discrete set of contact points on the structure, whose weights include the shell thickness or the cable cross-sectional area so that surface and line quadrature is converted into volume quadrature~\cite{kamensky2018a}.
If these coincide with the structural Gauss points, the contact resolution is tied to the quadrature rule, and enforcing no penetration between contact regions requires a denser quadrature, increasing the cost of the structural assembly.
Following~\citet{JOHNSON2021113960}, we instead distribute the contact points uniformly and decouple them from the quadrature, which allows more contact points to be added as needed to enforce no penetration without unnecessarily increasing the number of structural Gauss points.

\subsubsection{Penalty coupling for TEER}\label{s:penalty_TEER}
MV TEER induces highly non-physiological focal stress concentrations on the leaflets and alters the hemodynamics in the left heart.
To model TEER within the IMGA model, we first identify the area on the leaflet surface in which the clip device is implanted and apply a clip-equivalent force to replicate the grasping action, rather than explicitly modeling the clipping device for computational efficiency and simplicity.

We construct a pair of clip patches at the intended leaflet-grasping locations as illustrated in Figure~\ref{f:overview_figure}e. 
These patches serve solely as geometric templates and introduce no additional structural degrees of freedom. 
A prescribed point-to-point correspondence is established between the two auxiliary patches, and each pair of corresponding points is projected onto the closest points of the respective leaflet surfaces in the reference configuration. 
The resulting sets of projected leaflet material points define the two clipped leaflet regions, $\mathcal{S}_{\mathrm{clip}}^{A}$ and $\mathcal{S}_{\mathrm{clip}}^{B}$.
Then the correspondence between the two clip areas can be defined as a one-to-one mapping: 
\begin{align}
\Phi_{\mathrm{clip}}: \mathcal{S}_{\mathrm{clip}}^A \rightarrow \mathcal{S}_{\mathrm{clip}}^B\text{ ,}
\quad \Phi_{\mathrm{clip}} \left( \mathring{\mathbf{x}}^{\mathrm{sh},A} \right) = \mathring{\mathbf{x}}^{\mathrm{sh},B}\text{ ,}
\end{align}
in which $\mathring{\mathbf{x}}^{\mathrm{sh},A} \in \mathcal{S}_{\mathrm{clip}}^A$ and $\mathring{\mathbf{x}}^{\mathrm{sh},B} \in \mathcal{S}_{\mathrm{clip}}^B$ are paired leaflet material points in the reference configuration.
Note that the resulting correspondences are held fixed throughout the simulation.

The clipping constraint is introduced by a penalty potential:
\begin{align}
\label{eq:TEER_penalty_energy}
\Pi_{\mathrm{clip}} = \frac{1}{2} \int_{\mathcal{S}_{\mathrm{clip}}^A} \alpha_{\mathrm{clip}} \left( \|\mathbf{d}_{\mathrm{clip}}\| -g_{\mathrm{clip}} \right)^2 \,\mathrm{d}\mathcal{S}\text{ ,}
\end{align}
in which $\alpha_\text{clip}$ is a penalty parameter for clipping, $g_{\mathrm{clip}}$ is the target gap between two clip surfaces, 
\begin{align}
\mathbf{d}_{\mathrm{clip}} \left( \mathring{\mathbf{x}}^{\mathrm{sh},A} \right) = \mathbf{x}^{\mathrm{sh},A} \left( \mathring{\mathbf{x}}^{\mathrm{sh},A} \right) - \mathbf{x}^{\mathrm{sh},B} \circ\Phi_{\mathrm{clip}} \left( \mathring{\mathbf{x}}^{\mathrm{sh},A} \right)
\end{align}
is the current separation between paired points, and $\circ$ denotes the function composition operator.
A nonzero value of $g_{\mathrm{clip}}$ is introduced to represent the finite thickness of the clip device components clamped between the grasped leaflets.
It also keeps the clipping constraint consistent with the leaflet-to-leaflet contact treatment, which holds the coapting leaflets at a finite distance.
Enforcing $\|\mathbf{d}_{\mathrm{clip}}\|=0$ would place the clipping penalty in direct competition with the contact force at the edge of the clipped region.

The first variation of Eq.~\eqref{eq:TEER_penalty_energy} with respect to the leaflet positions gives
\begin{align}
\delta\Pi_{\mathrm{clip}}=\int_{\mathcal{S}_{\mathrm{clip}}^A}\alpha_{\mathrm{clip}}\left(\left\|\mathbf{d}_{\mathrm{clip}}\right\|-g_{\mathrm{clip}}\right)\delta\left\|\mathbf{d}_{\mathrm{clip}}\right\|\,\mathrm{d}\mathcal{S}\text{ .}
\label{eq:TEER_penalty_variation}
\end{align}
Since $\Phi_{\mathrm{clip}}$ is fixed throughout the simulation, its variation vanishes, and the variation of the separation vector is
\begin{equation}
\delta\mathbf{d}_{\mathrm{clip}} = \mathbf{w}_{\mathrm{s}}^{\mathrm{sh},A} - \mathbf{w}_{\mathrm{s}}^{\mathrm{sh},B} \circ\Phi_{\mathrm{clip}}\text{ .}
\end{equation}
Using the chain rule for the Euclidean norm gives
\begin{equation}
\delta \left\|\mathbf{d}_{\mathrm{clip}}\right\| = \frac{\mathbf{d}_{\mathrm{clip}}} {\left\|\mathbf{d}_{\mathrm{clip}}\right\|} \cdot \delta\mathbf{d}_{\mathrm{clip}}\text{ .}
\end{equation}
Therefore, the following TEER penalty contribution is added to Eq.~\eqref{eq:bs-full}:
\begin{align}
\label{eq:TEER_penalty}
+\,\alpha_\text{clip} \int_{\mathcal{S}_{\mathrm{clip}}^A} \delta\mathbf{d}_{\mathrm{clip}} \cdot
\left( 1-\frac{g_{\mathrm{clip}}}{\|\mathbf{d}_{\mathrm{clip}}\|} \right) \mathbf{d}_{\mathrm{clip}}\,\mathrm{d}\mathcal{S}\text{ .}
\end{align}
Eq.~\eqref{eq:TEER_penalty} constrains the distance between the paired current leaflet positions to the prescribed value $g_{\mathrm{clip}}$.
The penalty parameter $\alpha_\text{clip}$ is determined by
\begin{align}
    \alpha_\text{clip} = c_{\mathrm{clip}} \frac{E^{\mathrm{sh}} \thickness}{ {A_{\mathrm{clip}}} \left(1-(\nu^{\mathrm{sh}})^2\right)}\text{ ,}
\end{align}
in which $c_{\mathrm{clip}}$ is a dimensionless coefficient and $A_{\mathrm{clip}}$ is the reference area of the clipping region.
In our numerical study, we set the dimensionless penalty coefficient to $c_\text{clip}=3.5\times10^4$.
Additionally, the clipping penalty is introduced gradually during the simulation to avoid instability arising from sudden coupling between the clip surfaces. 

Since the clipping constraint is defined locally through a pair of clip patches, the present formulation can model multiple clip configurations without modifying the leaflet mesh topology or adding structural degrees of freedom. 
Additionally, because the penalty potential depends only on the magnitude of the separation vector, it does not naturally distinguish the relative side of one leaflet with respect to the other. 
In principle, the paired leaflet regions can pass through one another and recover the target distance on the opposite side. 
To prevent such nonphysical penetration, leaflet-to-leaflet contact remains active within the clipped regions.
The contact formulation provides a nonpenetration constraint, while the clipping penalty maintains the prescribed finite separation.
To keep the two contributions from competing, the target gap is chosen such that $g_{\mathrm{clip}} > r_{\mathrm{max}}$, so that the contact force vanishes at the prescribed separation and is engaged only when a paired region is driven toward penetration.
Together, these two contributions preserve the relative arrangement of the grasped leaflet regions and prevent both complete overlap and penetration during deformation.


\subsection{Immersogeometric discretization}\label{subsec:dis_be}
In this section, we discuss the immersogeometric discretization of the FSI problem, in which the subproblems are discretized independently and coupled through interfacial constraint equations.
Additionally, to accurately couple the independently discretized fluid and solid subproblems, we propose a subcell quadrature of the structural mesh.
This allows us to select an appropriate mesh level for each domain while maintaining robust interaction through the interfacial constraint equations.

\subsubsection{Fluid and structure subproblems}
The ALE fluid subproblem is discretized by the variational multiscale (VMS) method~\cite{Hughes00a, Hughes01b, Bazilevs07b}, with some modifications of the stabilization parameters to improve mass conservation as described in~\citet[Section 4.4]{Kamensky15ch}.
The ALE--VMS formulation is particularly suitable for the current MV TEER FSI simulations, as it provides a stabilized discretization that can handle both laminar and turbulent flow regimes. 
Within the fluid domain, the mesh velocity $\hat{\mathbf{u}}$ is determined by solving a fictitious linear elastostatic problem with displacement boundary conditions imposed by the motion of the fluid wall~\cite{Tezduyar93a, Johnson94a, Stein01f, Stein03a, Bazilevs08a}. 
The shell and cable structures are discretized by IGA with a Galerkin approach, using at least $C^1$-continuous B-splines to represent the structural geometries.
The same B-splines also serve as basis functions for the structural discrete spaces, naturally providing the $C^1$ continuity between elements required by the Kirchhoff--Love shell formulation.

\subsubsection{Fluid--structure kinematic constraints}\label{s:kinematic_constraints}

The motion of the atrial and ventricular walls is prescribed and imposed directly as a displacement boundary condition on the ALE fluid mesh.
The leaflet surface, by contrast, arbitrarily intersects the fluid mesh, and the fluid--shell coupling employs both Lagrange multiplier and penalty methods to enforce the kinematic constraints.
Following the weak imposition of the no-slip condition at the fluid--shell interface~\cite{Bazilevs07c,BazHsu12b}, we formally eliminate the tangential component of the Lagrange multiplier $\pmb{\lambda}$, leaving only a penalty term to enforce the no-slip condition weakly.
Because of the large pressure gradients across the leaflets, we retain the normal component of the Lagrange multiplier to strengthen the enforcement of the no-penetration condition at the fluid--shell interface:
\begin{equation}
\lambda = \pmb{\lambda} \cdot \mathbf{n}^\text{sh}_\text{s}\text{ ,}
\end{equation}
in which $\mathbf{n}^\text{sh}_\text{s}$ denotes the shell midsurface normal.
The chordae tendineae are not directly coupled to the blood flow in the present formulation because their slender cross sections make the mean hydrodynamic loading small relative to the tensile loads transmitted from the leaflets.
However, omitting the fluid--cable interaction removes the viscous resistance that would damp the transverse motion of the chordae structures.
Therefore, a damping term is included in the cable formulation~\eqref{eq:b2-cable} as a phenomenological representation of this unresolved hydrodynamic dissipation, which suppresses nonphysical high-frequency chordal oscillations.

The discrete variational equation for this FSI system is summarized as follows: Find $\mathbf{u}^h_\text{f}\in\mathcal{S}_u^h$, $p^h\in\mathcal{S}_p^h$, $\mathbf{y}^h\in\mathcal{S}_y^h$, and the scalar multiplier $\lambda\in\mathcal{S}_\ell^h$ such that for all $\mathbf{w}_\text{f}^h\in\mathcal{V}_u^h$, $q^h\in\mathcal{V}_q^h$, $\mathbf{w}_\text{s}^h\in\mathcal{V}_y^h$, and $\delta\lambda\in\mathcal{V}_\ell^h$,
\begin{align}
    \label{eq:general-al-variational2}
    \nonumber &B_\text{f}^{\text{VMS}}\left(\{\mathbf{w}^h_\text{f},q^h\},\{\mathbf{u}^h_\text{f},p^h\}\right)-F_\text{f}^{\text{VMS}}\left(\{\mathbf{w}^h_\text{f},q^h\}\right) + B_\text{s}\left(\mathbf{w}^h_\text{s},\mathbf{y}^h\right)-F_\text{s}\left(\mathbf{w}^h_\text{s}\right) \\
    \nonumber&~~+ \int_{\mathcal{S}^\text{sh}}\left(\mathbf{w}_\text{f}^h-\mathbf{w}_\text{s}^{\text{sh},h}\right)\cdot\lambda\mathbf{n}_\text{s}^\text{sh}~\text{d}\mathcal{S} \\
    \nonumber&~~+\int_{\mathcal{S}^\text{sh}}\delta\lambda\left(\left(\mathbf{u}_\text{f}^h-\mathbf{u}_\text{s}^{\text{sh},h}\right)\cdot\mathbf{n}_\text{s}^\text{sh} - \frac{r}{\beta^\text{sh}_{\text{NOR}}}\lambda\right)~\text{d}\mathcal{S} \\
    \nonumber&~~+ \int_{\mathcal{S}^\text{sh}}\left(\mathbf{w}^h_\text{f}-\mathbf{w}_\text{s}^{\text{sh},h}\right)\cdot\beta^\text{sh}_{\text{NOR}}\left(\left(\mathbf{u}_\text{f}^h-\mathbf{u}_\text{s}^{\text{sh},h}\right)\cdot\mathbf{n}_\text{s}^\text{sh}\right)\mathbf{n}_\text{s}^\text{sh}~\text{d}\mathcal{S} \\
    &~~+ \int_{\mathcal{S}^\text{sh}}\left(\mathbf{w}_\text{f}^h-\mathbf{w}_\text{s}^{\text{sh},h}\right)\cdot\beta^\text{sh}_{\text{TAN}}\left(\left(\mathbf{u}_\text{f}^h-\mathbf{u}_\text{s}^{\text{sh},h}\right) - \left(\left(\mathbf{u}_\text{f}^h-\mathbf{u}_\text{s}^{\text{sh},h}\right)\cdot\mathbf{n}_\text{s}^\text{sh}\right)\mathbf{n}_\text{s}^\text{sh}\right)~\text{d}\mathcal{S} = 0\text{ ,}
\end{align}
in which $B_{\text{f}}^{\text{VMS}}$ and $F_{\text{f}}^{\text{VMS}}$ are the VMS discretizations of $B_\text{f}$ and $F_{\text{f}}$, respectively, with their explicit forms given in~\citet{Kamensky15ch}, the superscript $h$ denotes the corresponding variable in the discrete space, $\mathcal{S}^\text{sh}$ is the fluid--shell interface, $\beta^\text{sh}_\text{NOR}$ and $\beta^\text{sh}_\text{TAN}$ are penalty parameters for the normal and tangential components of the fluid--shell kinematic constraint, respectively, and $r \geq 0$ is a dimensionless constant.
Note that the $r$-related term regularizes the no-penetration constraint~\cite{Kamensky15gr}, avoiding the need to satisfy an inf--sup condition and corresponding to a limiting case of strongly consistent Barbosa--Hughes stabilization~\cite{Barbosa91}.

\subsubsection{Subcell quadrature}
\label{s:subcell_discretization}
The fluid--structure kinematic constraints must be enforced at appropriate discrete spatial locations to accurately evaluate the interfacial integrals over $\mathcal{S}^\text{sh}$ in Eq.~\eqref{eq:general-al-variational2}.
In conventional IMGA, the kinematic constraints between the fluid and the shell are enforced directly at the structural quadrature points, that is, the points used to integrate the shell formulation. 
For cubic shell elements, this corresponds to $4 \times 4$ in-plane Gauss--Legendre points per element.
However, the interfacial integrals are sampled only at these quadrature points, and the coupling terms therefore contribute only to the fluid degrees of freedom whose basis function supports contain those points.
Consequently, as the background fluid mesh is refined, the fixed spatial density of the structural quadrature points leaves an increasing number of fluid degrees of freedom near the interface without any contribution from the coupling terms.
This resolution mismatch compromises the no-penetration condition, leading to unphysical fluid leakage through the valve leaflets or an inaccurate transfer of fluid tractions to the shell surface.

\begin{figure}[!t]\centering
  \centering
    \includegraphics[width=0.9\textwidth]{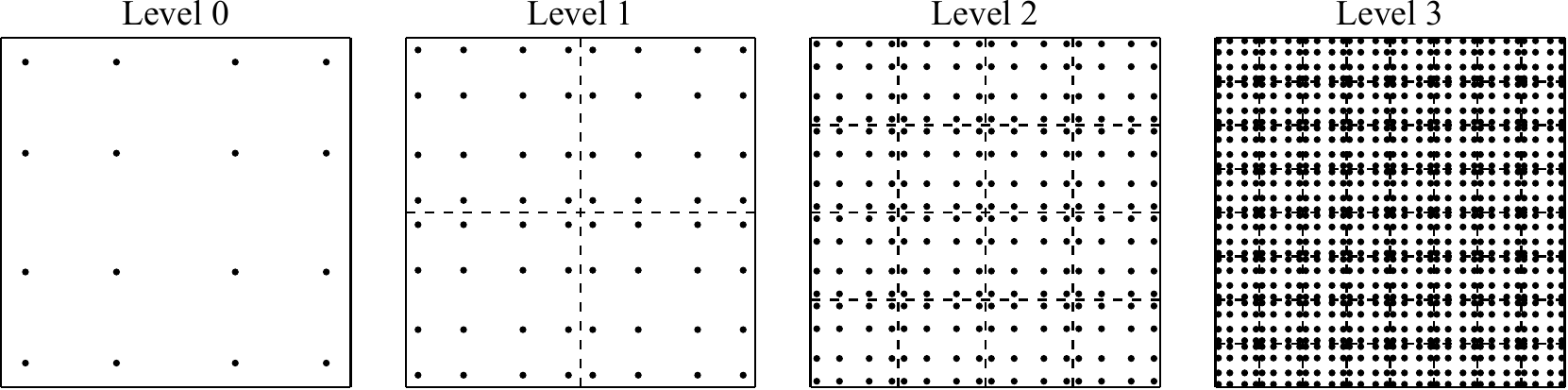}
  \caption{Subcell point distribution within a single element under recursive subdivision.}\label{f:Subcell}
\end{figure}

To address this challenge without over-refining the structural mesh, we propose a subcell quadrature rule for the immersed structure. 
As illustrated in Figure~\ref{f:Subcell}, each shell element is recursively subdivided to generate a localized distribution of dedicated subcell points. 
These subcell points serve only the fluid--structure coupling and are decoupled from the structural quadrature points, whose number is unchanged regardless of the subcell division level.
Since the fluid interacts exclusively with the subcell points, unphysical leakage through the interfaces can be mitigated without refining the structural mesh or increasing the number of structural quadrature points.
Once the kinematic constraint residuals are evaluated and the associated multipliers are updated at the subcell quadrature points, the resulting constraint tractions are directly assembled into the shell control-point force vector using the IGA basis functions evaluated at the parametric coordinates of those subcell points, without a separate projection step.

\subsection{Time integration}\label{s:timeINT}
The fluid and structure subproblems are discretized in time using the generalized-$\alpha$ method~\cite{Hul93,JanWhi99,Bazilevs08a}. 
Although a monolithic solution of the fluid, structure, and Lagrange multiplier unknowns is possible at each time step, we employ an iterative scheme~\cite{Kamensky15ch} for updating the multipliers, in which an unconstrained problem with a fixed multiplier field is solved one or more times in each time step.

The solution procedure is performed in two stages. 
In the first stage, the unknowns at the time step $n+1$ for the fluid, mesh displacement, shell structure, and cable structure are solved implicitly using a combination of quasi-direct and block-iterative coupling strategies~\cite{Tezduyar05a,BazilevsBook1}, while holding the Lagrange multiplier fixed at its current value from the previous time step $n$.
For each Newton iteration in this stage, the fluid subsystem is solved first, followed by solutions for the mesh-moving, shell, and cable subsystems. 
The most recent solutions from previous subsystems are used when solving each current subsystem. 
Residuals of all subsystems are converged to a specified tolerance before proceeding to the next step.

The second stage updates the fluid--shell interface Lagrange multiplier by adding the normal component of penalty forces obtained from the fluid and structure solutions of the first stage. 
In the discrete form, the Lagrange multiplier is represented by scalars stored at the quadrature points used for integrals over $\mathcal{S}^\text{sh}$, which are updated explicitly at time step  $n+1$  as
\begin{equation}
\lambda^{n+1} = \lambda^n + \beta^\text{sh}_{\text{NOR}} R^{n+\alpha_f}\text{ ,}
\end{equation}
in which $R^{n+\alpha_f}$ is the regularized no-penetration constraint residual,
\begin{equation}\label{eq:perturbed-constraint}
R^{n+\alpha_f} = \left(\mathbf{u}^{n+\alpha_f}_\text{f} - (\mathbf{u}^\text{sh}_\text{s})^{n+\alpha_f}\right)\cdot\mathbf{n}_\text{s}^\text{sh} - \frac{r}{\beta^\text{sh}_{\text{NOR}}}\lambda^{n+1}\text{ ,}
\end{equation}
and $n+\alpha_f$ denotes the intermediate time level of the generalized-$\alpha$ algorithm between steps $n$ and $n+1$. 
Stability and accuracy analyses of this dynamic augmented Lagrangian method are detailed in~\citet{Kamensky17gx} and \citet{Yu18ga}.

\section{Numerical results}\label{s:appl}
In this section, we apply the surface reconstruction pipeline of Section~\ref{s:MV_geo} and the immersogeometric FSI model of Section~\ref{s:simFk} to the human MV and its repair.
We first reconstruct in-vivo MV geometries from segmented 3D TEE data and evaluate the reconstruction accuracy against the segmented point cloud.
We then describe the simulation setup, establish the convergence of the structural, fluid, and subcell discretizations, and compare the pre- and post-operative hemodynamics of the prolapsed valve and its single- and double-clip repairs.

\subsection{Mitral valve geometric modeling}
To assess the robustness of the proposed surface reconstruction pipeline in Section~\ref{s:MV_geo}, we perform a surface generation directly from a segmented point cloud including missing data and noise, which are common in clinical imaging.
Then we construct a full human MV model by including functionally equivalent CT structures and CT origins, which are directly utilized in our IMGA FSI framework.

We evaluate the quality of a reconstructed in-vivo surface using a dimensionless error metric called the scaled Nearest Neighbor Distance (sNND)~\cite{moola2026valvefit}.
sNND at the point $\mathbf{p} \in \mathbf{P}$ is defined by
\begin{equation}
    \text{sNND} (\mathbf{p},\mathscr{S}) =\min_{\mathbf{s}\in\mathscr{S}}  \frac{ \| \mathbf{p} - \mathbf{s} \|}{\sqrt{\mathcal{A}\left(\mathscr{S} \right)}}\text{ ,}
\end{equation}
in which $\mathbf{s}$ is a point on the fitted surface and $\mathcal{A}\left(\mathscr{S} \right)$ is the surface area of $\mathscr{S}$, which allows for consistent comparisons across different anatomical sizes of patient-specific MV.

\subsubsection{Patient-specific mitral valve}
\begin{figure}[!t]\centering
    \raisebox{6pt}{\begin{subfigure}[t]{0.69\textwidth}\centering
    \includegraphics[width=\textwidth]{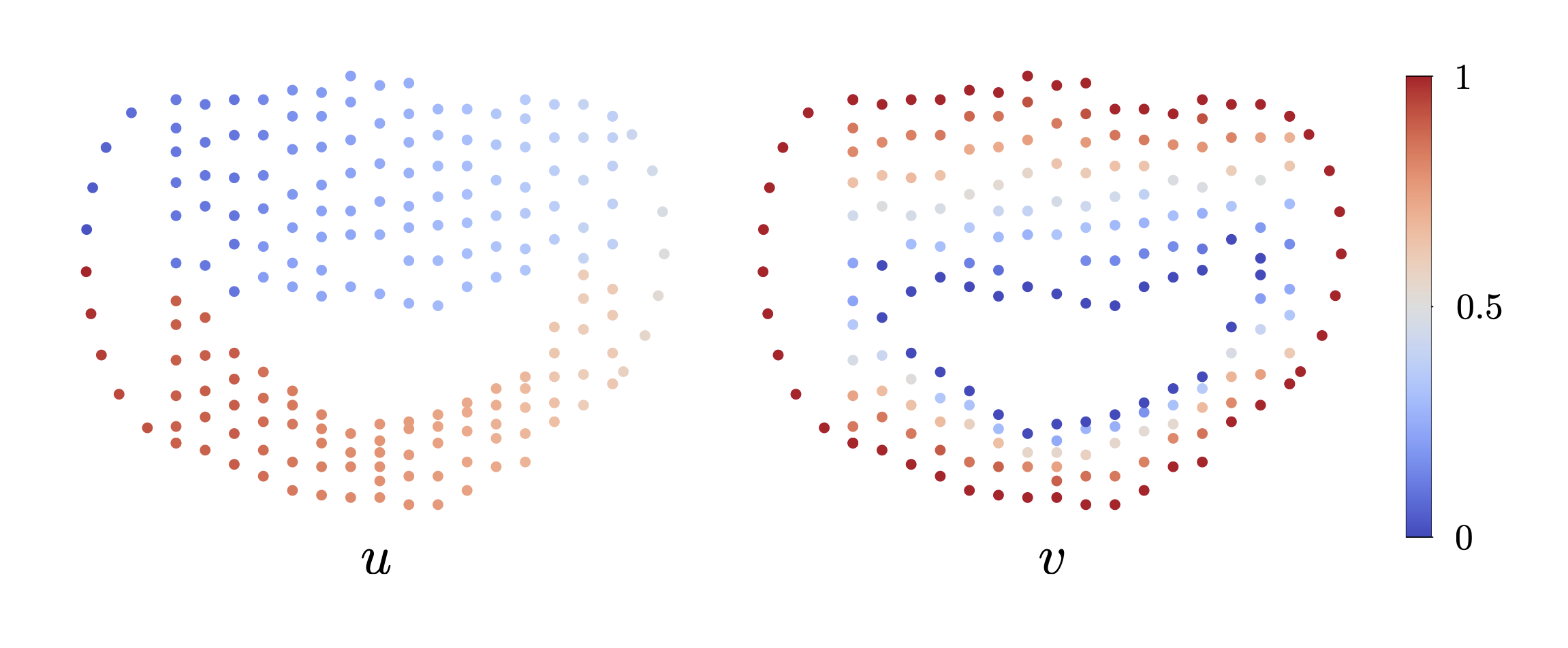}
    \caption{}\label{f:heat param}
    \end{subfigure}}
    \begin{subfigure}[t]{0.3\textwidth}\centering
    \includegraphics[width=\textwidth]{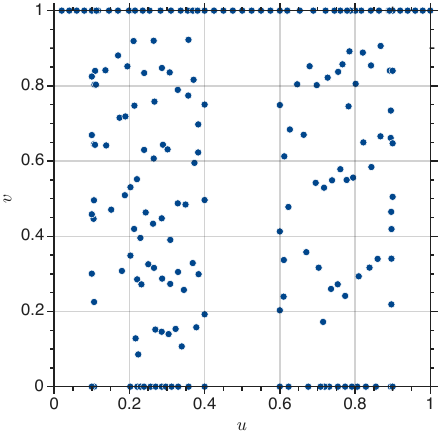}
    \caption{}\label{f:parametric domain}
    \end{subfigure}
  \caption{(a) Solutions of the Laplace equation on the top view of the segmented patient-specific MV point cloud in each parametric direction. (b) Parameterized segmented points in the parametric $(u,v)$-domain.}\label{f:parameterization}
\end{figure}

To demonstrate the reconstruction pipeline on a clinically relevant pathology, an MV point cloud is obtained directly from a patient-specific MV image featuring mitral regurgitation, as described in Section~\ref{s:segmentation}.
An example of a segmented patient-specific MV is shown in Figure~\ref{f:overview_figure}a.
Since our segmentation method produces two significant missing regions near the commissures, additional circumferential constraints are imposed on the separate anterior and posterior leaflet subproblems introduced in Section~\ref{s:parameterization}.
For the anterior leaflet subproblem, we apply additional constraints corresponding to $u = 0.1 \text{ and } 0.4$, in addition to $u = 0$ and $u = 0.5$ at the commissures.
Similarly, for the posterior leaflet, $u=0.6$ and $u=0.9$ are prescribed in addition to the commissural boundary values $ u =0.5$ and $ u =1$. 
Note that $u = 0$ and $u=1$ represent the same commissural seam as shown in Figure~\ref{f:overview_figure}b.
These particular constraints establish a robust parameterization that not only effectively defines the individual leaflet segments but also provides consistent parametric locations for chordal insertions.
We then parameterize the segmented point cloud based on the developed PDE-based parameterization in Section~\ref{s:parameterization}.
Circumferential and radial parameterized values are shown in Figure~\ref{f:parameterization}, and the parameterized values are used to obtain the matrix system for surface reconstruction as illustrated in Eq.~\eqref{eq:fitting equation} along with the periodic knot and open uniform knot vectors for the circumferential direction and radial direction, respectively.
The final matrix system is determined by a chosen number of control points in each direction and the degree of B-splines with corresponding knot vectors.
We use $n = 26$ and $p=3$ in the $u$-direction and $m = 4$ and $q=3$ in the $v$-direction.
Note that we add both the first- and second-order regularization terms to stabilize our direct surface reconstruction process.

\begin{figure}[!t]
    \centering
    \includegraphics[width=0.6\linewidth]{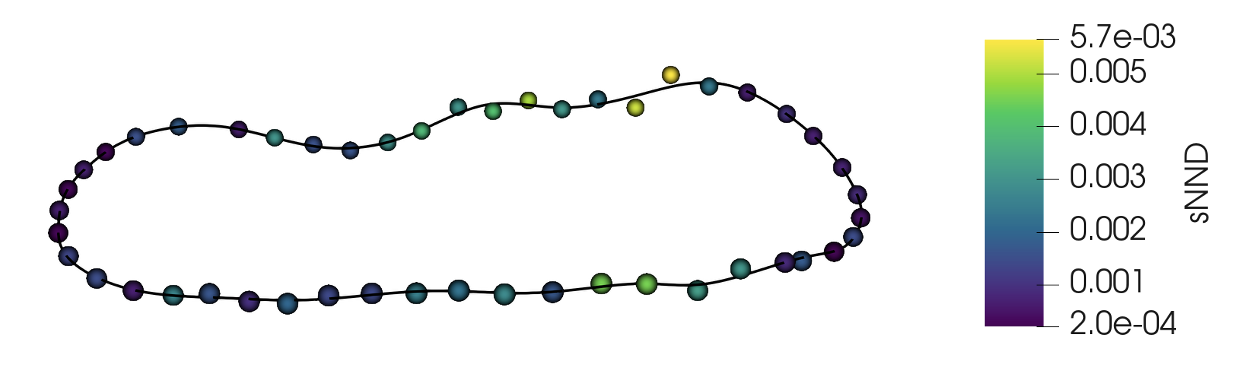}
    \caption{A reconstructed B-spline curve from the segmented annular points in Figure~\ref{f:overview_figure}a with colored pointwise sNND. We use $n=26$ and $p=3$. This serves as a constraint in the subsequent surface reconstruction. Note that the sNND reported here is defined as the shortest distance from each segmented point to the reconstructed curve normalized by the total arc length of the curve.
    }\label{f:annulus_fit}
\end{figure}

\begin{figure}[!t]
    \centering
    \begin{subfigure}[t]{\textwidth}
        \centering
        \includegraphics[width=0.39\linewidth]
            {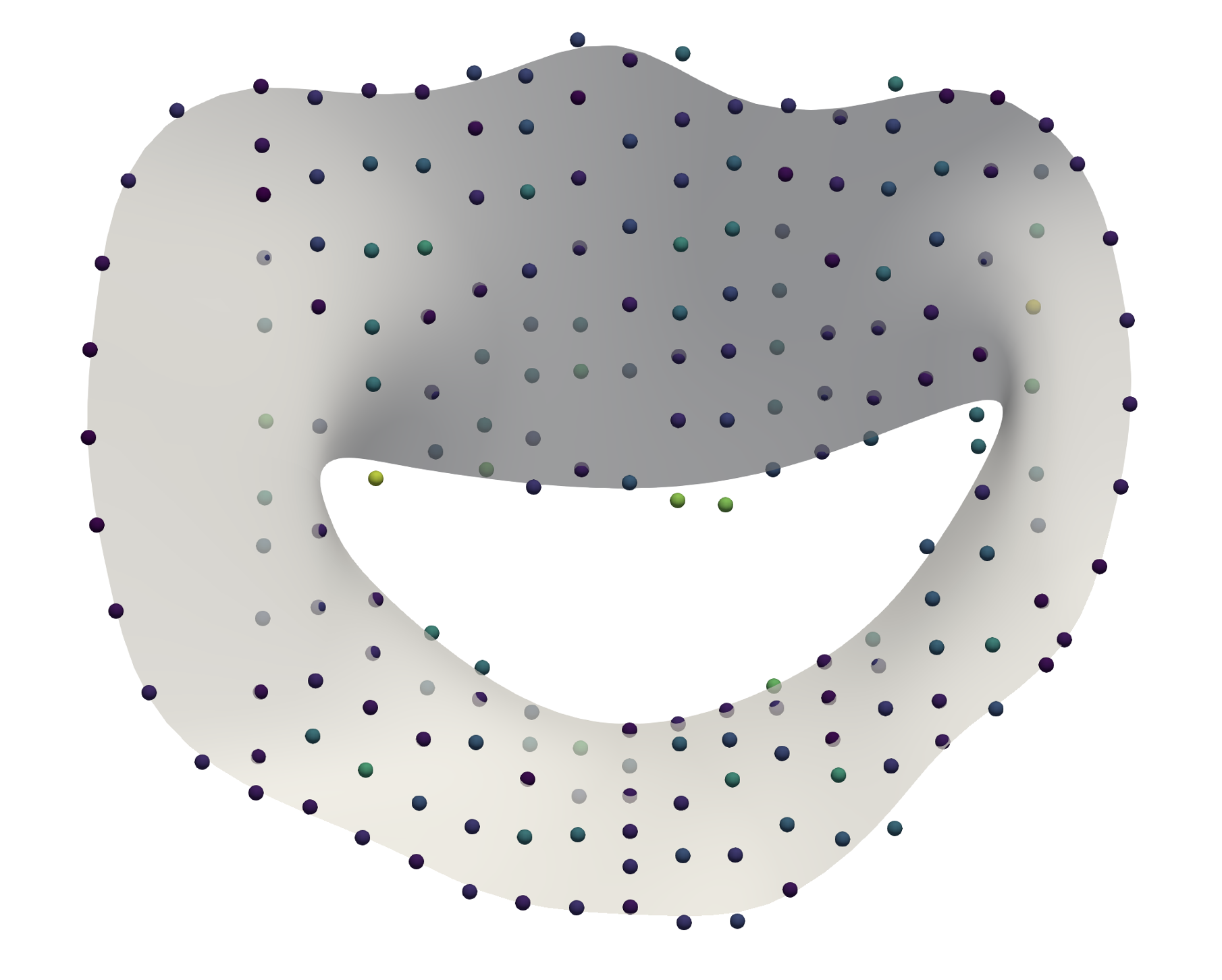}
        \raisebox{13pt}{%
            \includegraphics[width=0.59\linewidth]
                {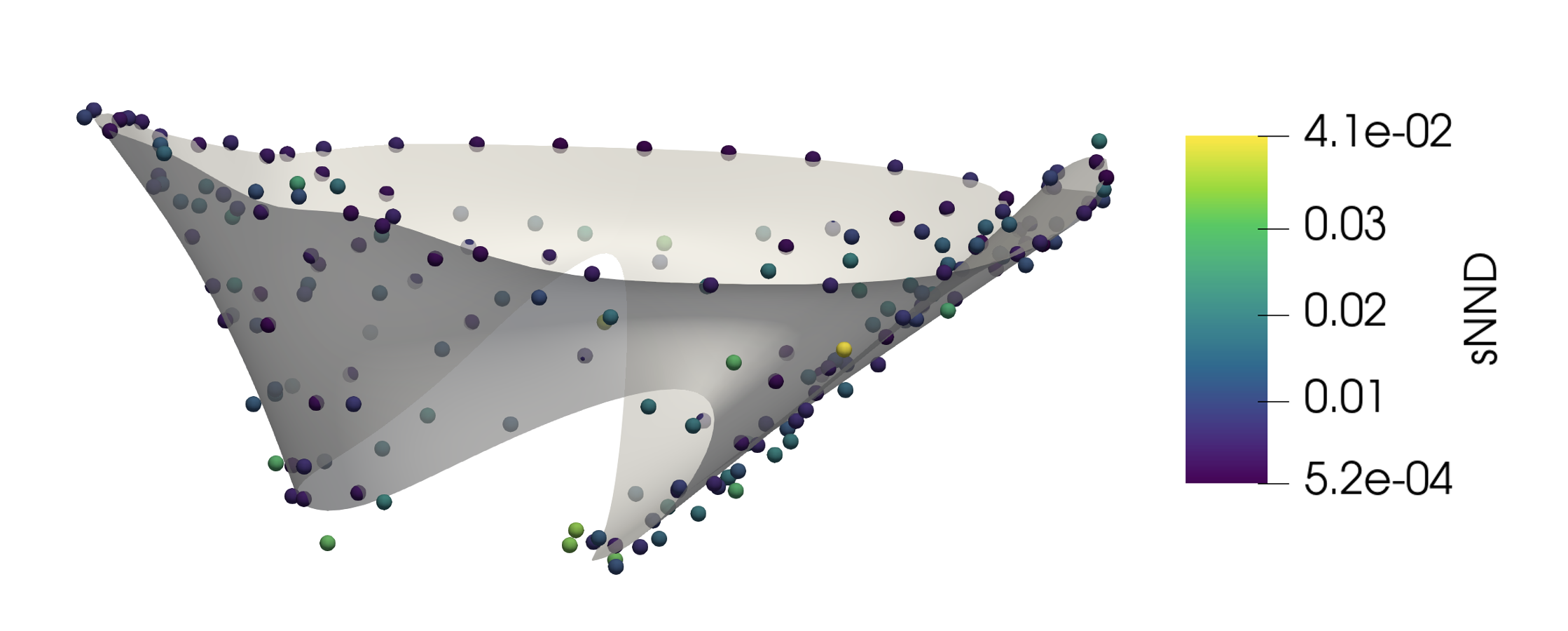}%
        }
    \end{subfigure}

    \caption{ A reconstructed B-spline surface from the segmented
    patient-specific MV point cloud in Figure~\ref{f:overview_figure}a
    with colored pointwise sNND. We use $n=26$ and $p=3$ for the
    circumferential direction and $m=4$ and $q=3$ for the radial
    direction.
    }
        \label{f:surface fitting}
\end{figure}

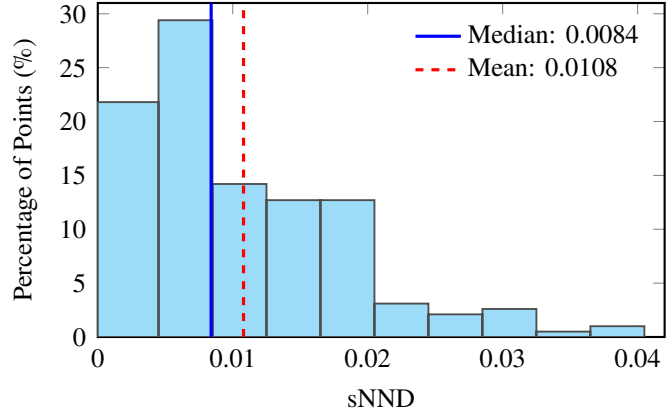
\begin{figure}[t!]
    \centering

    \begin{tikzpicture}
        \begin{axis}[
            width=0.55\textwidth,
            height=6cm,
            xlabel={sNND},
            ylabel={Percentage of Points (\%)},
            xmin=0, xmax=0.042,
            ymin=0, ymax=31,
            xtick={0,0.01,0.02,0.03,0.04},
            xticklabels={0,0.01,0.02,0.03,0.04},
            ytick={0,5,10,15,20,25,30},
            scaled x ticks=false,
            legend style={
                at={(0.98,0.97)},
                anchor=north east,
                draw=none,
                fill=none,
                font=\footnotesize
            },
            legend cell align={left},
            thick,
            label style={font=\footnotesize},
            tick label style={font=\footnotesize},
            legend style={font=\footnotesize},
        ]

        \addplot[
            ybar interval,
            fill=cyan!35,
            draw=black!70,
            line width=0.7pt,
            mark=none,
            forget plot
        ] coordinates {
            (0.0000,21.8)
            (0.0045,29.4)
            (0.0085,14.2)
            (0.0125,12.7)
            (0.0165,12.7)
            (0.0205,3.1)
            (0.0245,2.1)
            (0.0285,2.6)
            (0.0325,0.5)
            (0.0365,1.0)
            (0.0405,0.0)
        };

        \addplot[
            blue,
            very thick,
            solid,
            mark=none
        ] coordinates {
            (0.0084,0)
            (0.0084,31)
        };
        \addlegendentry{Median: 0.0084}

        \addplot[
            red,
            very thick,
            dashed,
            mark=none
        ] coordinates {
            (0.0108,0)
            (0.0108,31)
        };
        \addlegendentry{Mean: 0.0108}

        \end{axis}
    \end{tikzpicture}

    \caption{Histogram showing the percentage of points with respect to sNND reported in Figure~\ref{f:surface fitting}.}
    \label{f: sNND histogram}
\end{figure}

Figure~\ref{f:annulus_fit} shows an initial annular curve reconstruction for the subsequent surface reconstruction, and
Figure~\ref{f:surface fitting} shows a reconstructed human MV surface and corresponding pointwise sNND to the surface.
Approximately $70\%$ of the points have sNND values below the mean as shown in Figure~\ref{f: sNND histogram}. 
Table~\ref{t:fitting error} shows minimum, maximum, and mean sNND with several different choices of control points, demonstrating that the sNND error decreases as the number of control points increases.
We use $n = 26$ and $m = 4$ for the initial geometric modeling, which efficiently captures the overall in-vivo MV characteristics with a minimal number of control points.

\begin{table}[!t]\centering\small
\newcommand{\tabincell}[2]{\begin{tabular}{@{}#1@{}}#2\end{tabular}}
\caption{Minimum, maximum, and mean sNND errors between the reconstructed surface and the segmented point cloud with different choices of control points. We fix the number of control points in the $v$-direction and vary the number of control points in the $u$-direction.}
\label{t:fitting error}
    \begin{tabular}{|c|c|c|c|}\hline
      \tabincell{l}{$n$}   & \tabincell{l}{Minimum sNND}   & \tabincell{l}{Maximum sNND} & \tabincell{l}{Mean sNND}   \\ \hline
       $11$   &   $9.20 \times 10^{-4}$      &  $6.97 \times 10^{-2}$   &  $1.76 \times 10^{-2}$   \\
       $16$   &   $1.15 \times 10^{-3}$      &  $5.85 \times 10^{-2}$   &  $1.37 \times 10^{-2}$   \\
       $21$   &   $1.04 \times 10^{-3}$      &  $4.63 \times 10^{-2}$   &  $1.16 \times 10^{-2}$   \\
       $26$   &   $5.17 \times 10^{-4}$      &  $4.06 \times 10^{-2}$   &  $1.08 \times 10^{-2}$   \\
       $31$   &   $4.10 \times 10^{-4}$      &  $3.59 \times 10^{-2}$   &  $1.02 \times 10^{-2}$   \\
       $36$   &   $2.72 \times 10^{-4}$      &  $3.19 \times 10^{-2}$   &  $9.09 \times 10^{-3}$   \\\hline
    \end{tabular}
\end{table}

\subsubsection{Population-averaged mitral valve}
\label{s:MV template}
Patient-specific studies provide valuable insights into individualized valve functions and can capture clinically relevant anatomical details.
However, such a patient-specific model is limited in generalizability due to substantial inter-patient variability in valve morphology, and a single-subject study may not be representative. 
To investigate overall hemodynamic responses of the human MV after TEER, we construct a population-averaged MV surface for our numerical experiments rather than relying on a single patient-specific anatomy.

\begin{figure}[!t]
    \centering

    \begin{subfigure}[t]{0.65\textwidth}
        \centering
        \includegraphics[width=0.47\linewidth]
            {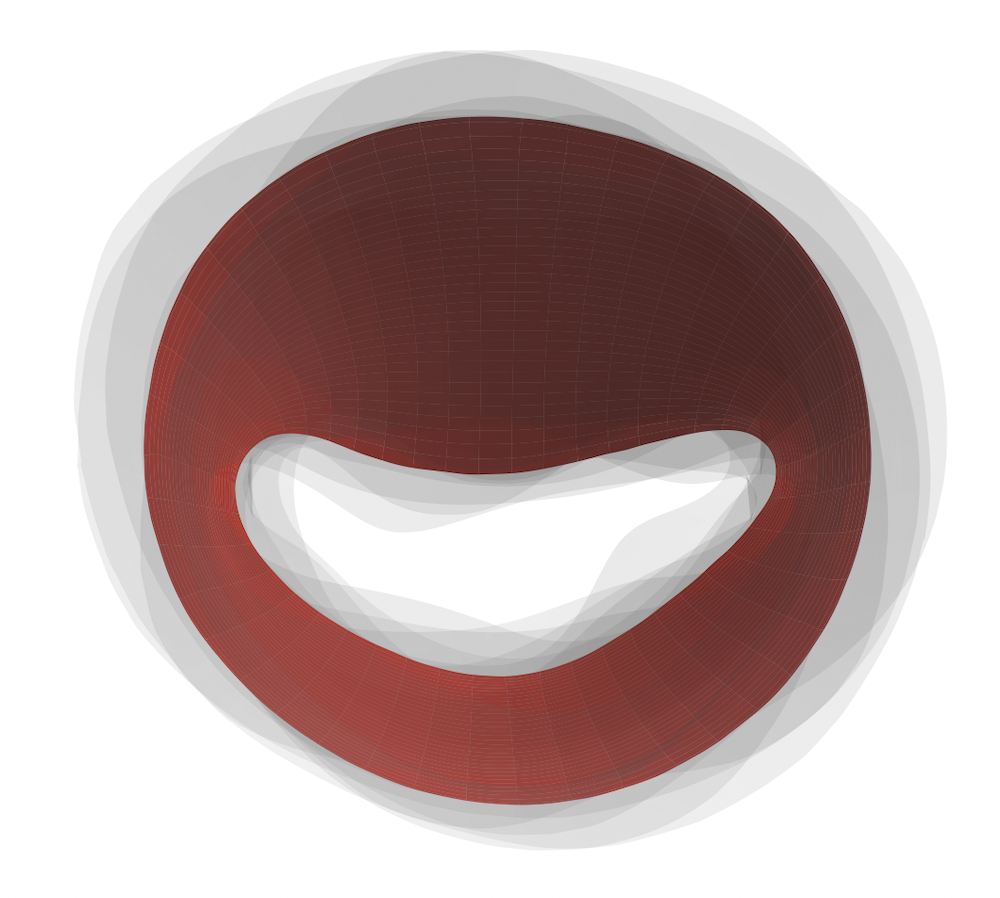}
        \raisebox{12pt}{%
            \includegraphics[width=0.47\linewidth]
                {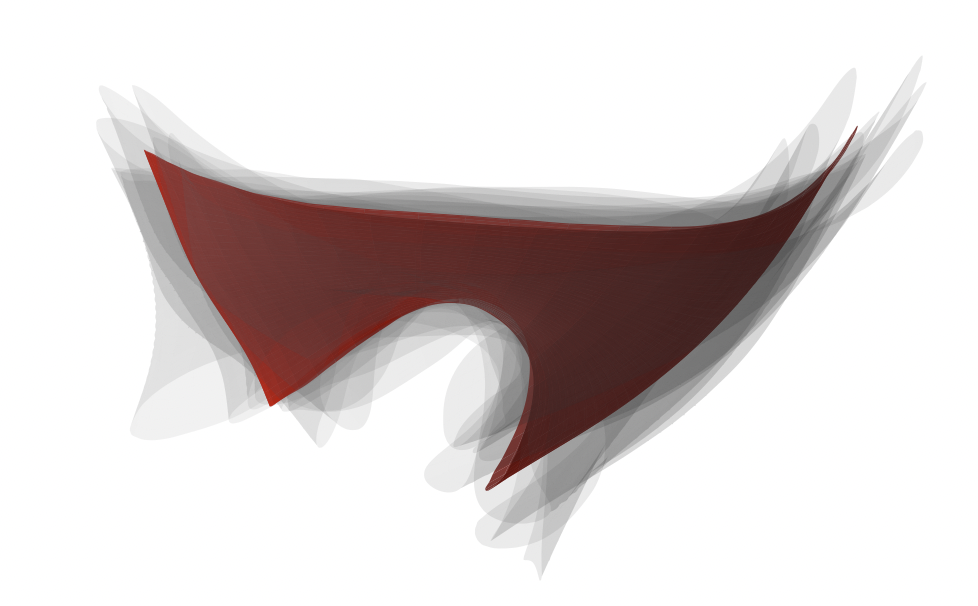}%
        }
        \caption{}
        \label{f:population averaged mv}
    \end{subfigure}
    \begin{subfigure}[t]{0.32\textwidth}
        \centering
        \includegraphics[width=\linewidth]
            {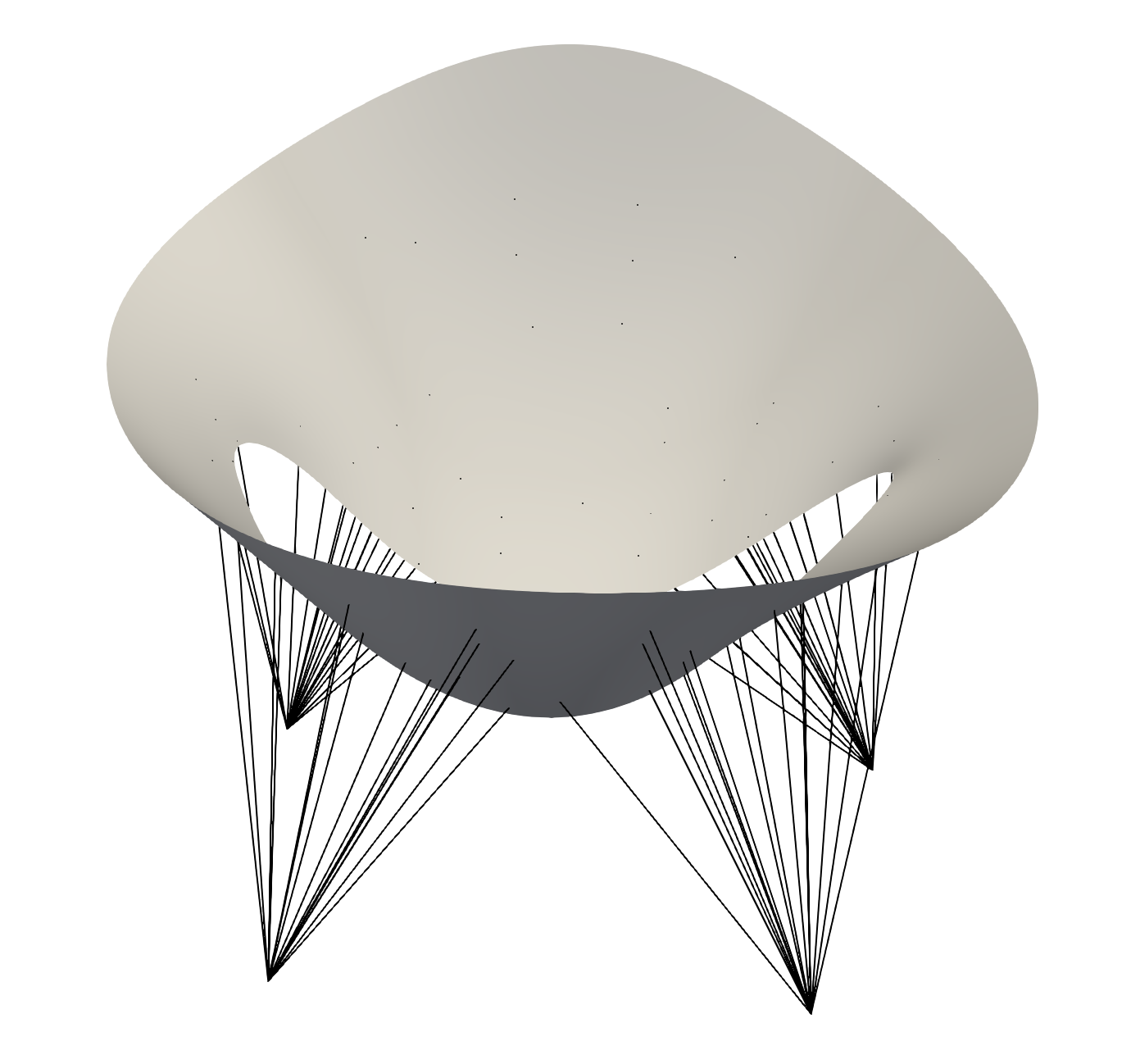}
        \caption{}
        \label{f:final_template_MV}
    \end{subfigure}

    \caption{%
    (a) Population-averaged MV surface (red) with
    10 patient-specific MVs (grey). All MV surfaces are centered
    at $(0,0,0)$ in a common coordinate system.
    (b) Population-averaged MV model representing the
    population-averaged MV leaflets along with 67 CTs based on
    the MVCT insertions and CT origins in Section~\ref{s:CT}.%
    }
    \label{f:population_averaged_mv_models}
\end{figure}

We obtain this surface by aligning and averaging ten patient-specific MV surfaces derived from an MR patient cohort through our proposed image-to-geometry pipeline.
Since our surface reconstruction process is performed entirely in parametric space, we can easily construct the population-averaged shape by taking the mean physical coordinates of corresponding parametric points across the valve.
Figure~\ref{f:population averaged mv} shows each patient-specific MV surface and the resulting population-averaged MV shape.
The obtained population-averaged surface captures the common anatomical features of the cohort, such as a distinct saddle shape and peaks near the middle of anterior and posterior leaflets, while smoothing out subject-specific irregularities.
It provides a stable and reproducible reference geometry for FSI simulations, allowing us to investigate hemodynamic responses that are broadly characteristic of the population rather than of a single patient.
A complete MV model is generated by combining the population-averaged leaflet surface with the population-averaged CT distribution map and the normalized CT origins introduced in Section~\ref{s:CT}.
Figure~\ref{f:final_template_MV} shows the resulting MV model, which consists of a single B-spline surface representing the mitral valve leaflets along with 67 CTs and 4 CT origins.

\subsection{Simulation setup}
\label{s:FSI setup}
We apply the proposed FSI methodology to model the population-averaged mitral valve within an idealized left ventricle undergoing prescribed motion. 
Fluid boundary conditions are pressure-driven with a physiologically realistic pressure waveform prescribed at the atrial boundary, while a three-element Windkessel and no-slip condition are dynamically applied to the arterial boundary based on hemodynamic quantities.
The pre-operative MV and its single- and double-clip repairs are simulated under identical boundary conditions.

\subsubsection{Pre-operative MV}
We use the population-averaged MV geometric model developed in Section~\ref{s:MV template} to simulate fully coupled fluid--structure interaction. 
The proposed geometric modeling pipeline efficiently captures the overall in-vivo MV characteristics with $n = 26$ and $m = 4$. 
For our valve simulations, the reconstructed geometry undergoes one level of refinement in the radial direction, yielding $m = 6$, which serves as the base mesh for subsequent simulation studies.
For the MV leaflet properties, the material parameters for the Lee--Sacks model are obtained by fitting to the MV leaflet model used in our previous work~\cite{Lee14a,liu2023computational}:
$c_0 = 10.11\,\text{kPa}$, $c_1 = 16.535\,\text{kPa}$, $c_2 = 1.0822$, $c_3 = 0.4316$, $w = 0.4231$, and the fiber direction $\mathbf{m}_0$ is set to the circumferential direction $\mathring{\mathbf{e}}_1$, that is, $(m_\text{c},m_\text{t},m_\text{n}) = (1,0,0)$.
The leaflet mass density is $1.0\,\text{g/cm}^3$, the thickness is $0.1\,\text{cm}$, and the in-vivo prestrain $\left(\lambda_\text{c},\lambda_\text{t}\right) = \left(1.3,0.934\right)$ is applied, adopting the calibrated values reported by~\citet{liu2023computational}.
We set an upper bound value of $\left(I_1^{\text{ub}}-3\right)= 0.8$ to improve the convergence of the nonlinear constitutive model and avoid constitutive responses in a non-physiological high-strain regime. 
For the CT properties, the Young's modulus is set to $E^{\text{ca}} = 2.45 \times10^8\,\text{dyn/cm}^2 = 24.5\,\text{MPa}$, the chordal density is $1.0\,\text{g/cm}^3$, and the chordal radius is $0.0225\,\text{cm}$~\cite{MILLINGTON-SANDERS1998,khalighi2017mitral}.
To replicate moderate-to-severe MR, we assign weakened chordal properties to the CTs in the P2 region indicated in Figure~\ref{f:MVCT}.
In particular, we set $E^{\text{ca}} = 1.45 \times10^7\,\text{dyn/cm}^2 = 1.45\,\text{MPa}$, which induces posterior leaflet prolapse characteristic of primary MR.
In addition, a mass-proportional damping coefficient of $c^\text{ca}=3.0\times10^{3}\,\mathrm{s}^{-1}$ is applied to the cable formulation to stabilize the chordal dynamics in the absence of direct fluid--cable coupling.
For the contact formulation, we set $r_\text{max} = 0.04\,\text{cm}$, $r_\text{min} = 0.008\,\text{cm}$, $\eta = 0.7$, $p_\text{c} =4$, and $k_\text{c} = 10^5\,\text{dyn/cm}^2$.

\subsubsection{Left heart chambers}
\begin{figure}[!t]\centering
    \centering
    \includegraphics[width=\textwidth]{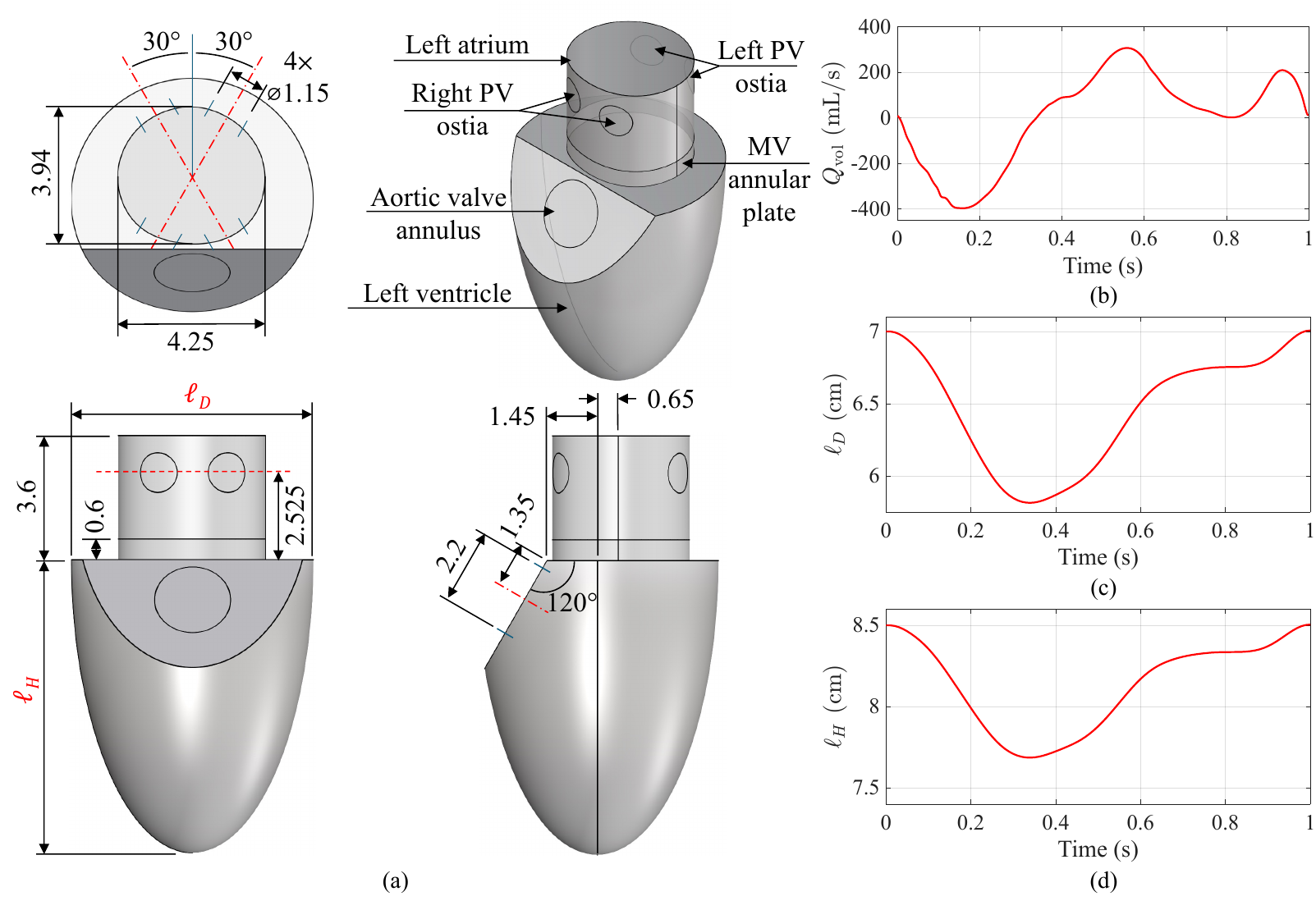}
    \caption{(a) Idealized geometry of the left atrium and ventricle. Dimensions ($\ell_D$, $\ell_H$) are parameterized to enable prescribed LV wall motion during the cardiac cycle. All units are in cm. (b) Prescribed volumetric flowrate ($Q_{\text{vol}}$). (c) Diameter and (d) height obtained by solving Eq.~\eqref{eqn:lv_ode}.}
    \label{fig:left_heart_geom}
\end{figure}
For FSI modeling of the mitral valve, we utilize an idealized computational fluid domain to represent the left heart chambers with a physiologically representative shape and scale. 
The fluid domain consists of the left atrium (LA) and the left ventricle (LV), connected through the MV annular plane, including specific openings for the pulmonary vein (PV) ostia and the aortic annulus as shown in Figure~\ref{fig:left_heart_geom}a. 
The LA is modeled as a vertical cylinder with a basal profile conforming to the MV annulus. 
The four PV ostia are located on the superior segment of the lateral face of the cylinder, symmetrically arranged in pairs to represent pulmonary venous inflow. 
The LV is represented by a truncated prolate ellipsoid that connects to the LA through the MV annular plane. 
The aortic annulus is included via a circular orifice on a planar cutout on the anterior basal region of the ventricle, which provides an outflow path during systolic ejection.

Since our primary focus is the FSI response of the MV leaflets, we use a prescribed ventricular wall motion to minimize the computational complexity of full structural dynamics of the left heart chambers while maintaining physiological realism. 
The wall motion is based on a simplified elastic membrane model that enforces a relationship between the LV volume variation and a clinically informed volumetric flowrate profile, $Q_{\text{vol}}(t)$, following the methodology established in~\cite{baccani2002vorte, tagliabue2017comple, tagliabue2017fluid}. 
The LV wall displacement is governed by the time-varying equatorial diameter $\ell_D(t)$ and base-to-apex height $\ell_H(t)$, and these parameters are obtained by solving the following system of ordinary differential equations (ODEs), which couples them to the volumetric flowrate:
\begin{subequations}
\label{eqn:lv_ode}
\begin{align}
\frac{d\ell_D}{dt} &= \frac{6Q_{\text{vol}}}{\pi} \frac{8\ell_H^2 - \ell_D^2}{20\ell_D\ell_H^3 - 2\ell_H\ell_D^3} \text{ ,} \\
\frac{d\ell_H}{dt} &= \frac{\ell_H}{\ell_D} \frac{d\ell_D}{dt} \frac{4\ell_H^2}{8\ell_H^2 - \ell_D^2} \text{ .}
\end{align}
\end{subequations}
The initial values are $\ell_D(0) = 7.0\,\text{cm}$ and $\ell_H(0) = 8.5\,\text{cm}$, corresponding to the end-diastolic stage of the cardiac cycle.
These parameters are derived based on the idealized ellipsoid model~\cite{tagliabue2017comple} and scaled to fit the population-averaged MV model constructed in Section~\ref{s:MV template}.
In particular, we calculate the pre-operative volumetric flowrate  $Q_{\text{vol}}(t)$ as described in Figure~\ref{fig:left_heart_geom}b from~\citet{caballero2020comprehensive}. 
The system of ODEs is solved using an explicit Runge--Kutta method~\cite{Bogacki1989runge} to obtain the time-varying $\ell_D$ and $\ell_H$ as presented in Figures~\ref{fig:left_heart_geom}c and \ref{fig:left_heart_geom}d. 
To map these parameters to the physical domain, the displacement of the LV wall $\mathbf{y}_\text{wall}$ is determined using an affine transformation:
\begin{align}
    \mathbf{y}_\text{wall} = \mathbf{A} (\mathring{\mathbf{x}} - \mathbf{x}_\text{e}) \text{ ,}
\end{align}
in which
\begin{align}
    \mathbf{A} =\begin{bmatrix}
        \frac{\ell_D}{\ell_{D}(0)}-1 & 0 & 0 \\
        0 & \frac{\ell_D}{\ell_{D}(0)}-1 & 0 \\
        0 & 0 & \frac{\ell_H}{\ell_{H}(0)}-1
        \end{bmatrix} \text{ ,}
\end{align}
$\mathring{\mathbf{x}}$ is the reference configuration, and $\mathbf{x}_\text{e}$ is the center of the ellipsoid representing the LV. 
Moreover, the mesh velocity of the LV wall is given by
\begin{align}
    \hat{\mathbf{u}}^h_\text{wall} = \frac{d\mathbf{A}}{dt} (\mathring{\mathbf{x}} - \mathbf{x}_\text{e}) \text{ .}
\end{align}
For the internal fluid domain, the mesh velocity is obtained by solving a linear elastostatics problem subject to the prescribed LV wall displacement as a boundary condition (BC)~\cite{Tezduyar93a, Johnson94a, Stein01f, Stein03a, Bazilevs08a}. 
For computational simplicity, the left atrium, the pulmonary vein ostia, the mitral valve annular plane, and the aortic annulus are treated as fixed boundaries. 
The resulting motion is localized to the LV, providing a physiologically informed contraction and relaxation of the ventricle. 
In our simulations, the fluid density is set to $1.0\,\text{g/cm}^3$, and the dynamic viscosity is set to $0.03\,\text{g}/(\text{cm}\cdot\text{s})$.

\subsubsection{Fluid boundary conditions}
To replicate physiological flow conditions without explicitly modeling the upstream pulmonary vasculature or downstream aortic tract, we employ cardiac state-dependent boundary conditions at the inlets (i.e., PV ostia) and the outlet (i.e., left ventricular outflow tract (LVOT)).

\begin{figure}[!t]\centering
	\begin{subfigure}{0.48\textwidth}\centering
        \includegraphics[width=\textwidth]{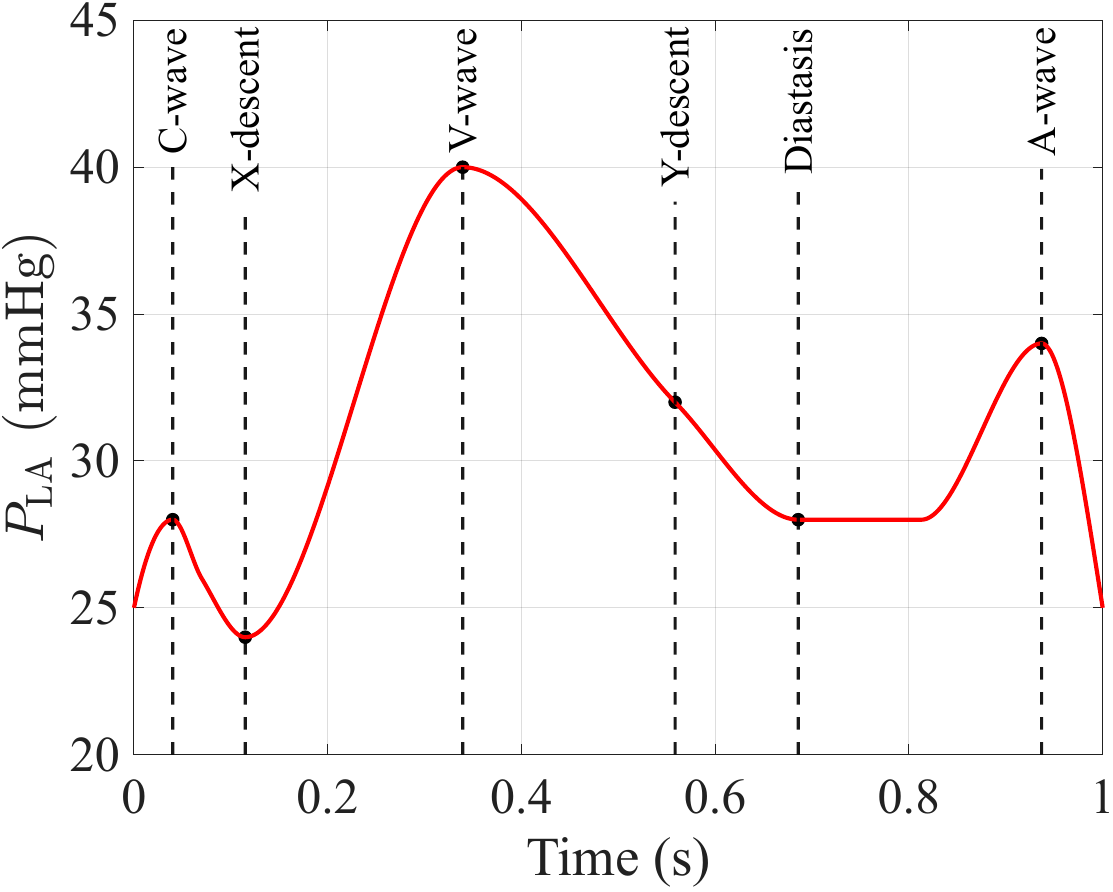}
        \caption{}
        \label{fig:fluidBC_atrial}
	\end{subfigure}
    \hspace{0.01\textwidth}
    \begin{subfigure}{0.48\textwidth}\centering
        \includegraphics[width=\textwidth]{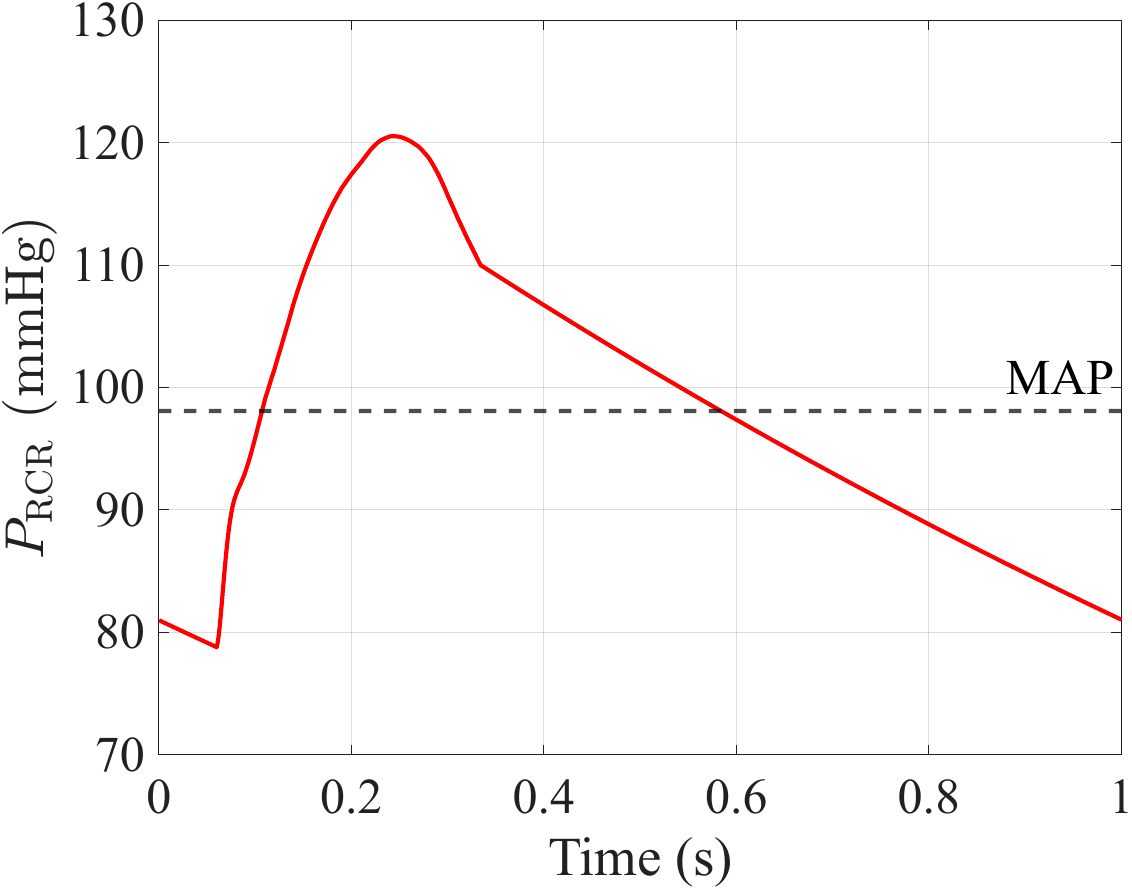}
        \caption{}
        \label{fig:fluidBC_RCR}
	\end{subfigure}
    \caption{(a) Prescribed atrial pressure driven by LV volumetric flowrate. (b) Aortic pressure generated by the Windkessel model with tuned parameters on a healthy prescribed flowrate.}
    \label{fig:fluidBC}
\end{figure}

A time-varying left atrial pressure waveform, shown in Figure~\ref{fig:fluidBC_atrial}, is constructed to be synchronized with the prescribed LV volumetric flowrate $Q_\text{vol}$ and applied at the PV ostia.
Mitral valve opening is identified as the instant at which $Q_\text{vol}$ first becomes positive following peak ejection, and the early E-wave and late A-wave are located over the subsequent diastolic interval (Figure~\ref{fig:left_heart_geom}b). 
These events anchor the timing of the V-wave, Y-descent, diastasis, and A-wave, while the C-wave and X-descent are placed within the isovolumetric contraction and early ejection intervals.
Characteristic pressure amplitudes are assigned to each event (V-wave: $40$~mmHg, A-wave: $34$~mmHg, X-descent minimum: $24$~mmHg). 
These values exceed the normal physiological range because, in severe mitral regurgitation, regurgitant flow into a non-compliant left atrium operating on the steep portion of its pressure--volume relation produces a disproportionately large pressure rise~\cite{pizzarello1984, kern1991vwave}.
A shape-preserving piecewise cubic Hermite interpolating polynomial is fitted through the event values to produce a continuous $P_\text{LA}(t)$ over the cardiac cycle.

Systemic arterial circulation is governed by a three-element Windkessel (RCR) lumped-parameter model coupled directly to the LVOT outlet face~\cite{Westerhof09, gao2017coupled}. 
The circuit consists of a proximal resistance ($R_\text{P}$), a capacitance ($C$), and a distal resistance ($R_\text{D}$), which collectively account for characteristic arterial compliance and total peripheral resistance. 
Parameter values ($R_\text{P} = 76.98\,\text{dyn}\cdot \text{s}/\text{cm}^5$, $C = 1.04 \times 10^{-3}\,\text{cm}^5/\text{dyn}$, and $R_\text{D} = 2111.47\,\text{dyn}\cdot\text{s}/\text{cm}^5$) are calibrated using a healthy mitral valve closure simulation to yield a physiological mean arterial pressure (MAP).
The obtained RCR parameters are then held fixed across all cases of prolapse and repair. 
Figure~\ref{fig:fluidBC_RCR} shows the resulting Windkessel pressure response computed offline for the prescribed healthy flow waveform. 
Note that this pressure waveform is not directly enforced as a boundary condition at the LVOT; rather, $P_\text{RCR}$ is dynamically computed through the implicit 0D--3D coupling and applied as a traction (Neumann) boundary condition on the outlet face~\cite{VFJT05, Esmaily13modul}. 

Additionally, to represent aortic valve dynamics throughout the cardiac cycle, the LVOT boundary condition switches between two states according to instantaneous pressure and flow criteria.
Systolic ejection begins when the left ventricular pressure exceeds the Windkessel pressure ($P_\text{LV} > P_\text{RCR}$).
During ejection, $P_\text{RCR}$ obtained from the 0D--3D coupling is applied as a traction boundary condition, representing the pressure load imposed by the systemic circulation.
Systolic ejection terminates when negative flow is detected at the LVOT ($Q_\text{LVOT} < 0 $), indicating aortic valve closure and the transition to diastole. 
Upon reaching this threshold, the outlet boundary condition transitions to a no-slip Dirichlet boundary condition ($Q_\text{LVOT} = 0$). 
During this phase, the Windkessel model continues to be solved in an isolated state with zero inflow, allowing $P_\text{RCR}$ to decay exponentially until $P_\text{LV} > P_\text{RCR}$ initiates the next systolic ejection phase.

\subsubsection{Repair modeling}
To construct the clipped leaflet regions and the fixed correspondence $\Phi_{\mathrm{clip}}$ introduced in Section~\ref{s:penalty_TEER}, we generate a pair of B-spline clip patches from the reconstructed leaflet geometry. 
Specifically, the clip patches are constructed from a surface interpolated from the structured leaflet control net over a two-dimensional parametric domain $(u, v) \in [0,1]^2$. 
These obtained patches serve only as geometric templates for identifying and pairing the clipped leaflet regions and introduce no additional structural degrees of freedom.

\begin{figure}[!t]\centering
	\begin{subfigure}{0.485\textwidth}\centering
        \includegraphics[width=.47\textwidth,valign=c]{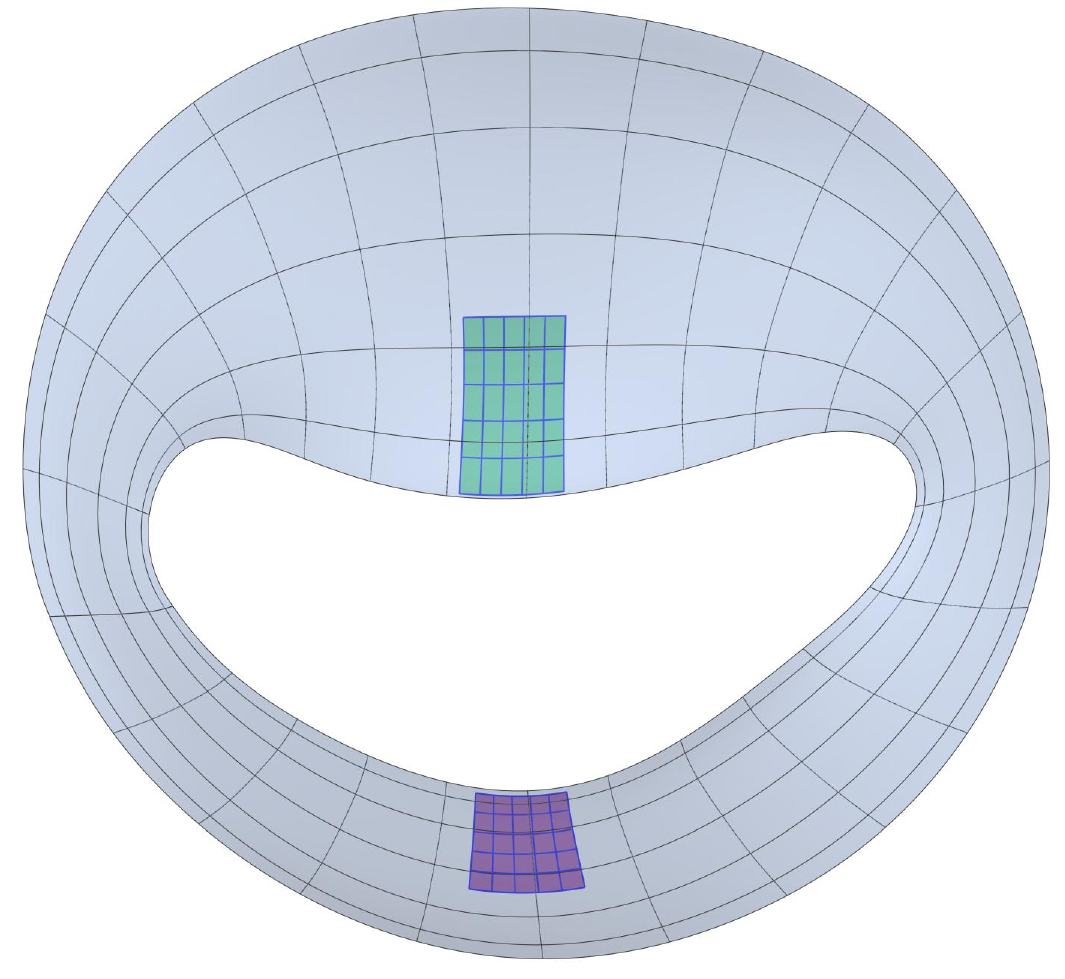}
        \includegraphics[width=.51\textwidth,valign=c]{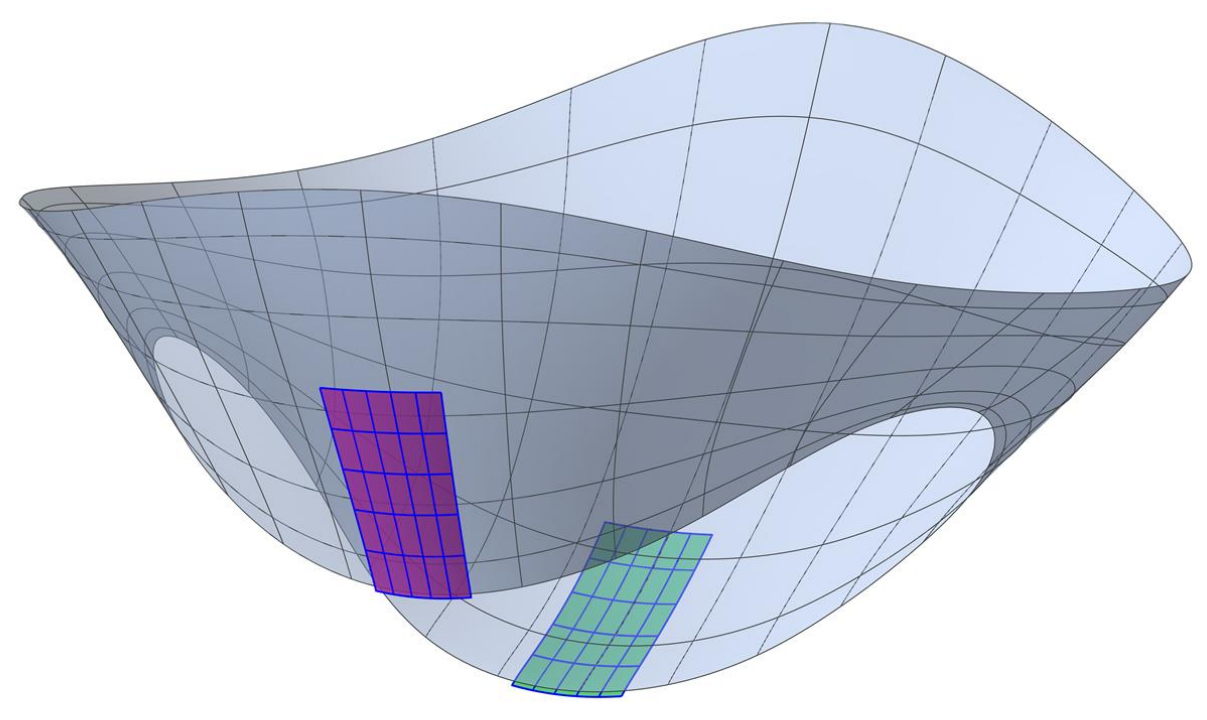}
        \caption{}
        \label{f:clip_patches_single}
	\end{subfigure}
    \hspace{0.01\textwidth}
    \begin{subfigure}{0.485\textwidth}\centering
        \includegraphics[width=.47\textwidth,valign=c]{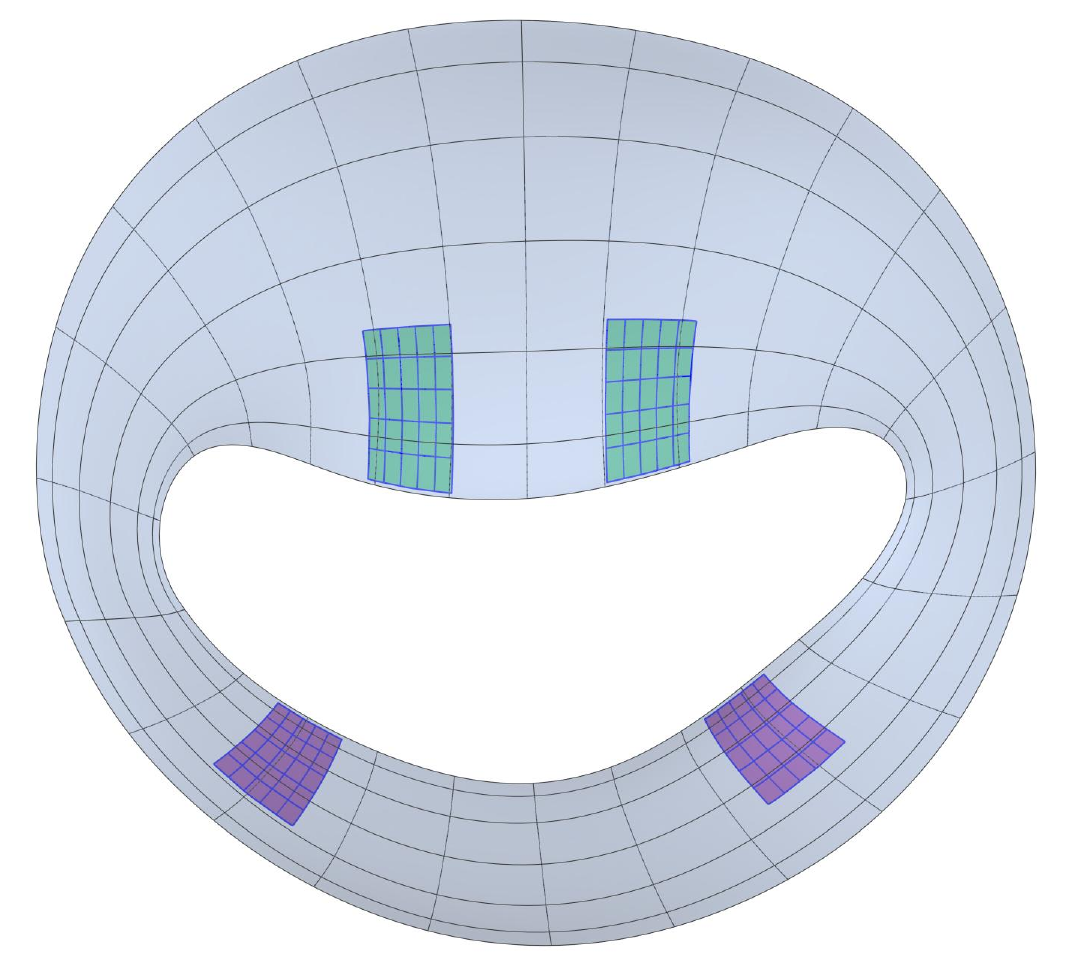}
        \includegraphics[width=.51\textwidth,valign=c]{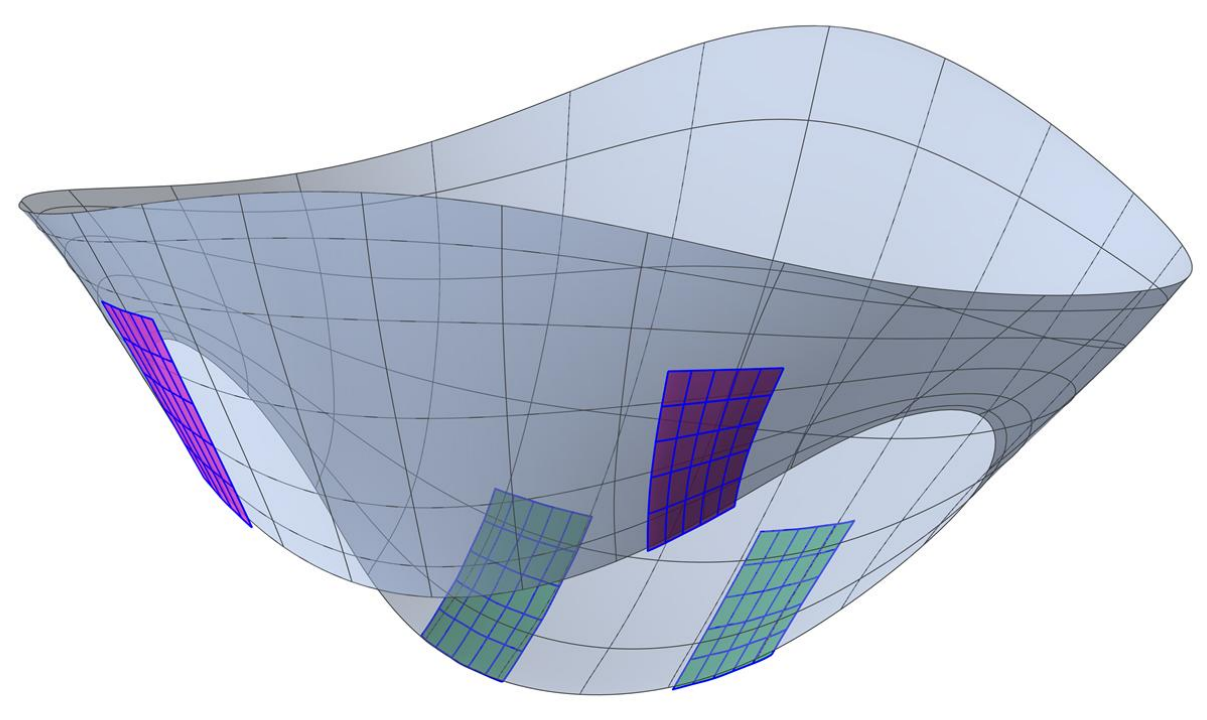}
        \caption{}
        \label{f:clip_patches_double}
	\end{subfigure}
    \caption{Clip patches on the population-averaged MV leaflet control mesh with $n = 26$ and $m = 6$. (a) A pair of clip patches representing a single central TEER device. (b) Two pairs of clip patches representing two TEER devices.}
    \label{f:clip_patches}
\end{figure}

The first clip patch is specified by two corner points, $(u_1, v_1)$ and $(u_2, v_2)$, and its bounding box $[u_{\min}, u_{\max}] \times [v_{\min}, v_{\max}]$ defines the parametric extent of the patch. 
The physical dimensions of this patch, $L_u$ and $L_v$, are evaluated by sampling the surface and computing discrete arc lengths along the $u$- and $v$-directions, respectively. 
The second patch is then placed at a specified center point $(u_\text{center}, v_\text{center})$ on the parametric domain with its physical size constrained to match that of the first patch. 
Since the surface is geometrically non-uniform, a direct parametric copy of the first patch cannot reproduce the same physical dimensions at a different location. 
We therefore employ an iterative solver that adjusts the parametric widths $\Delta u$ and $\Delta v$ of the second patch independently at each iteration by scaling them with the ratio of the target physical length to the current physical length.
This process is iterated until convergence,  yielding a second patch whose physical dimensions match those of the first to within a prescribed tolerance. 
Both patches are subsequently exported as B-spline surface patches with the same degree as the leaflet surface, thereby defining two separate clipping patches for TEER in our IMGA model. 
Figure~\ref{f:clip_patches} shows the resulting pair of clip patches on the population-averaged MV leaflet surface.

To ensure physiological realism in our repair simulations, the generated clip patches are constructed to match the $4\, \text{mm} \times 9\, \text{mm}$ leaflet-grasping area of a MitraClip device (Abbott)~\cite{Chakravarty20teer}.
Each obtained clip patch is discretized into a $3 \times 3$ element mesh using cubic B-spline basis functions, and a Gauss point distribution is constructed per element.
An explicit element-to-element and Gauss-point-to-Gauss-point mapping between the two patches is established to account for the periodic topology of the valve, preventing spurious twist in the coupling. 
Each clip-patch Gauss point is projected onto its closest point on the corresponding leaflet surface, and the resulting leaflet material points are paired according to this prescribed correspondence, which provides the discrete fixed correspondence for TEER.
Then the TEER penalty force defined in Eq.~\eqref{eq:TEER_penalty} is assembled directly into the degrees of freedom of the underlying leaflets, which ensures that the clip patches themselves do not introduce any additional structural degrees of freedom.
Additionally, the prescribed target gap prevents complete collapse of the paired leaflet regions, while the active contact penalty prevents interpenetration.
We set $g_\text{clip} = 0.05\,\text{cm}$ in our post-TEER simulations, which is larger than the contact cutoff distance $r_\text{max} = 0.04\,\text{cm}$.
This patch-wise coupling approach not only accurately reflects the actual size of the clipping device but also enables us to handle multiple clip configurations effectively, particularly in cases involving more than one clip device.

\subsection{Convergence study}
\label{s:convergence}
This section demonstrates the mesh convergence of our developed FSI model. 
We systematically isolate and verify individual modeling components before evaluating the fully coupled system.
An initial quasi-static structural convergence study of the valve determines the minimum mesh size required to accurately capture the structural deformation. 
Following this structural assessment, the fluid mesh convergence is performed along with the subcell quadrature across multiple background mesh resolutions.

\subsubsection{Structural convergence}
\label{s:structural_convergence}
\begin{figure}[!t]
    \centering
    \begin{subfigure}{0.3\textwidth}
        \centering
        \includegraphics[width=\textwidth]{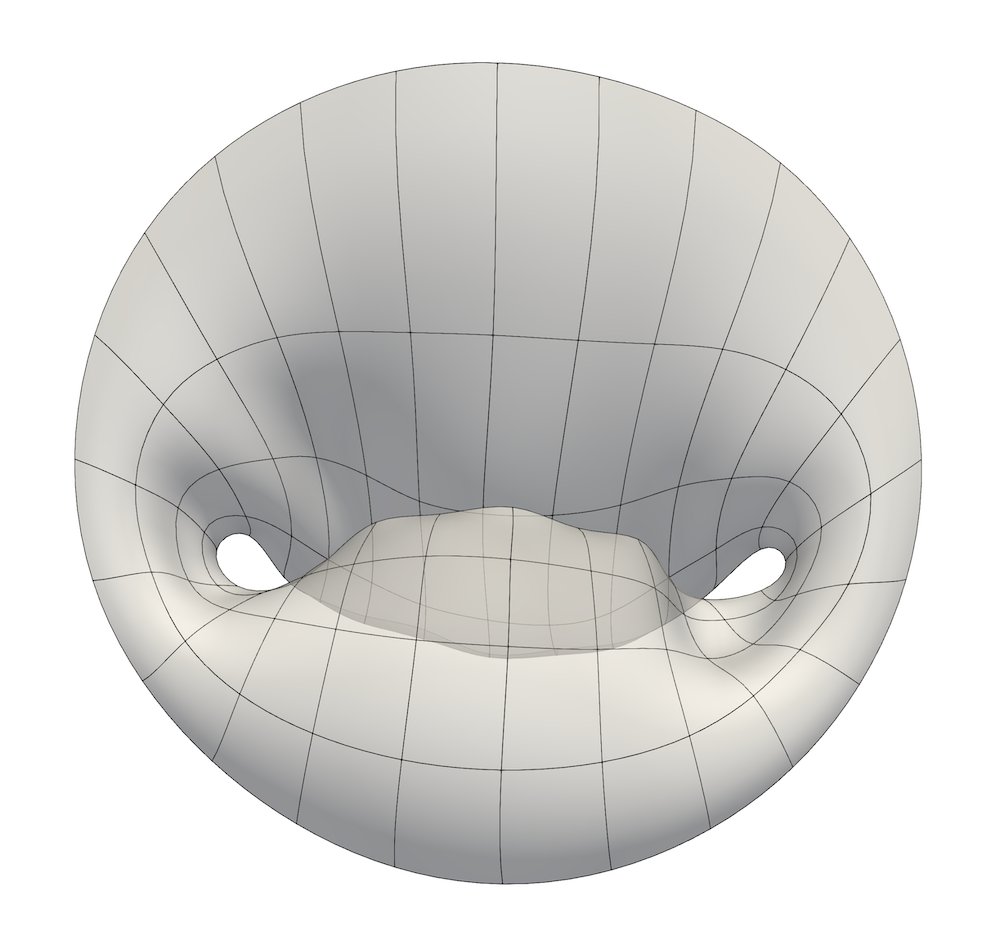}
        \caption{}
    \end{subfigure}
    \hfill 
    \begin{subfigure}{0.3\textwidth}
        \centering
        \includegraphics[width=\textwidth]{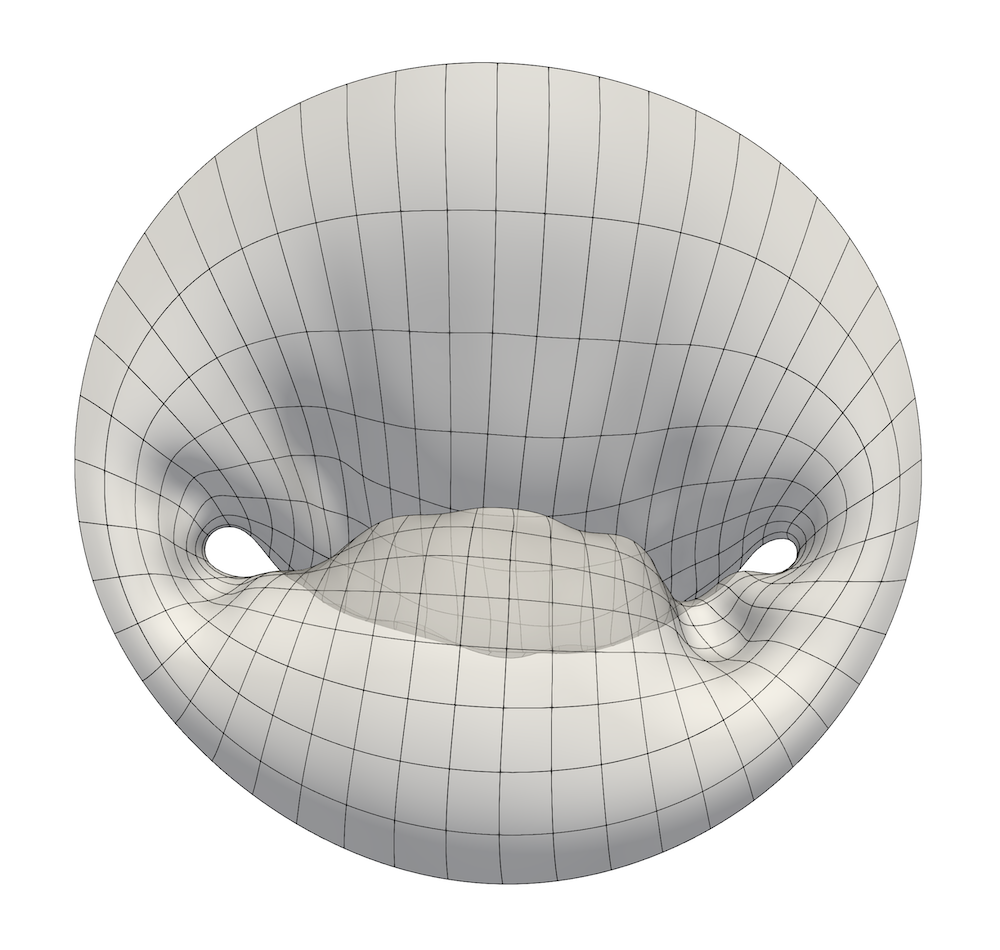}
        \caption{}
    \end{subfigure}
    \hfill
    \begin{subfigure}{0.3\textwidth}
        \centering
        \includegraphics[width=\textwidth]{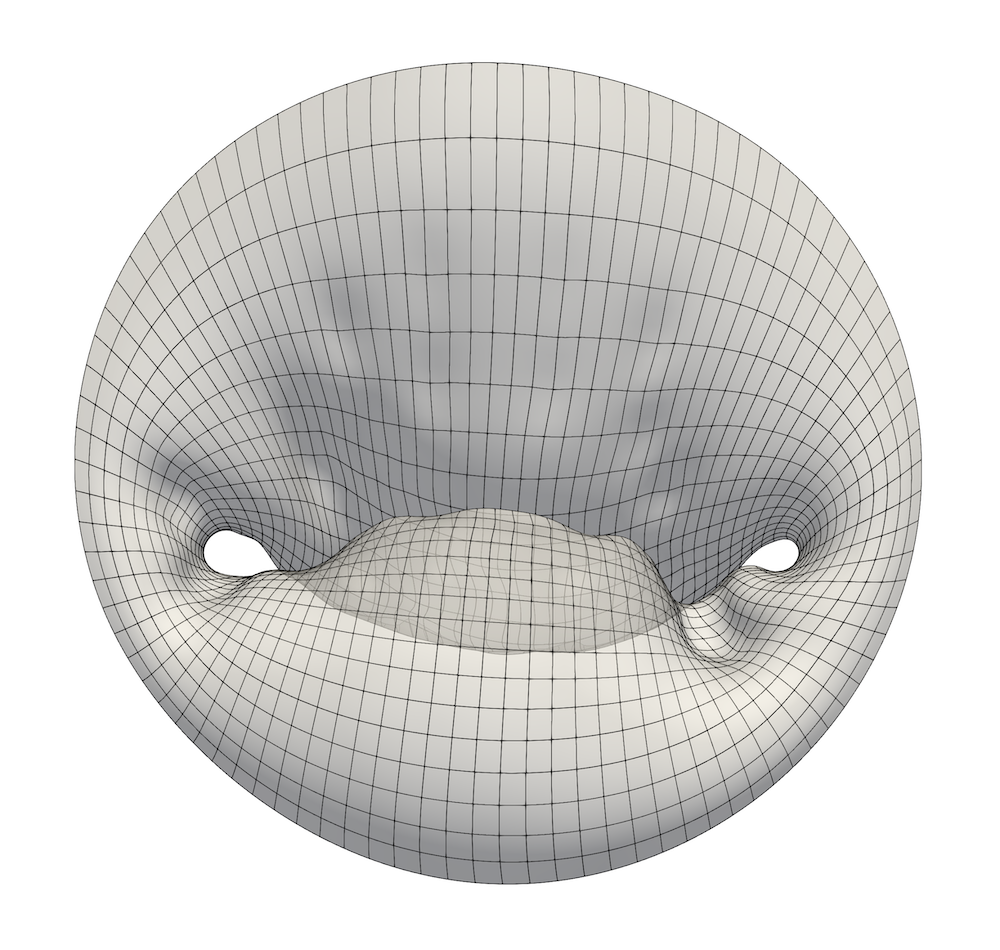}
        \caption{}
    \end{subfigure}

    \begin{subfigure}{\textwidth}
        \centering
        \includegraphics[width=\textwidth]{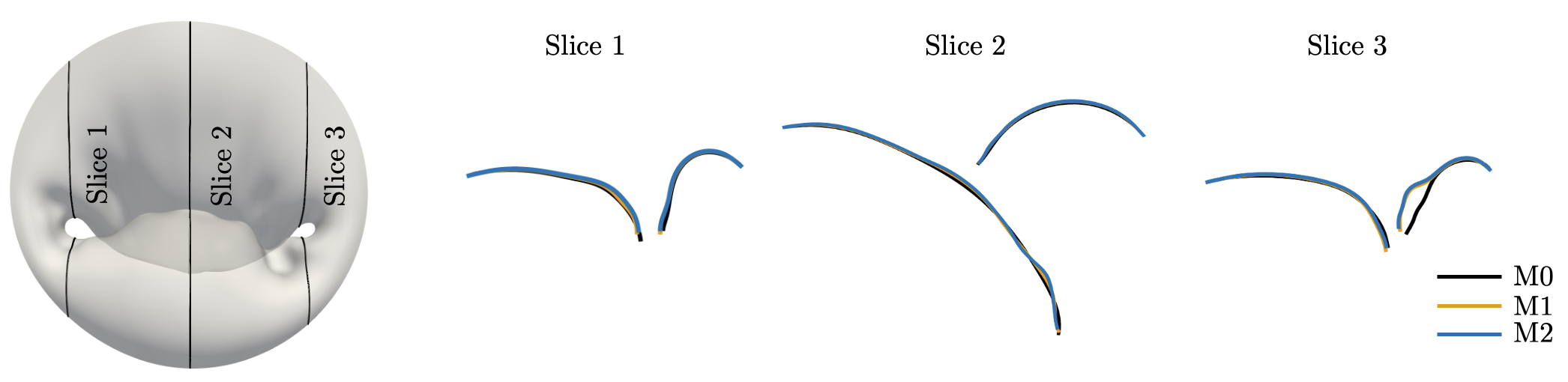}
        \caption{}\label{f:convergence_slices}
    \end{subfigure}
    
    \caption{Convergence of the fully closed configurations of the pre-operative MV under global $h$-refinement with a quasi-static loading of $100\, \text{mmHg}$, shown from the top view. (a) Mesh Level $0$ (M0) is represented by $24 \times4$ B-spline elements, (b) Mesh Level $1$ (M1) by $48\times8$ elements, and (c) Mesh Level $2$ (M2) by $96\times16$ elements. (d) Vertical slices for the closed configuration of the pre-operative MV under global $h$-refinement. The slices are extracted from the vertical planes indicated on the leaflet surface.}
    \label{fig:structural_convergence}
\end{figure}

For complex valvular geometries that undergo severe contact and localized buckling, it is challenging to achieve a consistent converged solution under mesh refinement. 
To demonstrate the structural mesh convergence of the pre-operative MV model within IMGA, we perform a structural analysis under quasi-static loading. 
A uniform surface-normal pressure load of $100\, \text{mmHg}$ is applied to the ventricular surface of the leaflets. 
Figure~\ref{fig:structural_convergence} illustrates the deformed configurations from the top view across three successive levels of global $h$-refinement: Mesh Level 0 (M0) with $24 \times 4$ B-spline elements, Mesh Level 1 (M1) with $48 \times 8$ elements, and Mesh Level 2 (M2) with $96 \times 16$ elements. 
To further verify the structural convergence, vertical slices are extracted from the defined planes on the leaflet surfaces (see Figure~\ref{f:convergence_slices}). 
Since we model a pathology of moderate-to-severe MR, significant prolapse of the posterior leaflet without any effective coaptation is observed.
The cross-sectional profiles demonstrate that the structural deformation converges under global $h$-refinement, establishing a sufficiently refined structural baseline for the subsequent FSI simulations. 

\begin{figure}[!t]
    \centering
    \begin{subfigure}{0.3\textwidth}
        \centering
        \includegraphics[width=\textwidth]{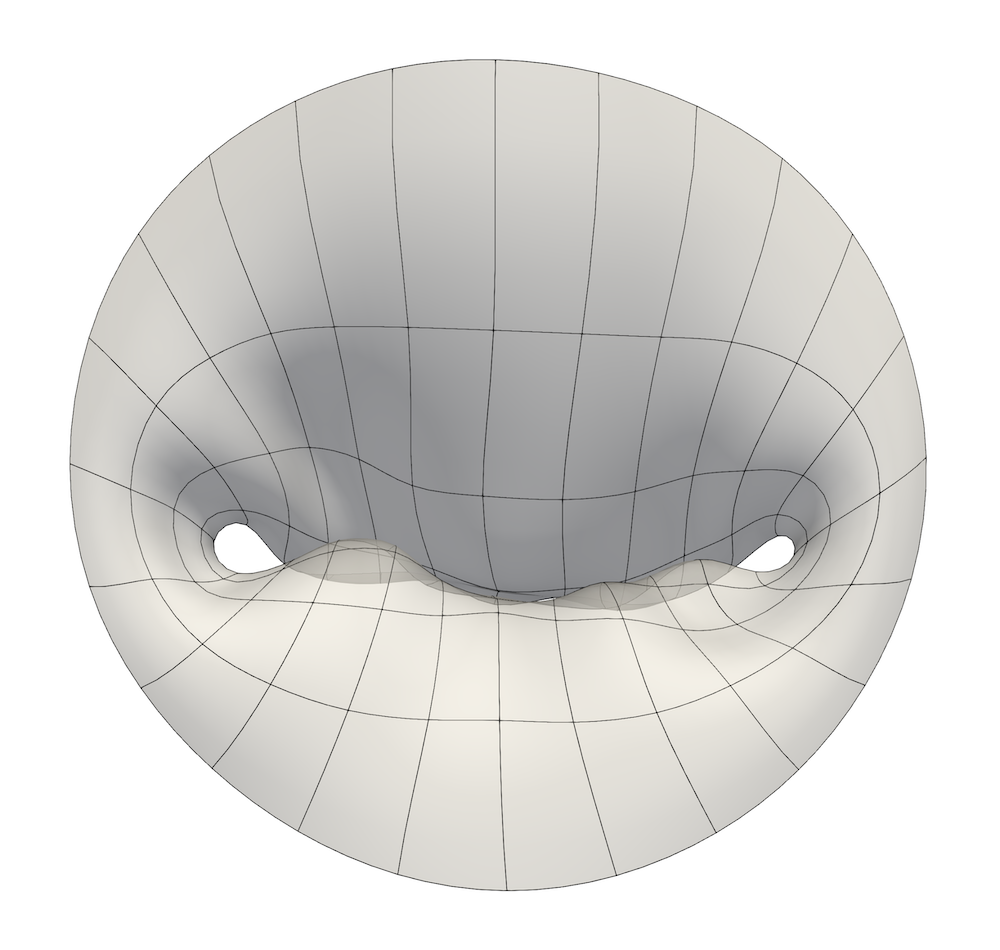}
        \caption{}
    \end{subfigure}
    \hfill 
    \begin{subfigure}{0.3\textwidth}
        \centering
        \includegraphics[width=\textwidth]{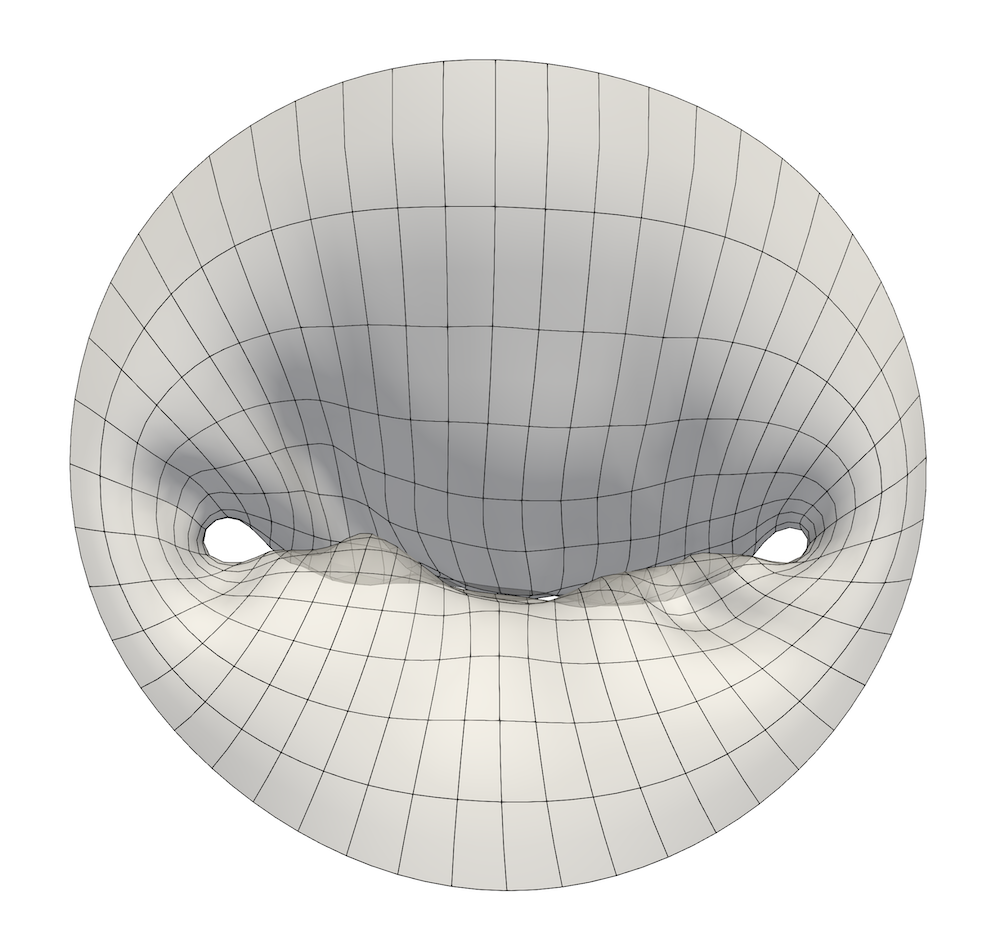}
        \caption{}
    \end{subfigure}
    \hfill
    \begin{subfigure}{0.3\textwidth}
        \centering
        \includegraphics[width=\textwidth]{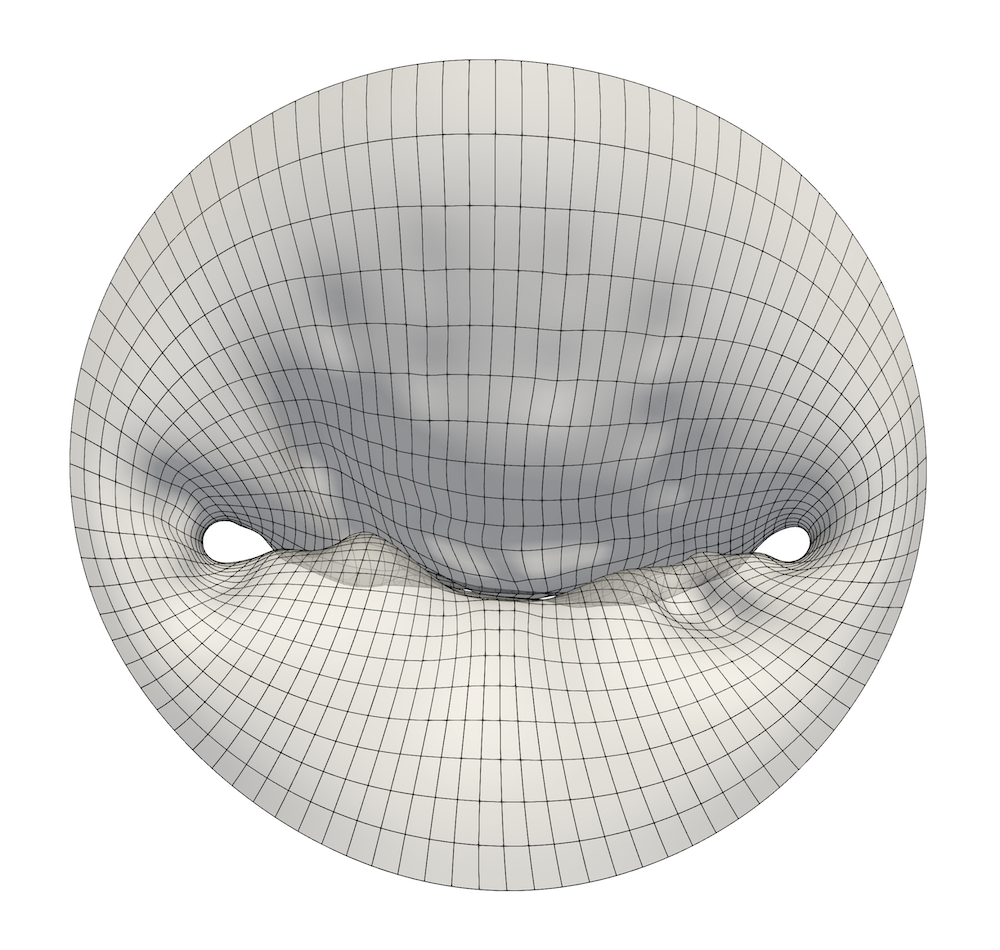}
        \caption{}
    \end{subfigure}

    \begin{subfigure}{\textwidth}
        \centering
        \includegraphics[width=0.65\textwidth]{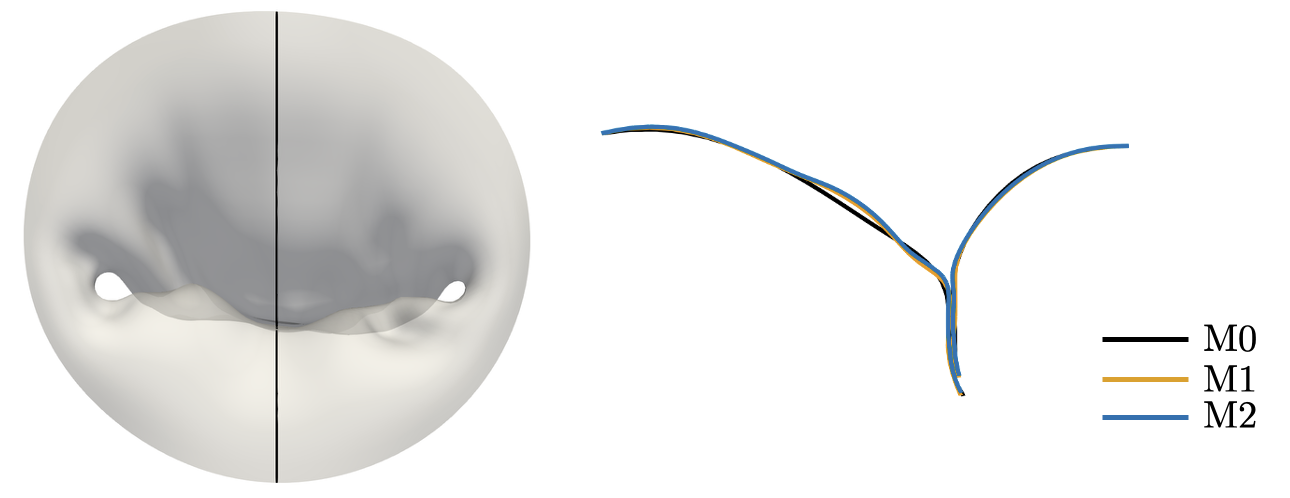}
        \caption{}\label{f:convergence_slices_post}
    \end{subfigure}
    
    \caption{Convergence of the fully closed configurations of the post-operative MV under global $h$-refinement with a quasi-static loading of $100\, \text{mmHg}$, shown from the top view. (a) M0, (b) M1, and (c) M2. (d) Vertical slices for the closed configuration of the post-operative MV under global $h$-refinement. The slices are extracted from the vertical plane indicated on the leaflet surface.}
    \label{fig:structural_convergence_post}
\end{figure}

In addition, we demonstrate the convergence of TEER modeling with the same mesh refinement strategy as in the pre-operative case.
In particular, we focus on a symmetric clipping scenario utilizing a single device as shown in Figure~\ref{f:clip_patches_single}.
Note that our parameterized pipeline for repair modeling allows us to effectively handle multiple clip configurations, including varying sizes and deployment positions.
Figure~\ref{fig:structural_convergence_post} shows the fully closed configurations of the post-operative MV under global $h$-refinement, and Figure~\ref{f:convergence_slices_post} further demonstrates the convergence of vertical slices of the repaired valve. 

\subsubsection{Subcell and fluid convergence}
Accurately enforcing the augmented Lagrangian constraints at the fluid--shell interface requires a sufficient spatial density of coupling points.
While this would usually involve global $h$-refinement of the structural mesh, we instead employ the subcell quadrature approach introduced in Section~\ref{s:subcell_discretization}, allowing us to enforce the constraints without over-refining the valve mesh.
Validating the subcell formulation across multiple background fluid mesh resolutions determines the appropriate subcell level for each background fluid mesh. 
Ultimately, this validation demonstrates independent fluid convergence and ensures robust enforcement of the fluid--structure coupling constraints along the interface.

Following the structural convergence results in Section~\ref{s:structural_convergence}, the shell mesh remains fixed at the M1 resolution, whereas the fluid mesh varies across three refinement levels. 
Table~\ref{t:mesh_stats} shows the corresponding fluid mesh statistics.
Since the dynamic RCR boundary condition relies on resolved flowrates that vary slightly during mesh refinement, it would introduce spurious variations into the fluid convergence analysis.
To maintain consistent boundary conditions across the different fluid meshes, we perform a quasi-static steady-state simulation with a transvalvular pressure gradient of $100\,\text{mmHg}$, enforcing  $120\,\text{mmHg}$ at the LVOT and $20\,\text{mmHg}$ at the PV ostia. 
All valvular properties and pathological conditions, including the localized weakening of CTs to induce prolapse, remain consistent with those in Section~\ref{s:structural_convergence}.

\begin{table}[!t]
    \centering\small
        \caption{Fluid mesh metrics with three different levels: coarse mesh level 0 (IM0), medium mesh level 1 (IM1), and fine mesh level 2 (IM2).}
    \begin{tabular}{lcccccc}
        \toprule
        \multirow{3}{*}[-1.2ex]{Mesh} & & & \multicolumn{3}{c}{Boundary layer} & \multirow{3}{*}[-1.2ex]{Total number of elements} \\
        \cmidrule(lr){4-6}
        & Base element & Near-valve & Number & First layer & Growth & \\
        & size (cm) & element size (cm) & of layers & height (cm) & rate & \\
        \midrule
        IM0 & $0.64$ & $0.32$ & $1$ & $0.2$ & $1.3$ & $64{,}416$ \\
        IM1 & $0.32$ & $0.16$ & $2$ & $0.1$ & $1.3$ & $245{,}041$ \\
        IM2 & $0.16$ & $0.08$ & $3$ & $0.05$ & $1.3$ & $1{,}444{,}320$ \\
        \bottomrule
    \end{tabular}
    \label{t:mesh_stats}
\end{table}

\begin{figure}[t!]
    \centering
    \begin{subfigure}[b]{0.48\textwidth}
        \centering
        \begin{tikzpicture}
            \begin{axis}[
                width=\linewidth,
                height=7cm,
                xlabel={Subcell Level},
                ylabel={Flowrate (mL/s)},
                xmin=-0.2, xmax=3.5,
                ymin=190, ymax=300,
                xtick={0, 1, 2, 3},
                xmajorgrids=true,
                grid style={dashed, black},
                legend style={
                    at={(1.00,0.85)},
                    anchor=east,
                    draw=none,
                    font=\footnotesize
                },
                thick,                
                label style={font=\footnotesize},
                tick label style={font=\footnotesize},
                legend style={font=\footnotesize},
            ]
            
            \addplot[color=blue, very thick, mark=*] coordinates {
                (0, 249.5) 
                (1, 248.39) 
                (2, 248.67) 
                (3, 248.11)
            };
            \addlegendentry{IM0}
            
            \addplot[color=red, very thick, mark=square*] coordinates {
                (0, 228.5) 
                (1, 205.73) 
                (2, 206.94) 
                (3, 205.46)
            };
            \addlegendentry{IM1}
            
            \addplot[color=green!50!black, very thick, mark=triangle*, mark size=3.0pt] coordinates {
                (0, 280.83) 
                (1, 217.03) 
                (2, 204.37) 
                (3, 204.10)
            };
            \addlegendentry{IM2}
            
            \end{axis}
        \end{tikzpicture}
        \caption{}
        \label{fig:convergence_subcell_flow}
    \end{subfigure}
    \hfill
    \begin{subfigure}[b]{0.48\textwidth}
        \centering
        \begin{tikzpicture}
            \begin{axis}[
                width=\linewidth,
                height=7cm,
                xmode=log,
                xlabel={$h_{\text{shell}}/h_{\text{fluid}}$},
                ylabel={Converged Subcell Level},
                xmin=1, xmax=8,
                xtick={1.5,3,6},
                xticklabels={1.5,3,6},
                xminorticks=false,
                ymin=-0.5, ymax=2.5,
                ytick={0, 1, 2, 3},
                ymajorgrids=true,
                xmajorgrids=true,
                grid style={dashed, black},
                thick,
                label style={font=\footnotesize},
                tick label style={font=\footnotesize},
                legend style={font=\footnotesize},
            ]
            \addplot[color=black, very thick, mark=square, mark size=3.0pt] coordinates {
                (1.5467, 0)
                (3.0932, 1)
                (6.1867, 2)
            };
            
            \end{axis}
        \end{tikzpicture}
        \caption{}
        \label{fig:convergence_mesh_ratio}
    \end{subfigure}
    
    \caption{(a) Outflow rate convergence across subcell level for different fluid mesh refinements. (b) Subcell level required to achieve flowrate convergence versus the mesh size ratio $h_\text{shell}/h_\text{fluid}$.}
    \label{fig:overall_convergence}
\end{figure}
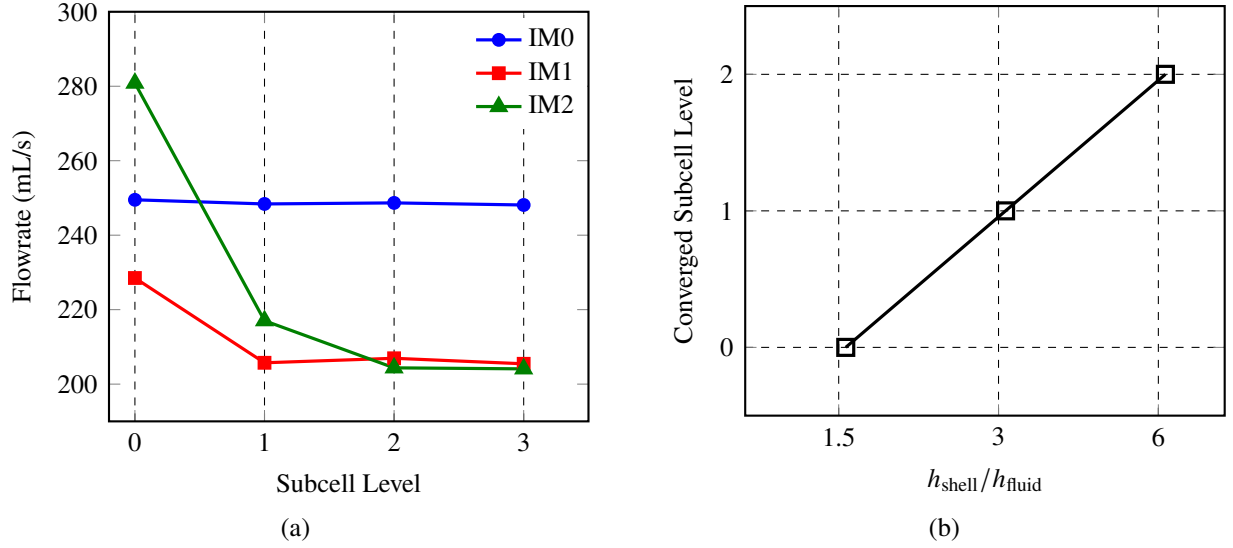

A total of twelve simulations are performed, spanning all combinations of four subcell levels (SC0, SC1, SC2, and SC3) and three fluid mesh levels (IM0, IM1, and IM2). 
Figure~\ref{fig:convergence_subcell_flow} shows the time-averaged volumetric flowrate through the atrial faces over 2 seconds. 
The flowrate obtained with the coarsest fluid mesh (IM0) is insensitive to subcell refinement, but it settles at a value well above the converged flowrates obtained with IM1 and IM2.
The required subcell level increases with fluid mesh refinement, with IM1 converging at SC1 and the finest fluid mesh (IM2) requiring SC2.
Figure~\ref{fig:convergence_mesh_ratio} illustrates the converged subcell level against the mesh size ratio $h_\text{shell}/h_\text{fluid}$, which establishes the necessary subcell density for a specific fluid--shell mesh pair.
These results demonstrate that converged hemodynamic quantities can be obtained under fluid mesh and subcell refinement while maintaining a single baseline structural resolution.
This allows us to use the coarsest converged valve mesh without refining it further or adding structural quadrature points to match the resolution of the background fluid mesh.

We further verify FSI convergence by utilizing the minimum converged subcell level for each grid: SC0 for IM0, SC1 for IM1, and SC2 for IM2 for both closure and opening. 
The closure case evaluates the geometric regurgitant orifice area (GROA) under a $100$~mmHg transvalvular pressure gradient, whereas the opening case evaluates the geometric orifice area (GOA) under a transvalvular pressure gradient of $4$~mmHg.
Simulations are run until the flow develops into a quasi-steady solution, after which we average the last 2 seconds of each simulation for both the fluid and the valve to obtain representative mean states.
We measure GROA and GOA by orienting the averaged valve model to see the maximum orifice opening.
Slices are taken along this view path and traced to find the location of greatest flow restriction, and the orifice area traced on this slice is taken as the GROA or GOA.
Table~\ref{t:metrics_comparison} shows that the key hemodynamic variables converge at the IM1 level, with negligible change from IM1 to IM2, and Table~\ref{t:GROA_convergence} shows the same behavior for the GROA and GOA under the closure and opening loads.
These results demonstrate convergence of both fluid and shell quantities and establish the structural, fluid, and subcell resolutions for pre- and post-operative FSI simulations.

\begin{table}[!t]
    \centering\small
    \caption{Mesh convergence of hemodynamics in the left heart.}
    \begin{tabular}{ccccc}
        \toprule
        Fluid mesh & Subcell level & $P_{\mathrm{LV}}$ & $P_{\mathrm{LA}}$ & $Q_{\mathrm{MV}}$ \\
        \midrule
        IM0     & SC0    & $118.66$~mmHg & $21.77$~mmHg  & $249.50$~mL/s  \\
        IM1     & SC1    & $118.99$~mmHg & $22.33$~mmHg  & $205.73$~mL/s  \\
        IM2     & SC2    & $118.99$~mmHg & $22.29$~mmHg  & $204.37$~mL/s  \\
        \bottomrule
    \end{tabular}
    \label{t:metrics_comparison}
\end{table}

\begin{table}[!t]
    \centering\small
    \caption{Mesh convergence of the geometric regurgitant orifice area (GROA) under a $100$~mmHg closing transvalvular pressure gradient and the geometric orifice area (GOA) under a $4$~mmHg opening transvalvular pressure gradient.}
    \begin{tabular}{cccc}
        \toprule
        Fluid mesh & Subcell level & GROA & GOA \\
        \midrule
        IM0     & SC0   & $1.88$~cm$^2$  & $5.03$~cm$^2$ \\
        IM1     & SC1   & $1.56$~cm$^2$  & $4.95$~cm$^2$ \\
        IM2     & SC2   & $1.56$~cm$^2$  & $4.95$~cm$^2$ \\
        \bottomrule
    \end{tabular}
    \label{t:GROA_convergence}
\end{table}

\subsection{FSI analysis}
In this section, we investigate the effects of TEER on left ventricular hemodynamics across three cases: the pre-operative prolapse, a central single-clip repair, and a double-clip repair with one clip placed on each side of the central position.
We first simulate a diseased mitral valve undergoing moderate-to-severe regurgitation over the cardiac cycle to establish a pathological baseline.
To evaluate the efficacy of the TEER procedure in reducing this MR, we then analyze the post-intervention hemodynamic changes.

All FSI simulations are conducted using consistent boundary conditions as described in Section~\ref{s:FSI setup}.
Based on our preliminary convergence studies in Section~\ref{s:convergence}, we utilize the converged discretization corresponding to the IM1 fluid domain, M1 structural valve mesh, and SC1 subcell level.
Each case is simulated over five cardiac cycles to achieve a periodic steady state. 
It is important to note that the RCR parameters are calibrated to a healthy physiological baseline and held constant across all cases.
When quantifying mitral valve performance, we use the guidelines from the American Society of Echocardiography report by~\citet{zoghbi2017} and a consensus document from the Mitral Valve Academic Research Consortium by~\citet{stone2015}.

\subsubsection{Pre-operative FSI}

\begin{figure}[!t]\centering
	\begin{subfigure}[b]{0.48\textwidth}\centering
        \includegraphics[width=\textwidth]{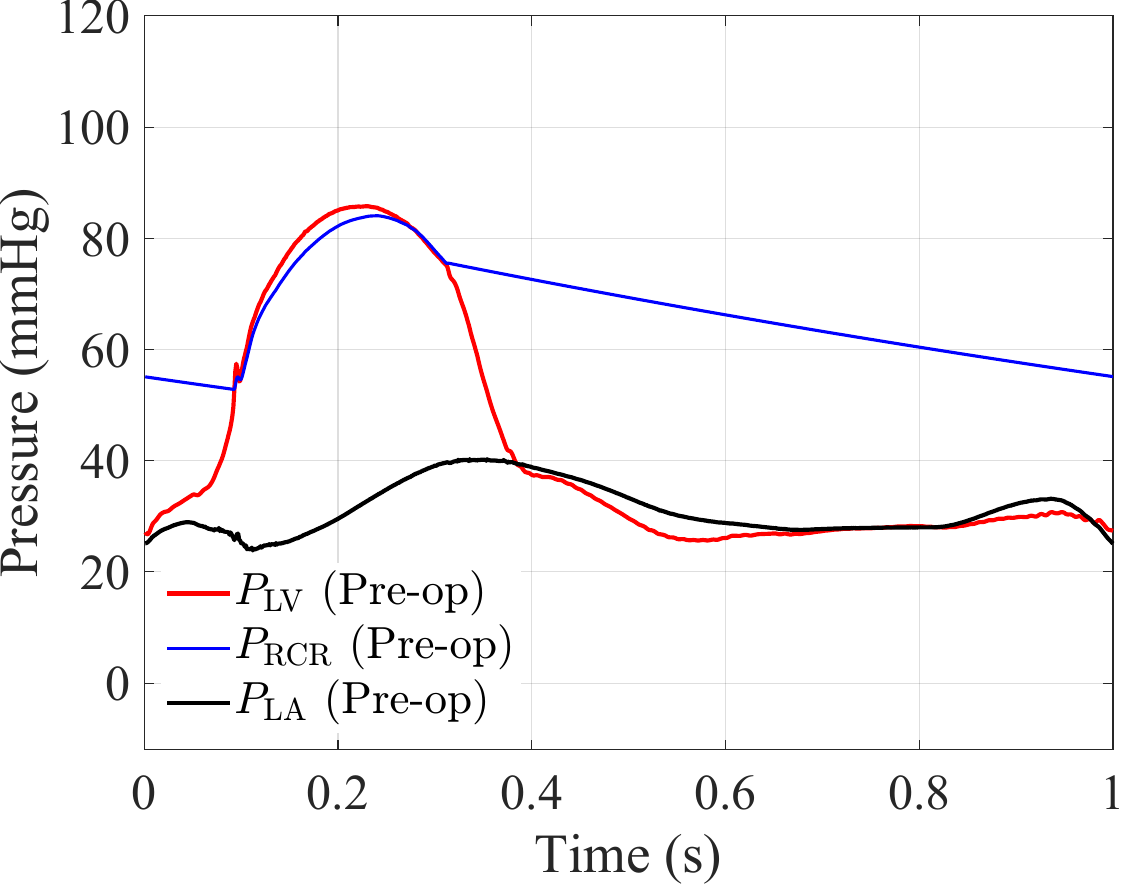}
        \caption{}
        \label{fig:prolapse_pressure}
	\end{subfigure}
    \hspace{0.01\textwidth}
    \begin{subfigure}[b]{0.48\textwidth}\centering
        \includegraphics[width=\textwidth]{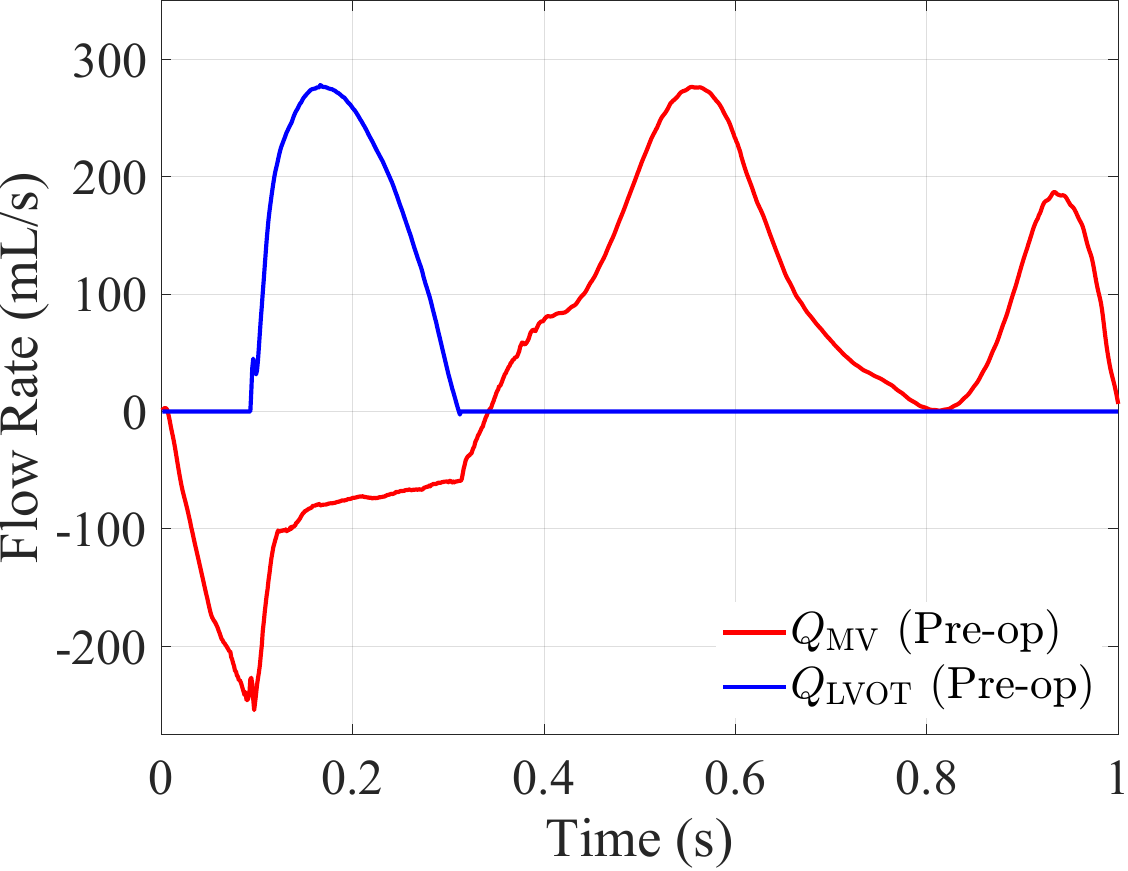}
        \caption{}
        \label{fig:prolapse_flowrate}
	\end{subfigure}
    \caption{(a) Pre-operative LV, LA, and RCR pressures over the cardiac cycle. (b) Pre-operative flowrate through the MV and LVOT over the cardiac cycle.}
    \label{fig:prolapse_graph}
\end{figure}

The pre-operative case represents a moderate-to-severe MR configuration characterized by a prolapsing posterior leaflet, serving as the pathological baseline against which all subsequent post-operative configurations are compared.
Figure~\ref{fig:prolapse_graph} illustrates the pre-operative hemodynamic quantities, and the pre-operative column of Table~\ref{t:PrePostOp_Summary} summarizes the hemodynamic and structural metrics at each phase over the cardiac cycle.

\begin{table}[!t]
\centering\small
\caption{Summary of pre-operative (prolapse) and post-operative hemodynamic and structural metrics over the cardiac cycle.}
\label{t:PrePostOp_Summary}
\begin{tabular}{@{}lccc@{}}
\toprule
\textbf{Metric} & \textbf{Pre-operative} & \textbf{Single-clip} & \textbf{Double-clip} \\
\midrule
\multicolumn{4}{@{}l}{\textit{Systole}} \\
Regurgitant Volume                                       & $32.9$~mL & $19.8$~mL & $16.3$~mL \\
LVOT Ejection Volume                                           & $40.5$~mL & $53.6$~mL & $57.1$~mL \\
Regurgitant Fraction                 & $44.8$\% & $27.0$\% & $22.2$\% \\
\midrule
\multicolumn{4}{@{}l}{\textit{Peak Mitral Regurgitation}} \\
Cycle Time                  & $0.15$~s & $0.13$~s & $0.13$~s \\
EROA                        & $0.46$~cm$^2$ & $0.22$~cm$^2$ & $0.18$~cm$^2$ \\
Regurgitant Jet Velocity$^a$    & $111$~cm/s & $123$~cm/s & $76$~cm/s \\
Average Jet Velocity$^a$        & $71$~cm/s & $65$~cm/s & $56$~cm/s \\
Transvalvular Pressure Gradient      & $54$~mmHg & $68$~mmHg & $72$~mmHg \\
Coaptation Area             & N/A$^b$ & $1.17$~cm$^2$ & $1.68$~cm$^2$ \\

\midrule
\multicolumn{4}{@{}l}{\textit{Diastole}} \\
EOA                         & $5.20$~cm$^2$ & $2.17$~cm$^2$ & $1.18$~cm$^2$ \\
Peak Transvalvular Velocity & $119$~cm/s & $220$~cm/s & $330$~cm/s \\
Mean Transvalvular Pressure Gradient      & $1$~mmHg & $5$~mmHg & $16$~mmHg \\
\midrule
\multicolumn{4}{@{}l}{\textit{Boundary Conditions}} \\
Outlet BC Switch (No-slip$\to$RCR)  & $0.090$~s & $0.064$~s & $0.064$~s \\
Outlet BC Switch (RCR$\to$No-slip)  & $0.32$~s & $0.34$~s & $0.34$~s \\
\bottomrule
\end{tabular}
\vspace{4pt}

{\footnotesize
$^a$See Remark~\ref{regurgitant_velocity} for details. 
$^b$Since prolapsed leaflets fail to achieve coaptation, coaptation area is not applicable for the pre-operative case. 
}
\end{table}

\begin{figure}[!t]\centering
  \centering
    \includegraphics[width=0.95\textwidth]{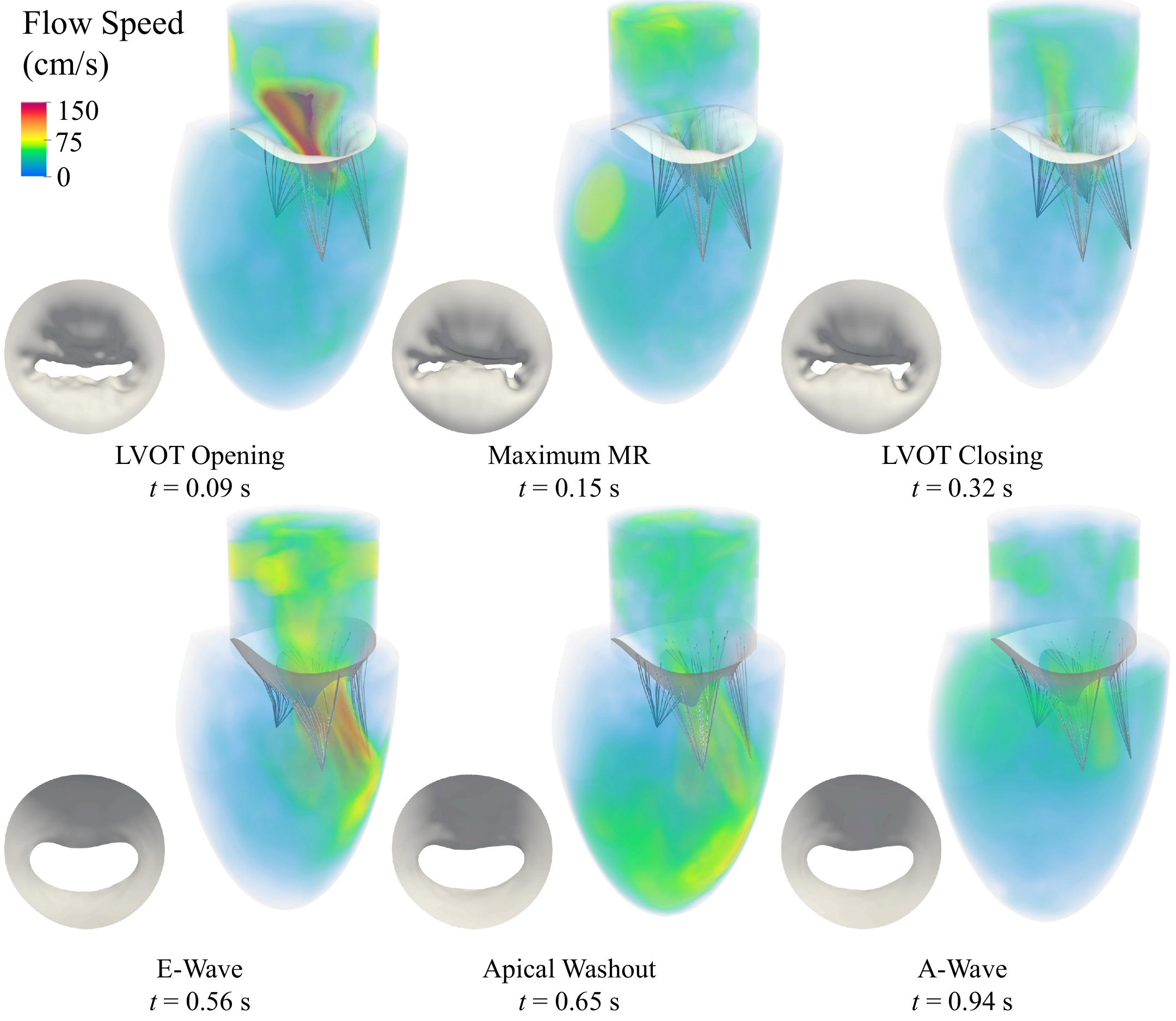}
  \caption{Visualization of the FSI results for the flow field in the left heart and MV motion over the complete cardiac cycle for mitral valve prolapse.}
  \label{f:Prolapse_Fluid}
\end{figure}

Figure~\ref{f:Prolapse_Fluid} visualizes the resulting flow field and valve kinematics throughout the cardiac cycle.
Ventricular contraction initiates at $t=0.0$~s, driving valve closure through increased left ventricular pressure before systolic ejection begins at $t=0.09$~s.
During early systole, we observe an initial closing volume, i.e., backflow as the valve closes, which is typically negligible in a healthy patient. 
However, the modeled valve and ventricle are dilated, which delays valve closure and results in a substantial closing volume. 
During the systolic ejection phase, flow leaves the ventricle through both the LVOT and the prolapsed valve due to incomplete closure of the MV.
Specifically, only $40.5$~mL is ejected forward while $32.9$~mL is lost to regurgitation, yielding a regurgitant fraction of $44.8\%$, consistent with moderate-to-severe MR.
Note that these values are calculated including the initial closing volume.

Effective regurgitant orifice area (EROA) and effective orifice area (EOA) are calculated over the systolic ejection and diastolic filling phases, respectively, using the formula and discharge coefficient given in Refs.~\cite{caballero2020comprehensive, chandran2012}, evaluated with root mean square flowrates and mean pressure gradients (see Remark~\ref{orifice_area}).
Transvalvular pressure gradients are measured at two planes normal to the valve axis, one positioned $1$~mm above the highest point of the annulus and the other $1$~mm below the lowest point of the free edge.
On each plane, one probe point is placed in the path of each orifice jet, giving as many probes per plane as there are orifices.
For instantaneous measurements, the pressures at the probes on each plane are averaged at a given time step.
For EROA and EOA calculations, the transvalvular pressure gradients are further averaged over the systolic ejection and diastolic filling phases, respectively, to approximate catheter-based pressure measurements.
We do not include the initial closing volume in the EROA calculation to isolate the leakage through the regurgitant orifices, which is the quantity TEER is intended to reduce, rather than excess backflow caused by a dilated valve and ventricle.

The regurgitant flow peaks at $t=0.15$~s, with a jet velocity of $111$~cm/s and a transvalvular pressure gradient of $54$~mmHg.
Over the systolic phase, the EROA is calculated to be $0.46$~cm$^2$, which is considered severe.
We also observe a low peak LV pressure of $86$~mmHg at around $t=0.23$~s, likely as a result of the moderate-to-severe MR, which delays and decreases systolic pressure build-up.
During the subsequent diastolic filling phase, the E-wave peaks at $t=0.56$~s, with a peak velocity of $119$~cm/s and an instantaneous transvalvular pressure gradient of $4$~mmHg. 
The E-wave pushes enough fluid toward the apex to produce apical washout, which can be seen at $t=0.65$~s. 
The A-wave follows shortly after at $t=0.94$~s, completing ventricular filling and concluding the cardiac cycle.
Over the diastolic phase, the EOA is calculated to be $5.20$~cm$^2$, which is far above the mitral stenosis threshold of $1.5$~cm$^2$~\cite{stone2015}.
This pre-operative case serves as the pathological baseline for the post-operative comparisons that follow.

\begin{remark}
\label{orifice_area}
The formulas used here to compute EROA and EOA follow the conventional Gorlin-type single-orifice model~\cite{gorlin1951hydraulic, chandran2012}, whose empirical discharge coefficient and assumed downstream pressure recovery are not strictly applicable to the multiple orifices produced by TEER.
We nonetheless report EROA and EOA because clinical guideline thresholds are defined in terms of these quantities~\cite{zoghbi2017, stone2015}, and evaluating them as in prior computational TEER studies~\cite{caballero2020comprehensive} keeps our values comparable to the literature.
\end{remark}

\begin{figure}[!t]\centering
	\begin{subfigure}[b]{0.48\textwidth}\centering
        \includegraphics[width=\textwidth]{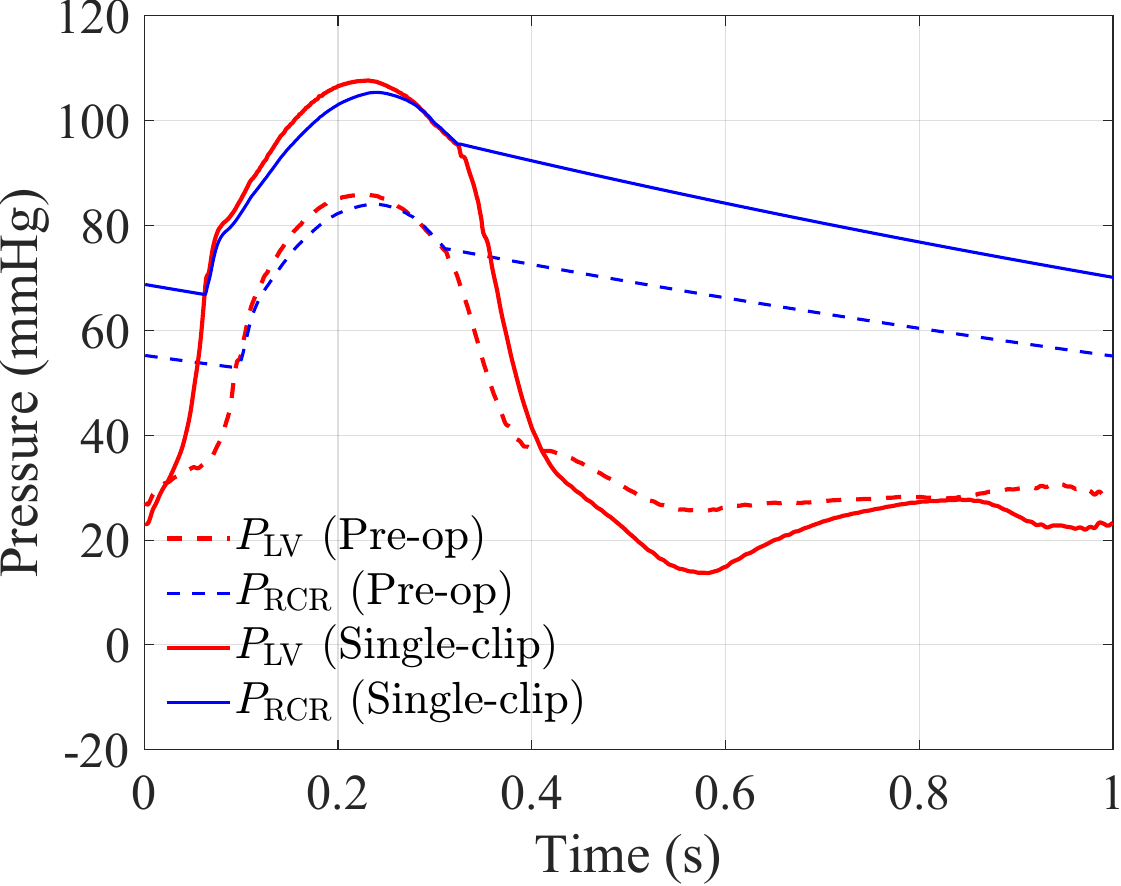}
        \caption{}
        \label{fig:1clip_pressure}
	\end{subfigure}
    \hspace{0.01\textwidth}
    \begin{subfigure}[b]{0.48\textwidth}\centering
        \includegraphics[width=\textwidth]{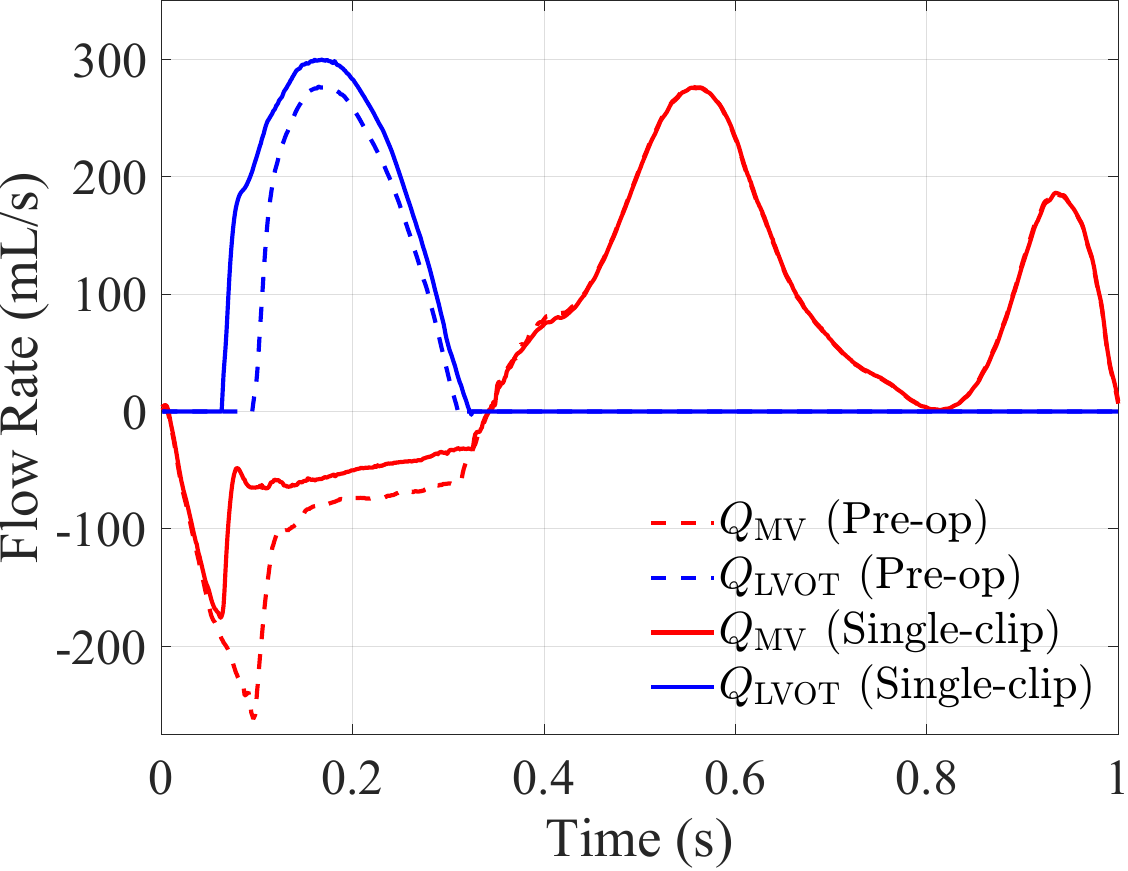}
        \caption{}
        \label{fig:1clip_flowrate}
	\end{subfigure}
    \caption{(a) Pressure in the left ventricle and RCR pressure over the cardiac cycle with a single clip. (b) Flowrate through the MV and LVOT over the cardiac cycle with a single clip. Pre-operative pressure and flowrate are shown as dashed lines for comparison.}
    \label{fig:1clip_graph}
\end{figure}

\begin{figure}[!t]\centering
	\begin{subfigure}[b]{0.48\textwidth}\centering
        \includegraphics[width=\textwidth]{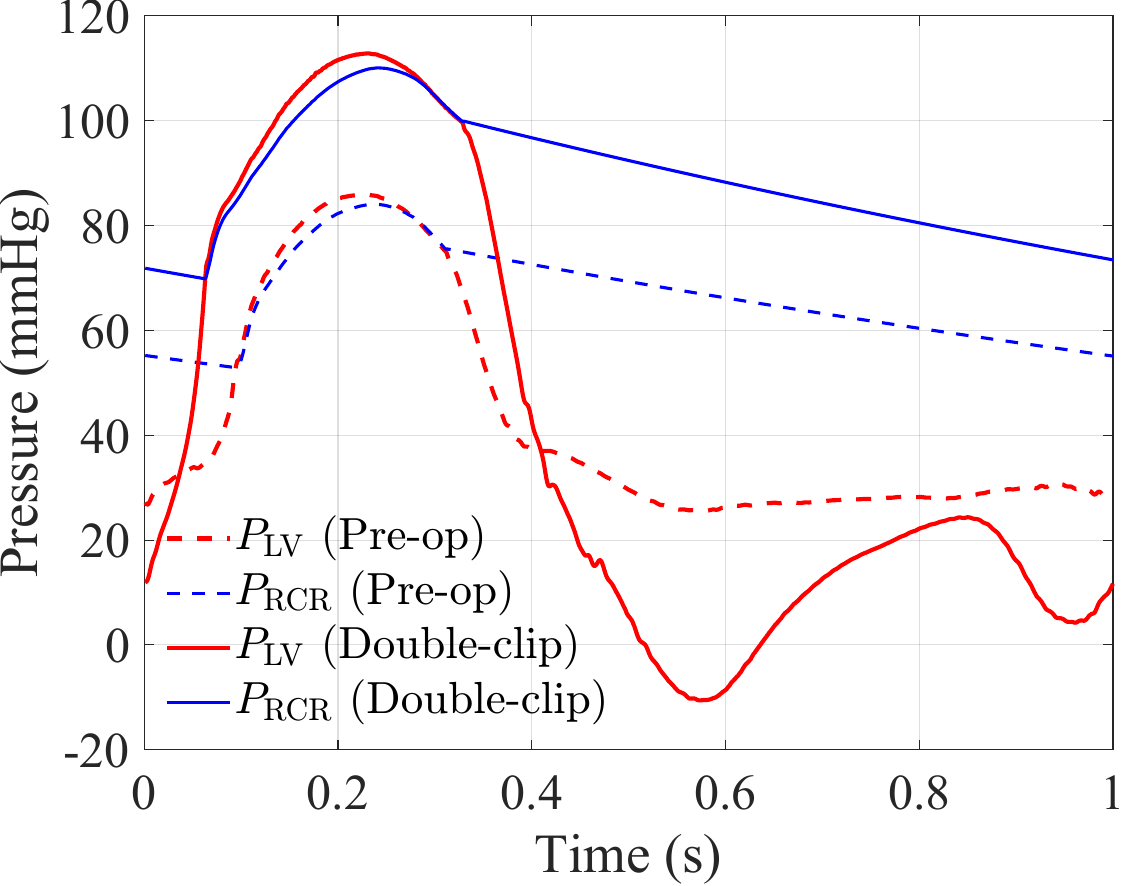}
        \caption{}
        \label{fig:2clip_pressure}
	\end{subfigure}
    \hspace{0.01\textwidth}
    \begin{subfigure}[b]{0.48\textwidth}\centering
        \includegraphics[width=\textwidth]{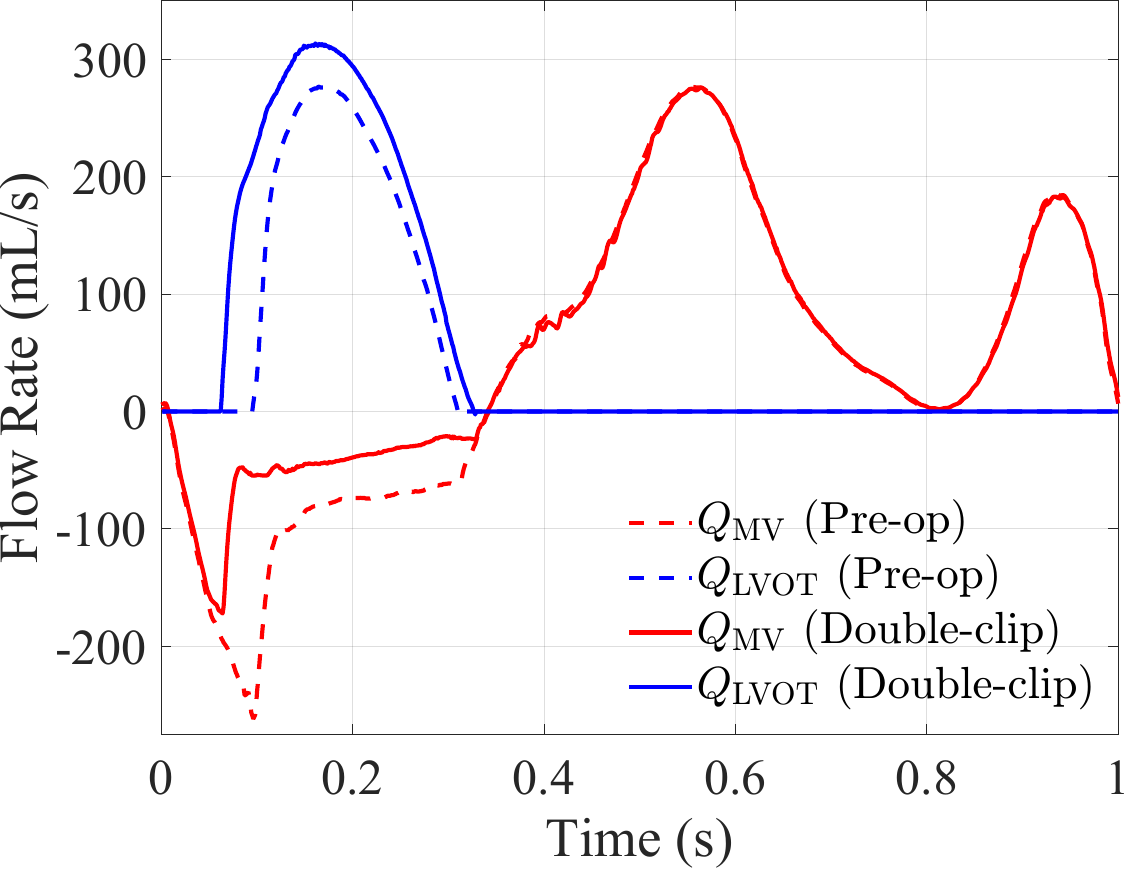}
        \caption{}
        \label{fig:2clip_flowrate}
	\end{subfigure}
    \caption{(a) Pressure in the left ventricle and RCR pressure over the cardiac cycle with double clips. (b) Flowrate through the MV and LVOT over the cardiac cycle with double clips. Pre-operative pressure and flowrate are shown as dashed lines for comparison.}
    \label{fig:2clip_graph}
\end{figure}

\subsubsection{Post-operative FSI}
We next consider the two post-operative cases: the single-clip repair and the double-clip repair.
This selection allows us to compare the hemodynamic efficacy of different TEER intervention strategies and assess how each affects both regurgitant flow and diastolic ventricular filling.

\begin{figure}[!t]\centering
  \centering
    \includegraphics[width=0.95\textwidth]{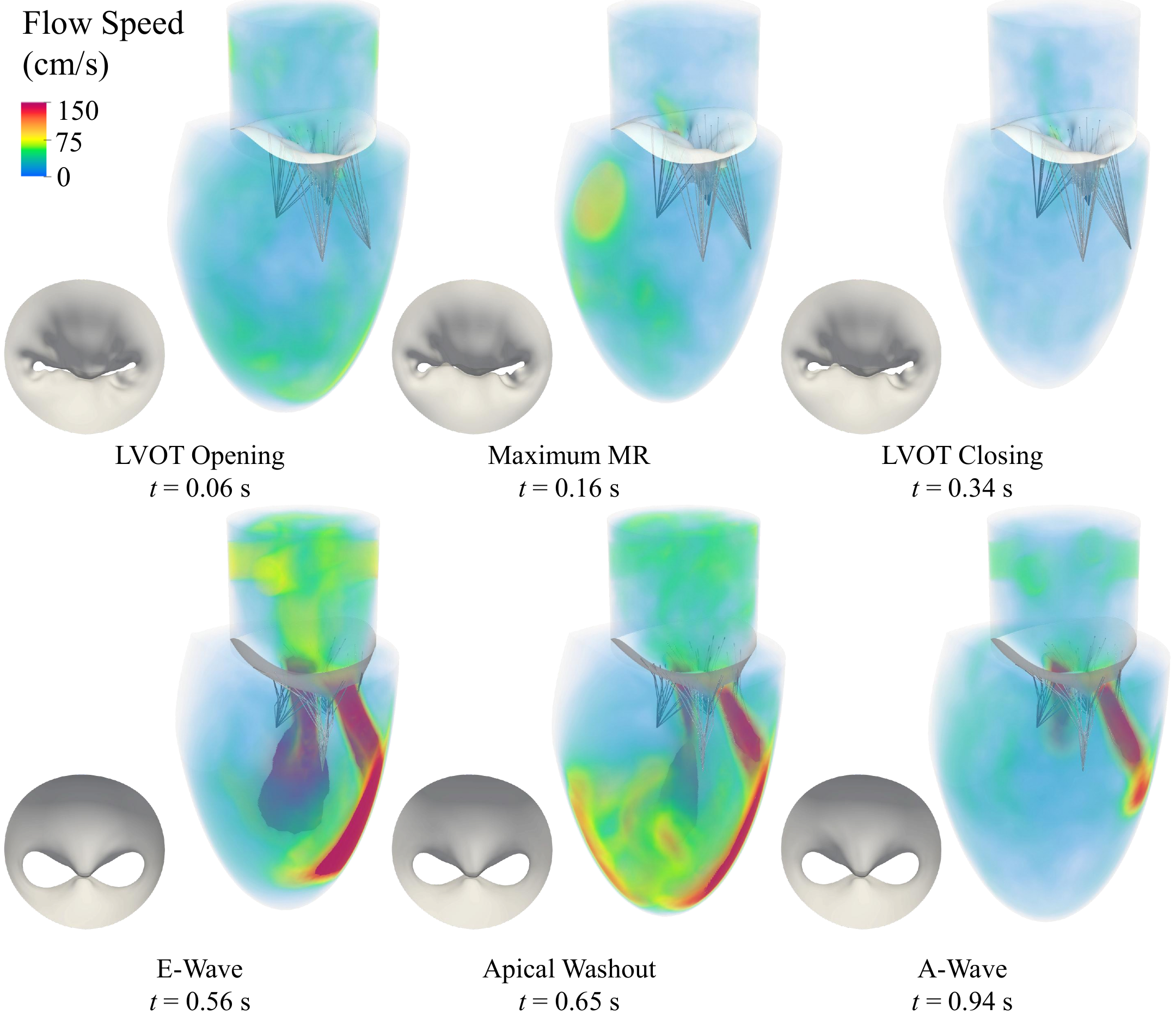}
  \caption{Visualization of the FSI results for the flow field in the left heart and MV motion over the complete cardiac cycle for single-clip TEER.}
  \label{f:1Clip_Fluid}
\end{figure}

\begin{figure}[!t]\centering
  \centering
    \includegraphics[width=0.95\textwidth]{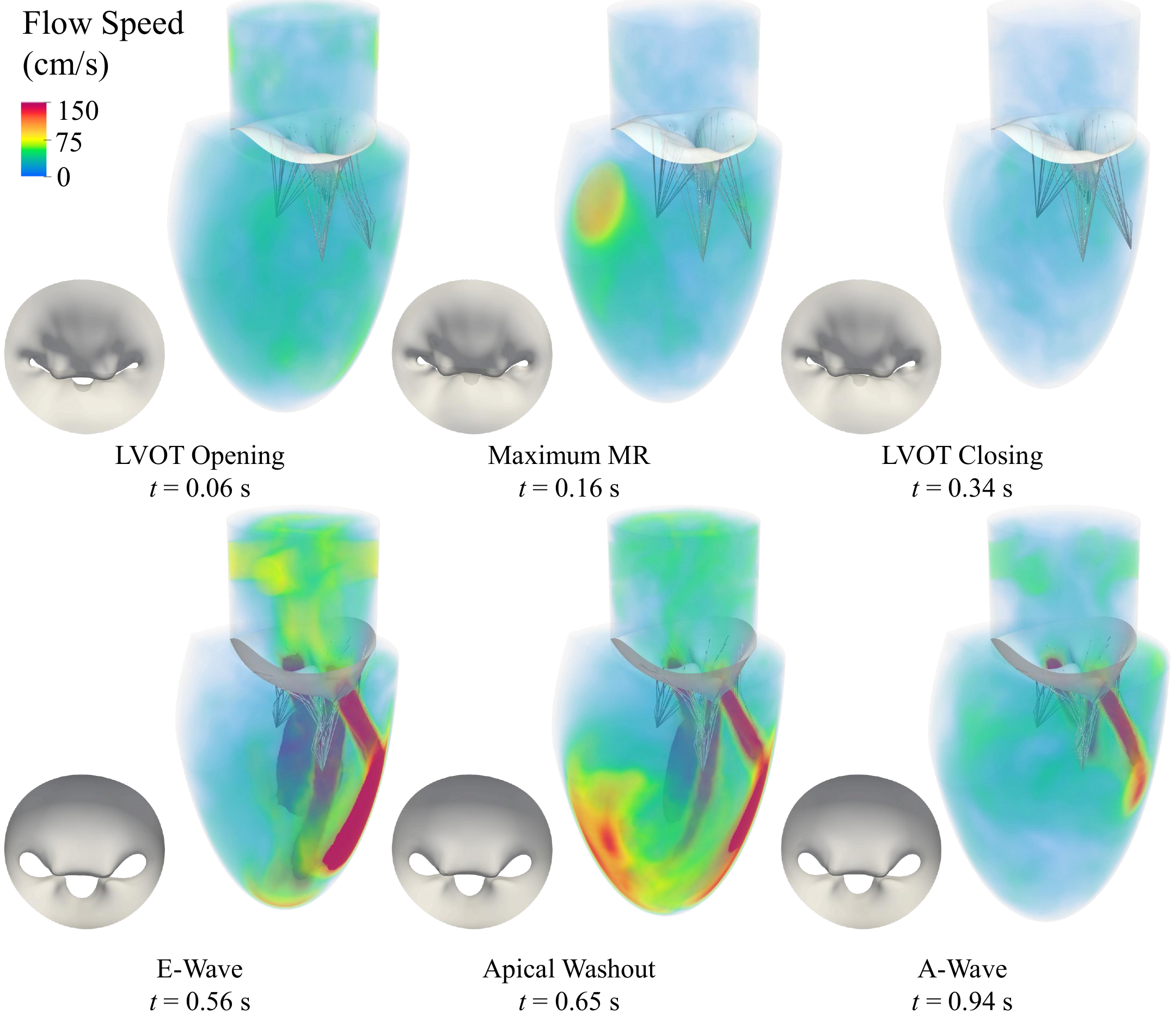}
  \caption{Visualization of the FSI results for the flow field in the left heart and MV motion over the complete cardiac cycle for double-clip TEER.}
  \label{f:2Clip_Fluid}
\end{figure}

\begin{figure}[!t]
  \centering
    \includegraphics[width=0.9\textwidth]{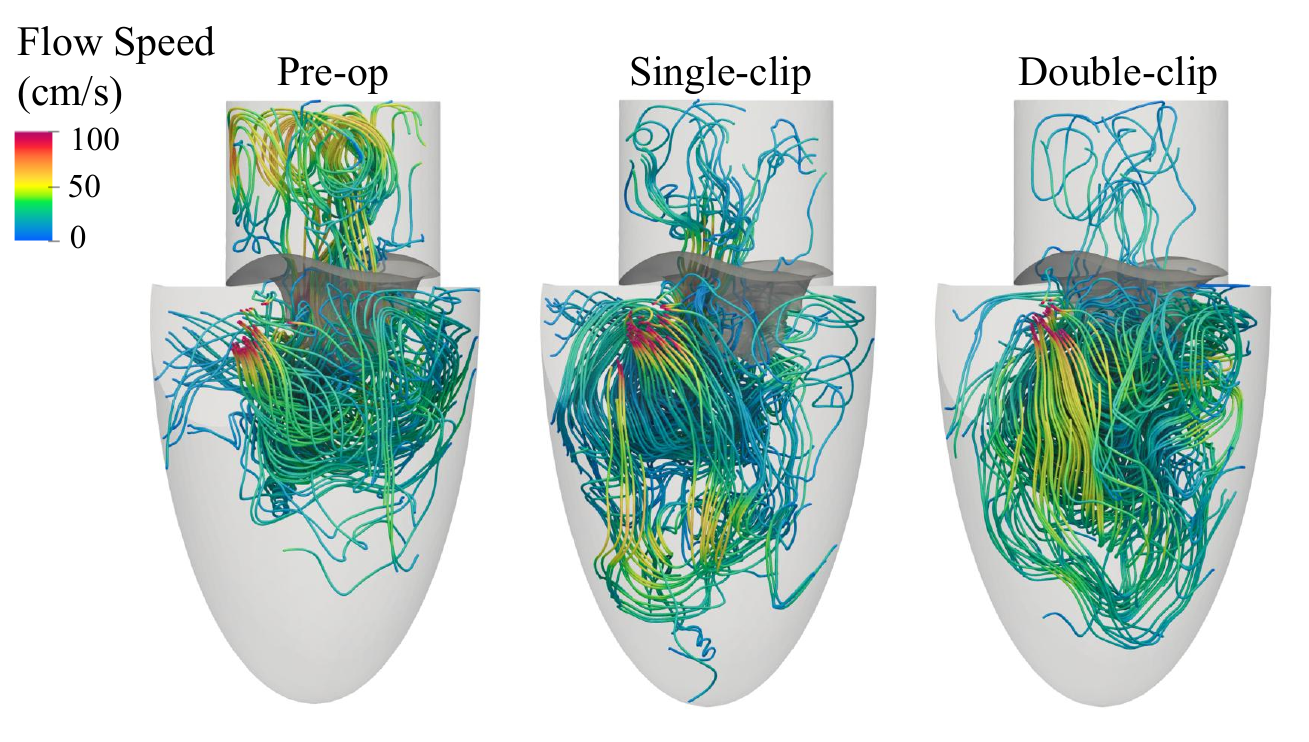}
  \caption{Visualization of streamlines at peak MR showing the jet behavior both leaving the LVOT and leaking through the open orifice of the mitral valve.}
  \label{f:streamlinesMR}
\end{figure}

Figures~\ref{fig:1clip_graph} and \ref{fig:2clip_graph} illustrate the temporal evolution of the hemodynamic quantities over the cardiac cycle for each respective case, and Figures~\ref{f:1Clip_Fluid} and \ref{f:2Clip_Fluid} show the corresponding flow fields and valve kinematics at representative phases throughout the cardiac cycle.
The corresponding hemodynamic and structural metrics for both post-operative configurations at each phase over the cardiac cycle are summarized in Table~\ref{t:PrePostOp_Summary}.
Compared to the prolapsed baseline, both repair strategies successfully improve leaflet coaptation and reduce the regurgitant fraction from $44.8\%$ to $27.0\%$ and $22.2\%$, respectively.
This demonstrates that TEER effectively mitigates systolic regurgitation. 
The substantial reduction in regurgitation is qualitatively visualized in Figure~\ref{f:streamlinesMR}.
The streamlines are visualized at the instant of peak MR to capture both the forward ejection through the LVOT and the leakage through the mitral valve, which confirms that the pre-operative case loses a substantially higher percentage of flow to regurgitation than the clipped cases. 

When measuring mitral regurgitant velocity, we report two values: regurgitant jet velocity and average jet velocity.
Functionally, when clinicians report regurgitant jet velocity, they do so through continuous-wave Doppler and report the highest velocity along a selected beam to measure the value.
This method is less effective post-TEER because the flow moves through multiple orifices and can be distributed unevenly among them, as seen in the Maximum MR panel of Figure~\ref{f:1Clip_Fluid}. 
This bias can heavily skew the measured regurgitant jet velocity.
Since our simulations provide the full three-dimensional velocity field at each time step, we are not limited to the measurements available to Doppler.
We therefore compute an average jet velocity at the time step associated with peak MR.
To calculate this value, we extract several slices of the fluid domain around the valve, compute the averaged velocity over the combined orifice areas for each slice, and report the highest value for each case.
Additional discussion of these techniques can be found in Remark~\ref{regurgitant_velocity}.

Leaflet coaptation area during closure further corroborates post-operative improvements in cardiac function. 
The prolapsed baseline exhibits no contact, whereas the single- and double-clip cases achieve coaptation of $1.17$~cm$^2$ and $1.68$~cm$^2$, respectively. 
This restored coaptation enables LV pressure to recover toward the physiological systolic peak of $120$~mmHg.
Furthermore, the EROA is substantially reduced for both cases, from $0.46$~cm$^2$ to $0.22$~cm$^2$ and $0.18$~cm$^2$, respectively.
Consequently, across all primary evaluated metrics for the overall systolic phase, both repair configurations successfully reduce the MR severity. 

During diastolic filling, as a consequence of TEER, the EOA is reduced from $5.20$~cm$^2$ to $2.17$~cm$^2$ and $1.18$~cm$^2$ for the single- and double-clip cases, respectively.
While the EOA remains above the clinically significant stenosis threshold of $1.5$~cm$^2$ for the single-clip case, the double-clip case falls below that threshold, which is a sign of moderate-to-severe stenosis.
The reduced orifice area also increases resistance to diastolic filling under the prescribed ventricular motion, producing a stronger suction effect and higher transvalvular flow velocities.
The diastolic transvalvular pressure gradients provide further evidence of stenosis. 
In both post-operative cases, the mean diastolic pressure gradient meets or exceeds the $5$~mmHg clinical pressure gradient threshold, reaching $5$~mmHg for the single-clip case and $16$~mmHg for the double-clip case.
These elevated gradients coincide with high-velocity jets through the reduced orifices, visible in the E- and A-wave panels of Figures~\ref{f:1Clip_Fluid} and \ref{f:2Clip_Fluid}.
Additionally, in the double-clip case, this restricted flow drives a negative LV pressure during diastolic filling, which further indicates impaired transmitral flow.

Note that our post-operative simulations reflect the immediate uncompensated hemodynamic response, showing the physiological behavior following MV repair in vivo~\cite{patel2014, gaasch2008}. The pre-operative case (Figure~\ref{fig:prolapse_pressure}) exhibits a drop in both LV and aortic pressure, reflecting hemodynamic response to acute MR reported in~\citet{patel2014}. Following MV repair, single- and double-clip cases show improved mean arterial pressure as shown in Figures~\ref{fig:1clip_pressure}~and~\ref{fig:2clip_pressure}, in agreement with the same study. 
However, the cardiovascular system also undergoes short- and long-term compensatory effects to preserve cardiac function, including tachycardic response, baroreflex, and ventricular and atrial remodeling~\cite{farandzha2023, yoran1983, shinbane1997, shoureshi2024}, which are not considered here.
Moreover, the negative LV pressure observed in the double-clip case during the E-wave suggests that the ventricle would need additional adaptation compared to the single-clip case to restore compensated behavior. 

\begin{figure}[!t]
  \centering
    \includegraphics[width=0.9\textwidth]{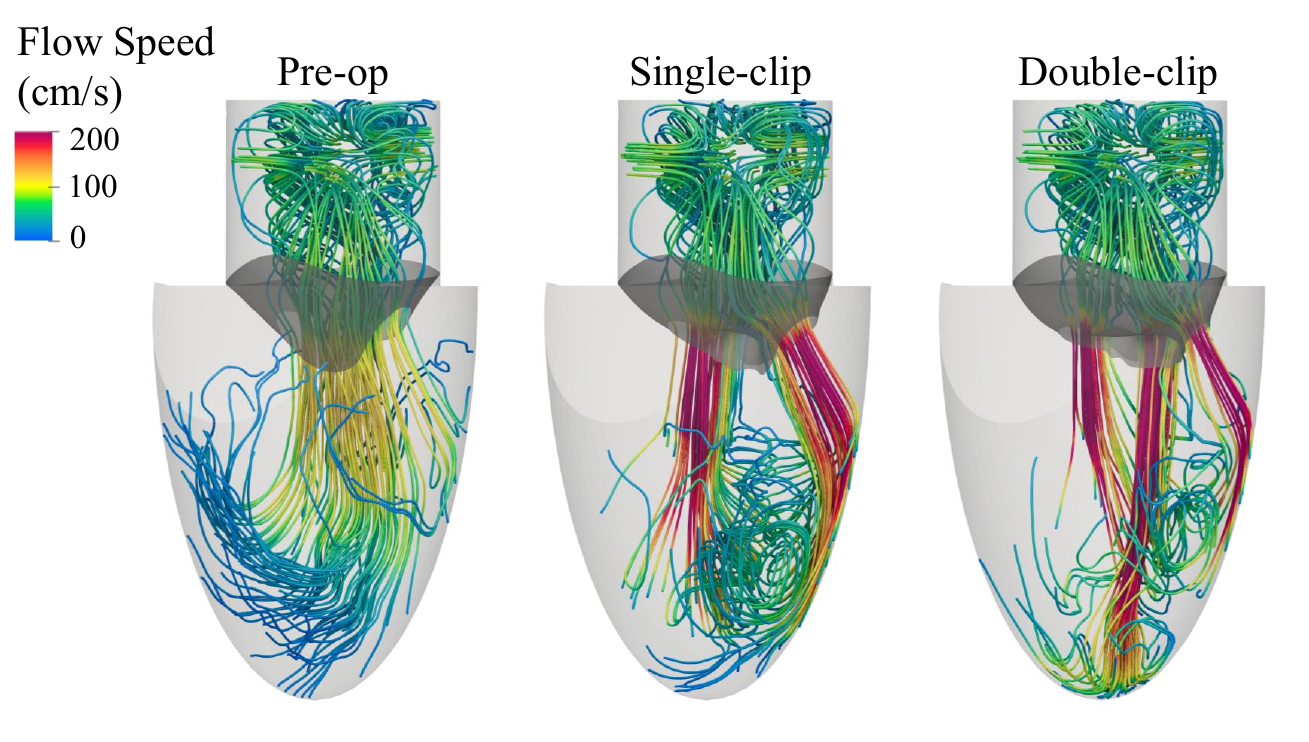}
  \caption{Visualization of streamlines at peak diastole showing the jet behavior leaving the valve as it travels to the ventricular apex.}
  \label{f:streamlines}
\end{figure}

Figure~\ref{f:streamlines} highlights the flow patterns through the mitral valve during peak diastole using streamlines. 
Instead of a single inflow jet seen in the pre-operative case, the post-operative cases generate multiple higher-velocity jets. 
In the single-clip case, no transvalvular flow passes through the central region of the mitral orifice because of the centrally positioned clip, and the flow is instead diverted toward the lateral ventricular walls. 
These diverted jets impinge on the ventricular walls before traveling toward the apex and creating a region of highly perturbed flow, which potentially reduces apical washout. 
Reduced apical washout is clinically relevant because it is associated with thrombus formation~\cite{osti_22403199, Delemarre1990, delewi2012}.
In the double-clip case, on the other hand, this effect is less pronounced because the central region is not obstructed by the clip devices and because the jets from the two lateral orifices are thinner and impinge on the ventricular walls at more basal locations.
Consequently, the region of perturbed flow does not extend as far toward the ventricular apex, so apical washout is potentially preserved.

\begin{remark}
\label{regurgitant_velocity}
Quantification methods used in pre-operative assessment, such as regurgitant jet velocity and 2D proximal isovelocity surface area (PISA), are difficult to apply post-operatively because of the presence of multiple orifices and jets.
Recent studies instead suggest using a time-averaged velocity over systole, area-based methods such as those that calculate EROA, or 3D PISA~\cite{marchetti2023, caballero2022, rottlander2025a}. 
The spatially averaged jet velocity metric used here is inspired by these approaches and leverages the computational nature of our simulations.
The difference in reported values between the two measurement methods can be seen in Table~\ref{t:PrePostOp_Summary}, where the regurgitant jet velocity for the single-clip case is higher than that of the pre-operative case, even though the single-clip repair reduces MR. 
The double-clip case also has a peak jet velocity $47$~cm/s slower than that of the single-clip case, which is much lower than expected given the relative regurgitant fraction between the single- and double-clip cases.
When using the average velocity to quantify regurgitant jets, we obtain a pre-operative case velocity of $71$~cm/s, a single-clip case velocity of $65$~cm/s, and a double-clip case velocity of $56$~cm/s.
Analyzing these values shows that the single- and double-clip velocity values are far closer to each other than before and more accurately reflect what is seen in Figures~\ref{f:1Clip_Fluid}~and~\ref{f:2Clip_Fluid}. 
This analysis gives us a more robust understanding of the mitral regurgitant jet behavior without overreliance on a single peak-velocity measurement.
\end{remark}

\section{Conclusions}
\label{s:discussion}
This paper develops an immersogeometric FSI model of the native human MV and its transcatheter edge-to-edge repair.
To obtain analysis-ready geometry from clinical images, we introduce a fully parameterized reconstruction pipeline that converts segmented 3D TEE data of a patient-specific MV into B-spline leaflet surfaces, using a PDE-based parameterization and a single regularized linear least-squares solve.
Functionally equivalent chordae are then generated on the reconstructed leaflets from population-averaged insertion and origin maps.
Because the leaflets of a native valve are under physiological load at the time of imaging, the imaged geometry is not stress-free. 
We therefore reformulate the Kirchhoff--Love shell to include an in-vivo prestrain tensor.
This allows the imaged geometry to serve directly as the reference configuration while the constitutive law and its consistent linearization remain compatible with the plane-stress and incompressibility reduction of the shell. 
To model the repair, we propose a patch-based penalty formulation for the clip that enforces leaflet coaptation over a parametrically defined region, which allows the number and placement of clips to be changed without modeling the device or modifying the leaflet discretization.

To simulate the complete MV apparatus, we formulate the IMGA FSI problem with weakly enforced fluid--shell kinematic constraints, a shell--cable penalty coupling at the chordal insertions, and a nonlocal contact formulation among all structural components.
Within this formulation, we introduce a subcell quadrature strategy that generates fluid--structure coupling points independently of the structural quadrature, which removes the over-refinement of the shell mesh otherwise required to prevent leakage as the background fluid mesh is refined.
To assess the robustness of the proposed framework, we evaluate structural, TEER penalty, fluid, and subcell convergence.
The mesh refinement study shows that converged valvular dynamics and hemodynamic quantities can be obtained with the structural and fluid mesh resolutions selected independently.
The framework is then applied to a prolapsed MV and its single- and double-clip repairs, with the same left heart model, boundary conditions, and calibrated lumped-parameter model used for all cases.
The resulting changes in valve function and left-heart hemodynamics are qualitatively consistent with the acute response to TEER reported clinically~\cite{siegel2011acute, herrmann2006mitral, gaemperli2013real}.

Several limitations remain before the framework can support patient-specific treatment planning.
The constitutive models for the leaflets and CTs use population-level parameters, and more advanced models with patient-specific material properties have yet to be incorporated.
The left heart is idealized with prescribed wall motion, and predicting its hemodynamics requires deformable walls with active contraction and a closed-loop lumped-parameter model, both fitted to the patient.
With these extensions, the framework presented here provides a path toward predictive simulation of TEER in individual patients.

\section*{Acknowledgments}
This work was supported in part by the National Heart, Lung, and Blood Institute of the National Institutes of Health under award number R01HL184128 and by the National Science Foundation under award number DMS-2436623. 
K.H. Kim acknowledges research support through the American Heart Association Postdoctoral Fellowship 25POST1366103.
This support is gratefully acknowledged. 
We also thank the Texas Advanced Computing Center (TACC) at The University of Texas at Austin for providing computational resources that have contributed to the research results reported within this paper.

\small
\bibliographystyle{unsrtnat}

\bibliography{ref}

\end{document}